\pdfoutput=1
\RequirePackage{silence}
\documentclass[a4paper,svgnames]{amsart}
\usepackage[hmarginratio={1:1},vmarginratio={1:1},lmargin=50.0pt,tmargin=45.0pt]{geometry}

\usepackage{booktabs}
\allowdisplaybreaks

\usepackage[T1]{fontenc}
\usepackage{latexsym,exscale,amsfonts,amssymb,mathtools}
\usepackage{amsmath,amsthm,amsfonts,amssymb,amscd,textcomp,bbm}
\usepackage{mathrsfs,stackrel}
\usepackage{etoolbox}
\usepackage{tikz,tikz-cd}
\usetikzlibrary{matrix,arrows.meta,decorations.markings,shapes.multipart,decorations.pathreplacing,decorations.pathmorphing,shapes.geometric}
\usepackage{aliascnt}
\usepackage{ytableau}
\usepackage{xparse}
\usepackage{cite}

\usepackage{dynkin-diagrams}
\usepackage{nicematrix}
\usepackage{tikz-3dplot}
\usepackage{mathrsfs}
\DeclareMathAlphabet{\mathscrbf}{OMS}{mdugm}{b}{n}

\usepackage{array}
\newcolumntype{C}{>{$}c<{$}}
\newcolumntype{P}[1]{>{\centering\arraybackslash}p{#1}} 

\usepackage{xcolor,colortbl}

\definecolor{mygray}{gray}{0.6}
\definecolor{mygraydark}{gray}{0.4}
\definecolor{mygraylight}{gray}{0.85}
\definecolor{spinach}{RGB}{46,139,87}
\definecolor{tomato}{RGB}{255,99,71}
\definecolor{orchid}{RGB}{143,40,194}
\definecolor{neon}{RGB}{77,77,255}
\definecolor{lightneon}{RGB}{110,110,255}
\definecolor{pumpkin}{RGB}{224,180,80}
\definecolor{citron}{RGB}{190,180,90}

\definecolor{lava}{RGB}{207,16,32}
\definecolor{cream}{RGB}{255,253,208}
\definecolor{verdigris}{RGB}{67,179,174}
\definecolor{Black}{RGB}{0,0,0}
\definecolor{mydarkblue}{RGB}{10,10,170}
\definecolor{darkspinach}{RGB}{20,70,20}
\definecolor{darktomato}{RGB}{155,40,30}
\definecolor{darkorchid}{RGB}{50,10,100}
\definecolor{darklava}{RGB}{150,8,16}

\definecolor{zero}{RGB}{0,0,0}
\definecolor{one}{RGB}{255,0,0}
\definecolor{two}{RGB}{0,255,0}
\definecolor{three}{RGB}{0,0,255}

\usepackage{todonotes}

\usepackage{enumitem}
\setlist[enumerate]{itemsep=0.15cm,label=\emph{\upshape(\alph*)}}
\setlist[enumerate,2]{itemsep=0.15cm,label=\emph{\upshape(\roman*)}}
\setlist[enumerate,3]{itemsep=0.15cm,label=\emph{\upshape(\Alph*)}}

\let\emph\relax
\DeclareTextFontCommand{\emph}{\bfseries\em}

\renewcommand{\dots}{\text{...}}

\newcommand{\placeholder}{{}_{-}}
\newcommand{\mystrut}{\rule[-0.2\baselineskip]{0pt}{0.9\baselineskip}}

\newcommand{\ie}{\text{i.e.}}
\newcommand{\eg}{\text{e.g.}}
\newcommand{\cf}{\text{cf.}}
\newcommand{\etc}{\text{etc.}}

\usepackage[all]{xy}
\usepackage{tikz}
\usetikzlibrary{cd}
\usetikzlibrary{decorations}
\usetikzlibrary{decorations.markings}
\usetikzlibrary{decorations.pathreplacing}
\usetikzlibrary{decorations.pathmorphing}
\usetikzlibrary{arrows.meta,shapes,positioning,matrix,calc}
\usetikzlibrary{shapes.callouts}
\tikzset{anchorbase/.style={baseline={([yshift=-0.5ex]current bounding box.center)}},
tinynodes/.style={font=\tiny,text height=0.25ex,text depth=0.05ex},
smallnodes/.style={font=\scriptsize,text height=0.75ex,text depth=0.15ex},
usual/.style={line width=2.0,color=black},
crossline/.style={preaction={draw=white,line width=5.75pt,-}},
}
\tikzstyle directed=[postaction={decorate,decoration={markings,mark=at position #1 with {\arrow[line width=0.3mm, black]{>}}}}]

\newcommand{\C}{\mathbb{C}}
\newcommand{\R}{\mathbb{R}}

\newcommand{\N}{\mathbb{Z}_{\geq 0}}
\newcommand{\Q}{\mathbb{Q}}
\newcommand{\Z}{\mathbb{Z}}

\newcommand{\op}[1]{\operatorname{#1}}
\newcommand{\optt}[1]{\operatorname{\texttt{#1}}}

\newcommand{\sank}{n}

\def\NewTheorem#1{%
\newaliascnt{#1}{equation}%
\newtheorem{#1}[#1]{#1}%
\aliascntresetthe{#1}%
\expandafter\def\csname #1autorefname\endcsname{#1}%
}
\def\equationautorefname~#1\null{(#1)\null}

\numberwithin{equation}{subsection}

\NewTheorem{Proposition}
\NewTheorem{Theorem}
\NewTheorem{Corollary}
\AtEndEnvironment{Corollary}{\null\hfill$\square$}%
\NewTheorem{Lemma}
\NewTheorem{Conjecture}
\NewTheorem{Speculation}
\NewTheorem{Observation}
\NewTheorem{Assumption}
\theoremstyle{definition}
\NewTheorem{Definition}
\AtEndEnvironment{Definition}{\null\hfill$\Diamond$}%
\NewTheorem{Classification Problem}
\AtEndEnvironment{Classification Problem}{\null\hfill$\Diamond$}%
\NewTheorem{Notation}
\AtEndEnvironment{Notation}{\null\hfill$\Diamond$}%
\NewTheorem{Example}
\AtEndEnvironment{Example}{\null\hfill$\Diamond$}%
\NewTheorem{Examples}
\AtEndEnvironment{Examples}{\vskip-10mm\null\hfill$\Diamond$}%

\theoremstyle{remark}
\NewTheorem{Remark}
\AtEndEnvironment{Remark}{\null\hfill$\Diamond$}%
\NewTheorem{Question}
\AtEndEnvironment{Question}{\null\hfill$\Diamond$}%

\usepackage[hypertexnames=false]{hyperref}
\usepackage{bookmark}
\hypersetup{
pdftoolbar=true,
pdfmenubar=true,
pdffitwindow=false,
pdfstartview={FitH},
pdftitle={Big data comparison of quantum invariants},
pdfauthor={Daniel Tubbenhauer and Victor Zhang},
pdfsubject={},
pdfcreator={Daniel Tubbenhauer and Victor Zhang},
pdfproducer={Daniel Tubbenhauer and Victor Zhang},
pdfkeywords={},
pdfnewwindow=true,
colorlinks=true,
linkcolor=mydarkblue,
citecolor=teal,
filecolor=magenta,
urlcolor=orchid,
linkbordercolor=lava,
citebordercolor=teal,
urlbordercolor=orchid,
linktocpage=true
}

\def\makeautorefname#1#2{\csdef{#1autorefname}{#2}}
\makeautorefname{section}{Section}%
\makeautorefname{subsection}{Section}%
\makeautorefname{subsubsection}{Section}%

\begin{document}
\title[Big data comparison of quantum invariants]{Big data comparison of quantum invariants}
\author[D. Tubbenhauer and V. Zhang]{Daniel Tubbenhauer and Victor Zhang}

\address{D.T.: The University of Sydney, School of Mathematics and Statistics F07, Office Carslaw 827, NSW 2006, Australia, \href{http://www.dtubbenhauer.com}{www.dtubbenhauer.com}, \href{https://orcid.org/0000-0001-7265-5047}{ORCID 0000-0001-7265-5047}}
\email{daniel.tubbenhauer@sydney.edu.au}

\address{V.Z.: University of New South Wales (UNSW), School of Mathematics and Statistics, NSW 2052, Australia, \href{https://dustbringer.github.io/}{dustbringer.github.io} }
\email{victor.zhang3@student.unsw.edu.au}

\begin{abstract}
We apply big data techniques, including exploratory and topological data analysis, to investigate quantum invariants. More precisely, our study explores the Jones polynomial's structural properties and contrasts its behavior under four principal methods of enhancement: coloring, rank increase, categorification, and leaving the realm of Lie algebras.
\end{abstract}


\addtocontents{toc}{\protect\setcounter{tocdepth}{1}}

\maketitle

\tableofcontents

\section{Introduction}\label{S:Intro}

In this paper, we define a \emph{quantum invariant} $(\mathfrak{g},V,\epsilon)$ as a knot invariant constructed, as, for example, in \cite{ReTu-invariants-3-manifolds-qgroups,We-knot-invariants}, by selecting a complex semisimple Lie algebra $\mathfrak{g}$, or a similar object, and a simple complex representation $V$ of $\mathfrak{g}$. We let $\epsilon=0$ denote the uncategorified version and $\epsilon=1$ represents its categorified counterpart. A key example is the Jones polynomial, given by $(\mathfrak{sl}_{2},\C^{2},0)$, where $\C^{2}$ is the vector representation of $\mathfrak{sl}_{2}$. Using a big data approach on the set of prime knots with up to $16$ crossings, we examine how the different methods of enhancing a quantum invariant, by varying $\mathfrak{g}$, $V$, or $\epsilon$, compare and interact.
Our main players are:
\begin{enumerate}
\item We start with the \emph{Jones polynomial} or \emph{A1 invariant} $(\mathfrak{sl}_{2},\C^{2},0)$ (for the vector representation). This is our reference invariant.

\item We investigate the \emph{2-colored Jones polynomial} or \emph{B1 invariant} $(\mathfrak{sl}_{2},\mathrm{Sym}^{2}\C^{2},0)$ (for the simple three-dimensional representation). This is coloring.

\item We look at the \emph{A2 invariant} $(\mathfrak{sl}_{3},\C^{3},0)$ (for the vector representation). This is a rank increase.

\item We then look at \emph{Khovanov homology} or \emph{A1${}^{c}$ invariant} $(\mathfrak{sl}_{2},\C^{2},1)$ (for the vector representation). This is categorification.

\item Finally, we have the most classical knot polynomial, the \emph{Alexander polynomial} or \emph{isotropic A1 invariant}
$(\mathfrak{gl}_{1|1},\C^{1|1},0)$ (for the vector representation). Here we leave the realm of Lie algebras.

\end{enumerate}
Some of our arguments are more general; see, for example, \autoref{T:Drop} and the text below it.

\begin{Remark}\label{R:Code}
All data files and additional material
(such as interactive and higher resolution pictures, code, and a possible empty erratum {\etc}) can be found \textit{online} at \cite{TuZh-quantum-big-data-code}. Included in this are various \textit{interactive} materials that the reader can run themselves.
\end{Remark}

\subsection{Background and ideas}

The Jones polynomial \cite{Jo-jones-polynomial} is widely recognized as one of the most significant knot invariants of the twentieth century, and Jones' groundbreaking discovery uncovered profound and unexpected connections between diverse areas of mathematics and physics. The Jones polynomial serves as the foundation of a larger family of knot invariants derived from quantum groups, known as (Witten)--Reshetikhin--Turaev invariants \cite{ReTu-invariants-3-manifolds-qgroups}. These invariants also have categorified counterparts \cite{Kh-cat-jones,We-knot-invariants}, positioning them as remarkable objects at the intersection of multiple mathematical fields.

One might ask a naive question:
\begin{gather*}
\fcolorbox{orchid!50}{spinach!10}{\mystrut``How effective are these invariants as tools for distinguishing knots?''} 
\end{gather*}
Such a focus underestimates the rich interplay of ideas and deep insights to which these invariants contribute. However, it is precisely this question that forms the focus of this paper. It turns out that this may be the wrong question to ask: \emph{quantum invariants are not expected to be strong knot invariants}. (This seems to be folk knowledge; see, for example, \cite{St-number-polynomials,110646} or \autoref{T:Drop}.) Indeed, they satisfy local relations ({\eg} skein relations), which make them easy to study but not particularly well suited as invariants.

With this in mind, a better question is:
\begin{gather*}
\fcolorbox{orchid!50}{spinach!10}{\mystrut``How do quantum invariants compare in distinguishing knots?''} 
\end{gather*}
This question, along with its variations, forms the central focus of this paper.

Our methodology employs systematic data processing and analysis, often categorized under \emph{``big data''}, using techniques such as data visualization, exploratory data analysis (EDA), and topological data analysis (TDA). This approach is inspired by previous work in knot theory and representation theory, such as \cite{LeHaSa-jones, DlGuSa-mapper, XYZ, LaTuVa-big-data}, and stands in contrast to deep learning-based methods in these fields.

Our computation took $\approx$ four months on the servers of the University of New South Wales (UNSW), heavily exploiting parallel computing.
The programs used were the Mathematica package KnotTheory \cite{KnotTheory} and a Javascript program written specifically for this project; see the link in \autoref{R:Code}.

Two final comments:
\begin{enumerate}[label=(\roman*)]

\item The computational complexity of (most) quantum invariants is not polynomial in the number of crossings, but can vary quite a bit; see \autoref{SS:Algorithm} for a more detailed analysis. We will mostly ignore this factor in the comparisons.

\item A critical aspect of this study is the scale of the data, {\eg} in small datasets, there's a risk of overfitting.
For this paper there is no issue as the number of knots grows exponentially in the number of crossings. For example, for $k_{n}$ = \#\{prime knots with $n$ crossings\}, \cite{ErSu-growth-knots,We-number-knots} proved that $2.68\leq\liminf_{n\to\infty}\sqrt[n]{k_{n}}\leq\limsup_{n\to\infty}\sqrt[n]{k_{n}}\leq 13.5$.

\end{enumerate}

\begin{Notation}
In this paper, following \cite{LaTuVa-big-data}, we have \emph{conjectures} and \emph{speculations}. 
Conjectures are presented in their standard sense, while speculations refer to preliminary hypotheses that lack full support from the data.
We hope that both serve as an inspiration to prove or, \fbox{{\color{purple}equally exciting}}, disprove the corresponding statements.
\end{Notation}

\subsection{Future directions}

Beyond the questions we ask below and extending our study \emph{from knots to links}, the following might be interesting.
\begin{enumerate}[label=(\Roman*)]

\item \textit{Changing the Lie type.} One thing not addressed in this paper is what happens when one changes the Lie type, for example, using the quantum B2, G2 or F4 invariants as nicely and diagrammatically constructed in \cite{Ku-spiders-rank-2,SaWe-quantum-f4}.

\item \textit{More colorings.} As we shall see, coloring seems to be a good way to enhance the Jones polynomial. In particular, it should be interesting to study all colored Jones polynomials at once (maybe using symmetric webs \cite{RoTu-symmetric-howe}), and would provide understanding into related invariants (hyperbolic volume, A-polynomial), see {\eg} \cite{Gu-cs-and-the-a-polynomial}.

\item \textit{Study the homology, not the polynomial.} We actually do not study Khovanov homology itself, but rather its Hilbert--Poincar{\'e} polynomial. This does not make any difference when working over $\Q$, but we ignore the torsion invariants. We do not expect them to change the picture, but we have no data to support that. See \cite{MuPrSiWaYa-torsion} for some directions, however, for very special knots.

\item \textit{Nonsemisimple invariants.}
Under the slogan of modified traces, see {\eg}
\cite{GePaMiTu-modified-qdims,GeKuPaMi-gen-traces-modified-dimensions}, non-semisimple knot invariants arise, which one still could call quantum invariants.
(One should be able to use quantum Satake and Nhedral Soergel bimodules, see \cite{El-q-satake,MaTu-soergel,MaMaMiTu-trihedral,LaTuVa-nhedral}, to produce and potentially enhance these invariants.) Putting these into the picture would be a nice addition.

\item \textit{More exotic quantum invariants.} 
All of the invariants studied in this paper come from tensor categories associated with generic quantum groups. Although we expect other invariants coming from more exotic tensor categories to behave similarly, it would be better to have supporting data. Examples include invariants derived from \cite{FlLa-knots,AnWe-mixed-qgroup,SuTuWeZh-mixed-tilting,MaMaMiTuZh-soergel-2reps}, either explicitly or implicitly via the exotic tensor categories that appear in those papers.

\item \textit{Three manifolds.} Some of the quantum invariants also give rise to three manifold invariants, and it would be interesting to know their effectiveness as invariants. However, three manifolds are difficult to enumerate, so the data size is limited, see \cite{MR3431029} for some progress.

\end{enumerate}
However, for none of these do we expect the overall behavior to change, so our main interests arise again from their comparison.
\medskip

\noindent\textbf{Acknowledgments.}
We extend our gratitude to the organizers of the exceptional conference ``Diagrammatic Intuition and Deep Learning in Mathematics'' (York, July 2024), and to Radmila Sazdanovic, whose inspiring talk at the event provided the spark for this project. We also thank Abel Lacabanne and Pedro Vaz for many related discussions, and William Hobkirk for help with the references. Additionally, we were aided by ChatGPT in a variety of tasks, including coding, proofreading, generating conjectures, and assisting in the development of proofs, and we gratefully acknowledge the valuable support ChatGPT provided throughout the process.

This research includes computations using the computational cluster Katana (Katana computational cluster DOI: \url{https://doi.org/10.26190/669X-A286}) supported by Research Technology Services at UNSW Sydney.

DT acknowledges the support of the ARC Future Fellowship FT230100489 and reflects on 35 years of questionable life choices, with more undoubtedly to come.

\section{Background in a nutshell}\label{S:Notation}

We assume that the reader is familiar with the basics of knot theory and monoidal categories, who may refer to {\eg} \cite{Ad-knots,EtGeNiOs-tensor-categories,Tu-qt}. Moreover, in the interest of brevity, instead of recalling quantum invariants in detail we refer
the reader to the wide literature on them, see \cite{ReTu-invariants-3-manifolds-qgroups,Tu-qgroups-3mfds,Tu-qt} for the uncategorified variants and \cite{Kh-cat-jones,BaNa-categorification-jones} for Khovanov homology. Instead, let us rather summarize the conventions we used for the calculations.

\subsection{Comments on the nomenclature}

\begin{Notation}
The quantum invariants of interest in this paper are called, for short, \emph{A2}, \emph{Alexander} (A), \emph{B1}, \emph{Jones} (J), \emph{Khovanov} (K). All of these are Laurent polynomials in $\Z[q,q^{-1}]$ or $\N[q,q^{-1},t,t^{-1}]$ (with $q$ and $t$ formal variables), but we, as often in the literature, just call them polynomials. The $t$ is such that the Khovanov homology specialized at $t=-1$ gives the Jones polynomial. We also have the $t=1$ specialization of Khovanov homology, called \emph{KhovanovT1} (KT1), which lives in $\N[q,q^{-1}]$. \emph{Dequantization} is the process of setting $q=1$, which makes all these invariants trivial.
\end{Notation}

All of our main invariants, and all quantum invariants, are associated with a root system (a Dynkin diagram) with simple roots $\alpha_{1},\dots,\alpha_{r}$ (nodes) and a dominant weight $a_{1}\alpha_{1}+\dots+a_{i}\alpha_{i}+\dots+a_{r}\alpha_{r}$ in $\N\alpha_{1}+\dots+\N\alpha_{r}$ (tuples $(a_{1},\dots,a_{r})$). For example, J comes from the root system of type A1 and the dominant weight is just the unique simple root itself.
In this notation:
\begin{enumerate}

\item \emph{Coloring} is the process of replacing 
$(a_{1},\dots,a_{i},\dots,a_{r})$ with $(a_{1},\dots,a_{i}+1,\dots,a_{r})$.

\item \emph{Rank increase} is the process of replacing $r$ by $r+1$ without leaving the type, whenever applicable.

\item \emph{Categorification} is the process of adding a superscript, {\eg} the A1 invariant becomes the A1${}^{c}$ invariant. The mathematics is explained in \cite{We-knot-invariants}.

\item \emph{Leaving the realm of Lie algebras} is ill-defined. One possibility is taking a Dynkin diagram for a Lie superalgebra as in, for example, \cite[Proposition 2.5.6]{Ka-lie-super}.

\end{enumerate}

\begin{Remark}
Categorification, for us, includes Khovanov homology but not, for example, \emph{odd Khovanov homology} (an alternative categorification of the A1 invariant, see \cite{OzRaSz-odd}) or \emph{HFK} (a categorification of the isotropic A1 invariant, see \cite{HFK1,HFK2}). However, our observations hold for odd Khovanov homology as well, see \autoref{SS:Expect}, and also most likely for HFK. We expect the overall behavior to be essentially independent of the precise form of categorification that one studies.
\end{Remark}

In more detail, let us focus on the case where $\mathfrak{g}$ is a simple Lie algebra. Most of the following also work for Lie superalgebras and other objects, but we will not use them much in this paper.

Let $V$ be a simple highest weight representation of $\mathfrak{g}$.
The Lie algebra $\mathfrak{g}$ is of \emph{classical type} if its Weyl group is of type A, B, C or D. For these, increasing the rank means going to the simple Lie algebra $\mathfrak{h}$ of the same type by adding a vertex as follows.
\begin{gather*}
\mathfrak{g}\colon\dynkin A{4}\Rightarrow\mathfrak{h}\colon\dynkin A{5},
\quad
\dynkin B{4}\Rightarrow\dynkin B{5},
\quad
\dynkin C{4}\Rightarrow\dynkin C{5},
\quad
\dynkin D{4}\Rightarrow\dynkin D{5}.
\end{gather*}
Here, $\mathfrak{g}$ is always on the left, $\mathfrak{h}$ on the right.
This induces an injection of highest weights from $\mathfrak{g}$ to $\mathfrak{h}$ by $(a_{1},\dots,a_{r})\mapsto (0,a_{1},\dots,a_{r})$, and it makes sense to take the representation $V^{\prime}$ associated to the image of 
$(a_{1},\dots,a_{r})$ for $\mathfrak{h}$. Increasing the rank (rank for short) for an associated $Q_{\mathfrak{g},V,\epsilon}$ means going to $Q_{\mathfrak{h},V^{\prime},\epsilon}$.
Coloring (color for short) means going from $Q_{\mathfrak{g},V,\epsilon}$ to $Q_{\mathfrak{g},W,\epsilon}$, where $V$ is of highest weight $(a_{1},\dots,a_{i},\dots,a_{r})$ and $W$ of highest weight $(a_{1},\dots,a_{i}+1,\dots,a_{r})$, and categorification (cat for short) means going from $\epsilon=0$ to $\epsilon=1$.

\begin{Remark}
Note that $(\mathfrak{sl}_{2},\mathrm{Sym}^{2}\C^{2},0)$ is the quantum invariant associated with $SO_{3}$, see {\eg}, \cite[Lemma 5A.7]{Tu-web-reps} for an explanation. However, we count this as a coloring instead of changing the Lie type. Still, we will call it the B1 invariant.
\end{Remark}

\begin{Remark}\label{R:Alexander}
The root system of type $\mathfrak{gl}_{1|1}$ consists of one isotropic root (isotropic roots can only appear for Lie superalgebras), hence the name isotropic A1 invariant.
The Alexander polynomial is thus a ``special'' quantum invariant when compared to the others, but it can also be constructed from an R-matrix; see, for example, \cite{Sa-alexander} for a nice summary.
\end{Remark}

Since quantum invariants behave well with respect to direct sums of representations, restricting to simple $\mathfrak{g}$-representations is a weak restriction.

\subsection{Prime knots and PD presentations}\label{SS:List}

We used a \emph{list of prime knots} up to $16$ crossings that can be found in \cite{knotinfo,TuZh-quantum-big-data-code}.

Throughout, let $n$ be the \emph{number of crossings}.

\begin{Remark}
The sequences $(k_{n})_{n=3}^{\infty}$ and $(\sum_{i=3}^{n}k_{i})_{n=3}^{\infty}$ are
\begin{align*}
\fcolorbox{orchid!50}{white}{\mystrut1,1,2,3,7,21,49,165,552,2176,9988,46972,253293,1388705},
&\quad
\text{for \# of crossings from $3$ to $16$},
\\
\fcolorbox{orchid!50}{white}{\mystrut1,2,4,7,14,35,84,249,801,2977,12965,59937,313230,1701935},
&\quad
\text{for \# of crossings from $\leq 3$ to $\leq 16$},
\end{align*}
see \cite[A002863]{Oeis}. Starting with the trefoil, every knot is listed explicitly. 
This is also the list we used for the big data comparison, {\ie} we do not distinguish a knot from its mirror.
\end{Remark}

Let us summarize the conventions.

\begin{enumerate}

\item The list contains \emph{planar diagrams} (PD for short) for each prime knot. 
Recall that a PD presentation of a knot diagram labels all of its edges with non-repeating numbers $\{1,\dots,r\}$, where $r$ is the number of edges. Each edge is then adjacent to two crossings, which induces a labeling of the crossings. We remember the crossings as symbols $X[i,j,k,l]$, where $i$, $j$, $k$ and $l$ are the labels of the edges around that crossing, starting from the incoming lower edge and proceeding counterclockwise. Explicitly, the first two knots on the list and their PD presentations are:
\begin{gather*}
PD[{\color{spinach}X[1, 5, 2, 4]},{\color{orchid}X[3, 1, 4, 6]},{\color{tomato}X[5, 3, 6, 2]}]
\leftrightsquigarrow
\begin{tikzpicture}[anchorbase]
\node at (0,0) {\includegraphics[width=0.2\textwidth]{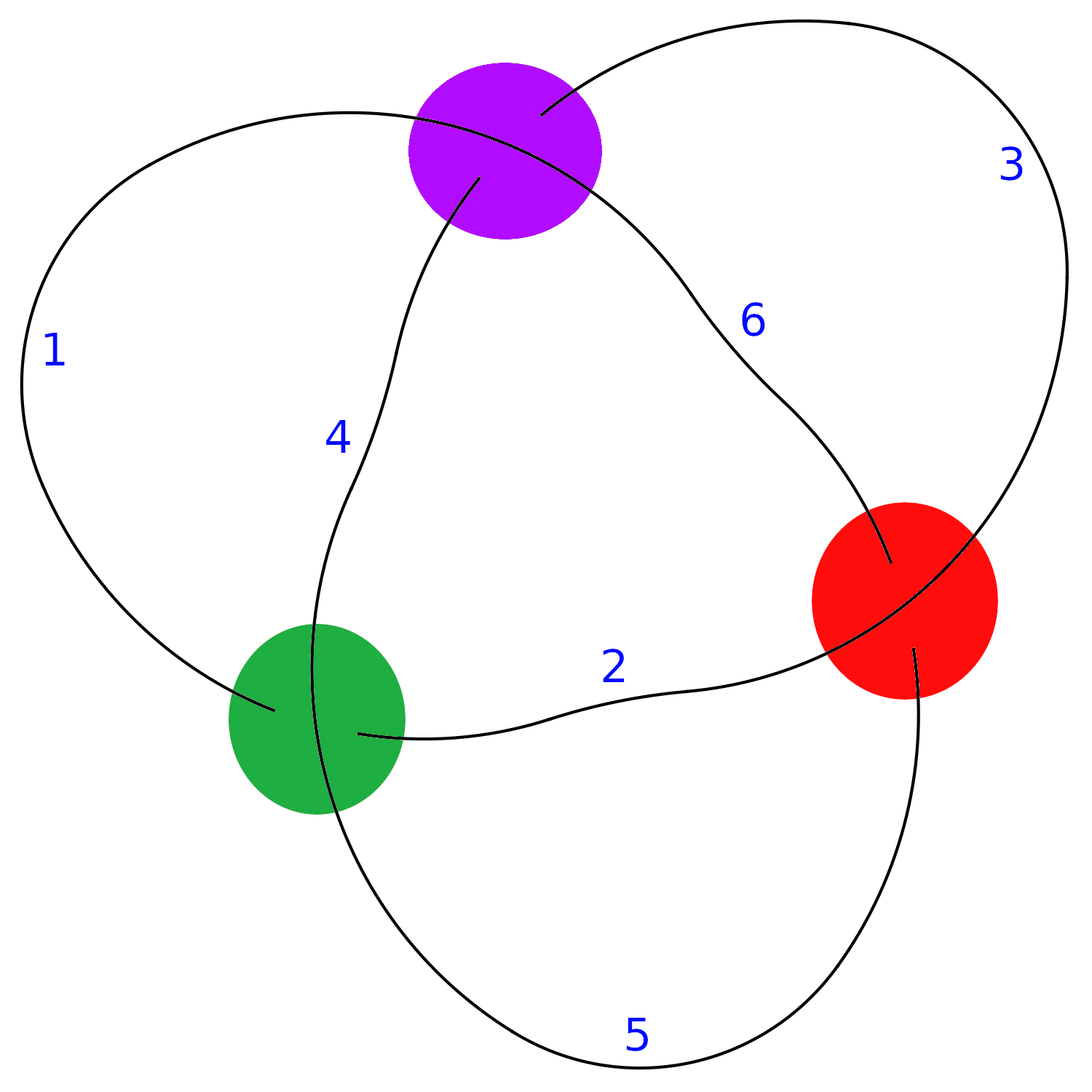}};
\end{tikzpicture}
,
\\
PD[X[4, 2, 5, 1], X[8, 6, 1, 5], X[6, 3, 7, 4], X[2, 7, 3, 8]]
\leftrightsquigarrow
\begin{tikzpicture}[anchorbase]
\node at (0,0) {\includegraphics[width=0.2\textwidth]{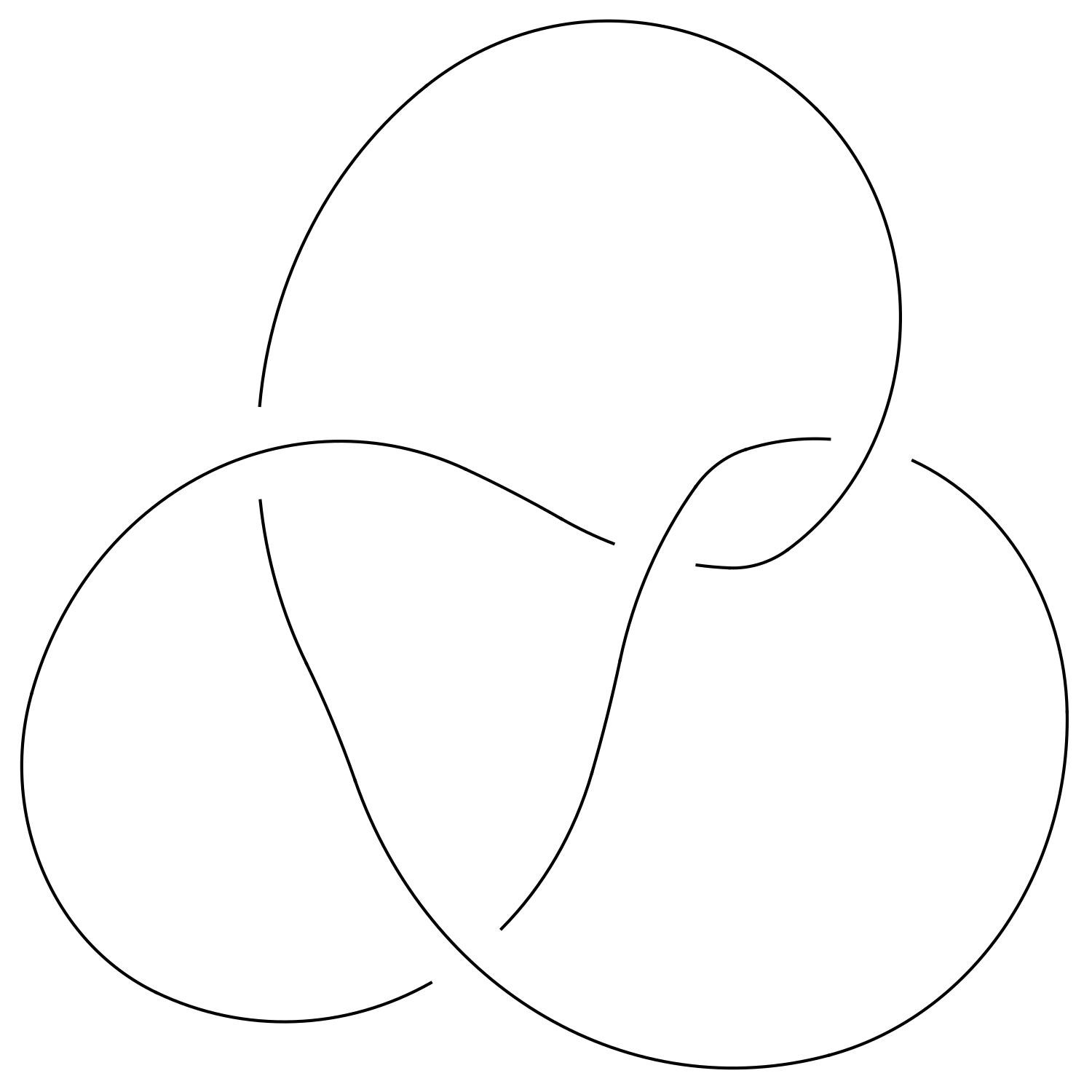}};
\end{tikzpicture}
.
\end{gather*}
We only added the labels and highlighted the crossings for the \emph{trefoil} knot at the top. The bottom knot is the \emph{figure eight knot}.

\item In particular, the mirror image of the trefoil is not on the list. In general, the list identifies knots and their mirror images and has only the PD presentation for one of them.

\item Note that {\eg}, the Jones polynomial $J$ satisfies $J(K\# L)=J(K)
J(L)$, so in some sense the running over prime knots suffices for our purposes. 

\end{enumerate}

\subsection{Notation and Conventions}\label{SS:Notation}

We now briefly explain the setting we use throughout, that is the same as in \cite{LaTuVa-big-data}. 
First, we collect the various polynomials as 
\emph{lists} (ordered, repetitions allowed, notation $[\placeholder]$).
If not stated otherwise, we will identify a 
polynomial with a vector as follows. We use some upper bound $d$ for the absolute degree of polynomials appearing. Then
\begin{gather*}
a_{-d}q^{-d}+\dots+a_{d}q^{d}
\leftrightsquigarrow
[a_{-d},\dots,a_{d}]\in\Z^{2d+1}\subset\R^{2d+1},
\end{gather*}
where we assign zero to the coefficient of included terms not appearing in the polynomial. In this way, all polynomials are in the same vector space.
Similarly, for a two-variable polynomial we get a matrix containing the various coefficients, and we \emph{flatten} it into a vector.

\begin{Notation}
For a vector $v$, we write $\sum v$ for the sum over all entries of a vector, $\optt{max}(v)$ for the maximal entry, $\optt{abs}(v)$ for the componentwise absolute value of the vector, and $\optt{roots}(v)$ for the set of roots of the corresponding normalized polynomial. Additionally we write $\optt{span}(v)$ for the number of entries ignoring padded zeros, {\ie} zeros at the beginning and end of a vector. This corresponds to the difference between the maximal and minimal degree of a given polynomial. Naturally, $\optt{span}$ only applies to single variable polynomials.
\end{Notation}

Let $Q$ denote a quantum invariant and let
\begin{gather*}
\mathcal{K}_{n}=[K|K\text{ a prime knot with $\leq$ n crossings}],
\end{gather*}
be the (ordered) \emph{list of prime knots} as in \autoref{SS:List}. 
Denoting by $Q(K)$ the value of $Q$ on the knot $K$, we also study
\begin{gather*}
\mathcal{Q}_{n}=[Q(K)|K\in\mathcal{K}_{n}],
\end{gather*}
the \emph{list of polynomials} for $n$ or fewer crossings.

\begin{Example}
To be completely explicit, let $n=4$. Then
\begin{gather*}
\mathcal{K}_{4}=\big[
PD[X[1, 5, 2, 4],X[3, 1, 4, 6],X[5, 3, 6, 2]]
,
PD[X[4, 2, 5, 1], X[8, 6, 1, 5], X[6, 3, 7, 4], X[2, 7, 3, 8]]
\big]
\end{gather*}
as in \autoref{SS:List} above. For the Jones polynomial $Q=Q_{\mathfrak{sl}_{2},\C^{2},0}$
we then have (for $d=4$)
\begin{gather*}
\mathcal{Q}_{4}=
\big[
[0,0,0,0,0,1,0,1,-1]
,
[0,0,1,-1,1,-1,1,0,0]
\big],
\end{gather*}
as one may easily check (or find in \cite{TuZh-quantum-big-data-code} or many other places).
\end{Example}

Throughout, we are interested in \emph{large $n$ behavior}. We fix some notation that we will use.

\begin{Notation}
For functions $f,g\colon\N\to\R_{>0}$ we use:
\begin{gather*}
\begin{aligned}
f\in\Omega(g)
&\;\Leftrightarrow\;
\exists C>0,\,\exists n_{0}
\;\text{ such that }\;
\fbox{$|f(n)|\geq C\cdot g(n)$},
\;\forall n>n_{0}
,
\\
f\in O(g)
&\;\Leftrightarrow\;
\exists C>0,\,\exists n_{0}
\;\text{ such that }\;
\fbox{$|f(n)|\leq C\cdot g(n)$},
\;\forall n>n_{0}
.
\end{aligned}
\end{gather*}
We use these, in order, as \emph{asymptotic lower and upper bounds}.
\end{Notation}

\begin{Notation}
To analyze functions that behave like $a^{n}$ for some $a\in\R_{>1}$ we use
\emph{successive quotients}. For a function $f\colon\N\to\R_{>0}$, this is the function
$n\mapsto \frac{f(n+1)}{f(n)}$, whose value effectively approximates $a$ for $f(n)=a^{n}$.
\end{Notation}

As a final note, particularly for the data below,
we sometimes use a \emph{log plot} for displaying data, which means that we use a logarithmic scale on the $y$-axis, and when we show floating point numbers, they are \emph
{truncated} ({\ie} a floor) to the appropriate decimal precision. For example, the numbers $0.6669$ and $0.6661$ are written as $0.666$ in three-decimal precision.

\subsection{Computation of quantum invariants}

We now explain the different computational techniques that we have used to compute the quantum invariants.

\smallskip
\noindent{\color{mygray}\rule[1.0ex]{\linewidth}{0.5pt}}

\textit{For A2, A, B1, J; used for A2, B1, J.}
The first method is using \emph{R-matrices}.

We work over $\C(q)$ for a formal variable $q$ that potentially also has roots like $q^{1/k}$.
The uncategorified quantum invariants $Q=Q_{\mathfrak{g},V,0}$ that we use, except the Alexander polynomial, can be defined and computed as follows. For a knot, fix a Morse presentation of the knot, arranged vertically. 

For simplicity, assume that $V$ is self-dual, which works for B1 and J. In the Morse presentation, we have four basic pieces and an identity, that we name as
\begin{gather*}
R=
\begin{tikzpicture}[anchorbase,scale=1]
\draw[usual] (0.5,0) to[out=90,in=270] (0,0.5);
\draw[usual,crossline] (0,0) to[out=90,in=270] (0.5,0.5);
\end{tikzpicture}
,\quad
R^{-1}=
\begin{tikzpicture}[anchorbase,scale=1]
\draw[usual] (0,0) to[out=90,in=270] (0.5,0.5);
\draw[usual,crossline] (0.5,0) to[out=90,in=270] (0,0.5);
\end{tikzpicture}
,\quad
cap=
\begin{tikzpicture}[anchorbase,scale=1]
\draw[usual] (0,0) to[out=90,in=180] (0.25,0.25) to[out=0,in=90] (0.5,0);
\end{tikzpicture}
,\quad
cup=
\begin{tikzpicture}[anchorbase,scale=1]
\draw[usual] (0.5,0) to[out=270,in=0] (0.25,-0.25) to[out=180,in=270] (0,0);
\end{tikzpicture}
,\quad
id=
\begin{tikzpicture}[anchorbase,scale=1]
\draw[usual] (0.5,0) to[out=90,in=270] (0.5,0.5);
\end{tikzpicture}
\ .
\end{gather*}
We associate these to linear maps (matrices upon choice of basis) denoted with the same symbols
\begin{gather}\label{Eq:maps}
R,R^{-1}\colon V_{q}\otimes V_{q}\to V_{q}\otimes V_{q}
,\quad
cap\colon V_{q}\otimes V_{q}\to\C(q)
,\quad
cup\colon \C(q)\to V_{q}\otimes V_{q}
,\quad
id\colon V_{q}\to V_{q},v\mapsto v
.
\end{gather}
Here, $V_{q}$ is a representation of a quantum group associated with $\mathfrak{g}$ that dequantizes to the $\mathfrak{g}$ representation $V$, and all the above maps are intertwiners for the quantum group that dequantize to, respective to above, the flip maps, the pairing and coparing of $\mathfrak{g}$ representations, and the identity.

Horizontal composition is then the tensor product (Kronecker product upon choice of basis).
Write $id\otimes\dots\otimes id$ with $k$ factors simply as $k$.
Then for the figure eight knot this is could be
\begin{gather*}
\begin{tikzpicture}[anchorbase,scale=1]
\draw[usual] (0.5,0) to[out=90,in=270] (0,0.5);
\draw[usual,crossline] (0,0) to[out=90,in=270] (0.5,0.5);
\draw[usual] (1,0) to[out=90,in=270] (1,0.5);
\draw[usual] (0,0.5) to[out=90,in=270] (0,1);
\draw[usual] (0.5,0.5) to[out=90,in=270] (1,1);
\draw[usual,crossline] (1,0.5) to[out=90,in=270] (0.5,1);
\draw[usual] (0.5,1) to[out=90,in=270] (0,1.5);
\draw[usual,crossline] (0,1) to[out=90,in=270] (0.5,1.5);
\draw[usual] (1,1) to[out=90,in=270] (1,1.5);
\draw[usual] (0,1.5) to[out=90,in=270] (0,2);
\draw[usual] (0.5,1.5) to[out=90,in=270] (1,2);
\draw[usual,crossline] (1,1.5) to[out=90,in=270] (0.5,2);
\draw[usual] (1,2) to[out=90,in=180] (1.25,2.25) to[out=0,in=90] (1.5,2) to (1.5,0) to[out=270,in=0] (1.25,-0.25) to[out=180,in=270] (1,0);
\draw[usual] (0.5,2) to[out=90,in=180] (1.25,2.5) to[out=0,in=90] (2,2) to (2,0) to[out=270,in=0] (1.25,-0.5) to[out=180,in=270] (0.5,0);
\draw[usual] (0,2) to[out=90,in=180] (1.25,2.75) to[out=0,in=90] (2.5,2) to (2.5,0) to[out=270,in=0] (1.25,-0.75) to[out=180,in=270] (0,0);
\end{tikzpicture}
\leftrightsquigarrow
\scalebox{0.65}{$\begin{tikzcd}[ampersand replacement=\&,column sep=0.15em,row sep=0.5em]
0\otimes cap \otimes 0
\ar[d,-,"\circ"]
\\
1\otimes cap\otimes 1
\ar[d,-,"\circ"]
\\
2\otimes cap\otimes 2
\ar[d,-,"\circ"]
\\
1\otimes R^{-1}\otimes 3
\ar[d,-,"\circ"]
\\
0\otimes R\otimes 4
\ar[d,-,"\circ"]
\\
1\otimes R^{-1}\otimes 3
\ar[d,-,"\circ"]
\\
0\otimes R\otimes 4
\ar[d,-,"\circ"]
\\
2\otimes cup\otimes 2
\ar[d,-,"\circ"]
\\
1\otimes cup\otimes 1
\ar[d,-,"\circ"]
\\
0\otimes cup \otimes 0
\end{tikzcd}$}
,
\end{gather*}
with composition, for example from bottom to top.

The maps in \autoref{Eq:maps} are not uniquely determined from what we wrote above. There are some choices involved, but they are mostly irrelevant and just rescale the quantum invariant. To be completely explicit, we used the conventions used in \cite{KnotTheory}. For example, in coordinates,
\begin{gather*}
R=
\begin{pmatrix}
q^{1/2} & 0 & 0 & 0
\\
0 & 0 & q & 0
\\
0 & q & q^{1/2}-q^{3/2} & 0
\\
0 & 0 & 0 & q^{1/2}
\end{pmatrix}
,
\end{gather*}
is the choice of $R$-matrix for the Jones polynomial.

The Alexander polynomial can be computed using R-matrices, with the slight catch that one needs to renormalize the result appropriately as one gets zero, {\cf} \cite{Sa-alexander}, which is a common phenomena when working with isotropic roots.

\smallskip
\noindent{\color{mygray}\rule[1.0ex]{\linewidth}{0.5pt}}

\textit{For A2, A, B1, J, K, KT1; used for B1, K, KT1.}
Another way of computing quantum invariants is a \emph{skein theory approach}.

For Jones this, for example, uses the relations:
\begin{gather*}
\begin{tikzpicture}[anchorbase,scale=1]
\draw[usual] (0.5,0) to [out=90,in=270] (0,0.5);
\draw[usual,crossline] (0,0) to [out=90,in=270] (0.5,0.5);
\end{tikzpicture}
=
q^{1/2}\cdot
\begin{tikzpicture}[anchorbase,scale=1]
\draw[usual] (0.5,0) to [out=90,in=270] (0.5,0.5);
\draw[usual] (0,0) to [out=90,in=270] (0,0.5);
\end{tikzpicture}
+
q^{-1/2}\cdot
\begin{tikzpicture}[anchorbase,scale=1]
\draw[usual] (0,0.5) to [out=270,in=180] (0.25,0.3) to [out=0,in=270] (0.5,0.5);
\draw[usual] (0,0) to [out=90,in=180] (0.25,0.2) to [out=0,in=90] (0.5,0);
\end{tikzpicture}
,
\quad
\begin{tikzpicture}[anchorbase,scale=1]
\draw[usual] (0,0) to [out=90,in=270] (0.5,0.5);
\draw[usual,crossline] (0.5,0) to [out=90,in=270] (0,0.5);
\end{tikzpicture}
=
q^{-1/2}\cdot
\begin{tikzpicture}[anchorbase,scale=1]
\draw[usual] (0.5,0) to [out=90,in=270] (0.5,0.5);
\draw[usual] (0,0) to [out=90,in=270] (0,0.5);
\end{tikzpicture}
+
q^{1/2}\cdot
\begin{tikzpicture}[anchorbase,scale=1]
\draw[usual] (0,0.5) to [out=270,in=180] (0.25,0.3) to [out=0,in=270] (0.5,0.5);
\draw[usual] (0,0) to [out=90,in=180] (0.25,0.2) to [out=0,in=90] (0.5,0);
\end{tikzpicture}
,\quad
\begin{tikzpicture}[anchorbase,scale=1]
\draw[usual] (0,0) to [out=270,in=180] (0.25,-0.25) to [out=0,in=270] (0.5,0);
\draw[usual] (0,0) to [out=90,in=180] (0.25,0.25) to [out=0,in=90] (0.5,0);
\end{tikzpicture}
=
-(q+q^{-1}).
\end{gather*}
One can then replace all crossings, as usual, and get the desired polynomial. We exemplify this here for the Hopf link:
\begin{align*}
\begin{tikzpicture}[anchorbase,scale=1]
\draw[usual] (0.5,0) to [out=90,in=270] (0,0.5);
\draw[usual,crossline] (0,0) to [out=90,in=270] (0.5,0.5);
\draw[usual] (0.5,0.5) to [out=90,in=270] (0,1);
\draw[usual,crossline] (0,0.5) to [out=90,in=270] (0.5,1);
\draw[usual] (0.5,1) to[out=90,in=180] (0.75,1.25) to[out=0,in=90] (1,1) to (1,0) to[out=270,in=0] (0.75,-0.25) to[out=180,in=270] (0.5,0);
\draw[usual] (0,1) to[out=90,in=180] (0.75,1.5) to[out=0,in=90] (1.5,1) to (1.5,0) to[out=270,in=0] (0.75,-0.5) to[out=180,in=270] (0,0);
\end{tikzpicture}
&=
q\cdot
\begin{tikzpicture}[anchorbase,scale=1]
\draw[usual] (0,0) to (0,0.5);
\draw[usual] (0.5,0) to (0.5,0.5);
\draw[usual] (0,0.5) to (0,1);
\draw[usual] (0.5,0.5) to (0.5,1);
\draw[usual] (0.5,1) to[out=90,in=180] (0.75,1.25) to[out=0,in=90] (1,1) to (1,0) to[out=270,in=0] (0.75,-0.25) to[out=180,in=270] (0.5,0);
\draw[usual] (0,1) to[out=90,in=180] (0.75,1.5) to[out=0,in=90] (1.5,1) to (1.5,0) to[out=270,in=0] (0.75,-0.5) to[out=180,in=270] (0,0);
\end{tikzpicture}
+
\begin{tikzpicture}[anchorbase,scale=1]
\draw[usual] (0,0) to (0,0.5);
\draw[usual] (0.5,0) to (0.5,0.5);
\draw[usual] (0,0.5) to [out=90,in=180] (0.25,0.7) to [out=0,in=90] (0.5,0.5);
\draw[usual] (0,1) to [out=270,in=180] (0.25,0.8) to [out=0,in=270] (0.5,1);
\draw[usual] (0.5,1) to[out=90,in=180] (0.75,1.25) to[out=0,in=90] (1,1) to (1,0) to[out=270,in=0] (0.75,-0.25) to[out=180,in=270] (0.5,0);
\draw[usual] (0,1) to[out=90,in=180] (0.75,1.5) to[out=0,in=90] (1.5,1) to (1.5,0) to[out=270,in=0] (0.75,-0.5) to[out=180,in=270] (0,0);
\end{tikzpicture}
+
\begin{tikzpicture}[anchorbase,scale=1]
\draw[usual] (0,0) to [out=90,in=180] (0.25,0.2) to [out=0,in=90] (0.5,0);
\draw[usual] (0,0.5) to [out=270,in=180] (0.25,0.3) to [out=0,in=270] (0.5,0.5);
\draw[usual] (0,0.5) to (0,1);
\draw[usual] (0.5,0.5) to (0.5,1);
\draw[usual] (0.5,1) to[out=90,in=180] (0.75,1.25) to[out=0,in=90] (1,1) to (1,0) to[out=270,in=0] (0.75,-0.25) to[out=180,in=270] (0.5,0);
\draw[usual] (0,1) to[out=90,in=180] (0.75,1.5) to[out=0,in=90] (1.5,1) to (1.5,0) to[out=270,in=0] (0.75,-0.5) to[out=180,in=270] (0,0);
\end{tikzpicture}
+
q^{-1}\cdot
\begin{tikzpicture}[anchorbase,scale=1]
\draw[usual] (0,0) to [out=90,in=180] (0.25,0.2) to [out=0,in=90] (0.5,0);
\draw[usual] (0,0.5) to [out=270,in=180] (0.25,0.3) to [out=0,in=270] (0.5,0.5);
\draw[usual] (0,0.5) to [out=90,in=180] (0.25,0.7) to [out=0,in=90] (0.5,0.5);
\draw[usual] (0,1) to [out=270,in=180] (0.25,0.8) to [out=0,in=270] (0.5,1);
\draw[usual] (0.5,1) to[out=90,in=180] (0.75,1.25) to[out=0,in=90] (1,1) to (1,0) to[out=270,in=0] (0.75,-0.25) to[out=180,in=270] (0.5,0);
\draw[usual] (0,1) to[out=90,in=180] (0.75,1.5) to[out=0,in=90] (1.5,1) to (1.5,0) to[out=270,in=0] (0.75,-0.5) to[out=180,in=270] (0,0);
\end{tikzpicture}
\\
&=
q(q+q^{-1})^{2}
-2(q+q^{-1})
+q^{-1}(q+q^{-1})^{2}
=
q^{3}+q+q^{-1}+q^{-3}.
\end{align*}
We never used this for the decategorified quantum invariants, except for B1. We elaborate in \cite{TuZh-quantum-big-data-code}, for this paper it is enough to know that for B1 we followed the conventions in \cite[Section 5]{Tu-web-reps}.

We do not know a general R-matrix strategy to compute categorified quantum invariants. For Khovanov we used \cite{KnotTheory}, which is based on \cite{BaNa-fast} and uses the conventions from \cite{BaNa-categorification-jones}. In particular, 
it uses a skein theory approach where every crossing is a complex, in a certain category, of the form
\begin{gather*}
\begin{tikzpicture}[anchorbase,scale=1]
\draw[usual] (0.5,0) to[out=90,in=270] (0,0.5);
\draw[usual,crossline] (0,0) to[out=90,in=270] (0.5,0.5);
\end{tikzpicture}
\leftrightsquigarrow
0\to
\underline{\begin{tikzpicture}[anchorbase,scale=1]
\draw[usual] (0,0) to[out=90,in=180] (0.25,0.2) to[out=0,in=90] (0.5,0);
\draw[usual] (0.5,0.5) to[out=270,in=0] (0.25,0.3) to[out=180,in=270] (0,0.5);
\end{tikzpicture}}
\xrightarrow{d}
q\cdot
\begin{tikzpicture}[anchorbase,scale=1]
\draw[usual] (0.5,0) to[out=90,in=270] (0.5,0.5);
\draw[usual] (0,0) to[out=90,in=270] (0,0.5);
\end{tikzpicture}
\to 0
,
\end{gather*}
for a certain differential $d$ and underlined part in homological degree zero. The whole construction is then tensored together, as complexes, over the crossings. The rest of the calculation is similar to the Jones calculation, just in complexes.

\smallskip
\noindent{\color{mygray}\rule[1.0ex]{\linewidth}{0.5pt}}

\textit{For A; used for A.}
The following is the \emph{determinant approach} which generally does not work for quantum invariants. It inspired ``pseudo quantum invariants'' (our choice of name) as in \cite{BaNaVe-polynomial-time-invariant}.

The arguably best approach to computing the A invariant is to use the \emph{Seifert matrix} $S(K)$ of a knot (which in turn collects the linking numbers of the Seifert circles). Precisely,
\begin{gather*}
A(K)
=
\det(S(K)-qS(K)^{T}).
\end{gather*}
The example to keep in mind is the trefoil, where the Seifert matrix is $\begin{psmallmatrix}1 & 0\\-1 & 1\end{psmallmatrix}$. Thus:
\begin{gather*}
A(PD[X[1, 5, 2, 4],X[3, 1, 4, 6],X[5, 3, 6, 2]])
=\det
\left(
\begin{psmallmatrix}1 & 0\\-1 & 1\end{psmallmatrix}
-q
\begin{psmallmatrix}1 & -1\\0 & 1\end{psmallmatrix}
\right)
=
\det\left(
\begin{psmallmatrix}1-q & q\\-1 & 1-q\end{psmallmatrix}
\right)
=1-q+q^{2}.
\end{gather*}
This only determines A up to a scalar, and we scale it so that the result is the same when one would use a skein approach, as in \cite{KnotTheory}. 

\subsection{Algorithmic complexity}\label{SS:Algorithm}

We will now give a rough estimate of the complexity of the above algorithms.
We go through the above list in slightly changed order.

\begin{Remark}
We are analyzing specific algorithms, not the problem of computing quantum invariants themselves, which might be very different. We also do not try to give the best upper bounds possible.
\end{Remark}

\smallskip
\noindent{\color{mygray}\rule[1.0ex]{\linewidth}{0.5pt}}

\textit{Determinant approach.} In this scenario the Seifert matrix $S(K)$ is of size $2g$-by-$2g$, where $g$ is the genus of the diagram of the knot $K$. Note that $g\leq n/2$, so that the Seifert matrix is at most of size $n$-by-$n$. The main computational complexity is then the computation of the determinant, which, using LU decomposition, is in $O(n^{3})$.
Thus:
\begin{gather*}
A\in O(n^{3}).
\end{gather*}
This is very cheap compared to what is up next.

\smallskip
\noindent{\color{mygray}\rule[1.0ex]{\linewidth}{0.5pt}}

\textit{R-matrix approach.} Let $N=\dim_{\C}V$, and let $p(n)$ denote some polynomial in $n$, and let $m$ be the carving width of the four-valent planar graph associated with a knot. With a bit of care, see \cite[Theorem 1.1]{Ma-complexity-quantum}, one can show that $m\in O(\sqrt{n})$ and then
\begin{gather*}
Q_{\mathfrak{g},V,0}\in O(p(n)N^{3m/2})=O(p(n)N^{3\sqrt{n}/2}).
\end{gather*}
Thus, computing quantum invariants using $R$-matrices is superpolynomial in $n$, with the leading factor determined by
the dimension of the underlying representation. This is not surprising if one keeps in mind that the main difficulty in this computational approach is not the number of crossings, but rather the number of strands, since this corresponds to tensoring $V$, which, even in decompositions \cite{CoOsTu-growth,KhSiTu-monoidal-cryptography}, behaves exponentially in $N$.
The algorithm explained in \cite{Ma-complexity-quantum} is not quite the one used in KnotTheory \cite{KnotTheory}, but the one in KnotTheory should have roughly the same complexity (assuming that the algorithm that produces a Morse presentation is optimized; which is not the case for the program used by \cite{KnotTheory}).

\begin{Remark}
Finding a presentation of a knot with small carving width is surprisingly inexpensive, namely polynomial in $n$, see
\cite{SeTh-ratcatcher}.
\end{Remark}

\begin{Remark}
It is known that computing the Jones polynomial is {\#}P-hard, see, for example, \cite[Section 6.3]{We-complexity}. Thus, under the standard assumptions that complexity classes do not collapse, the R-matrix approach is among the fastest possible.
\end{Remark}

\smallskip
\noindent{\color{mygray}\rule[1.0ex]{\linewidth}{0.5pt}}

\textit{Skein theory approach.}
Let $M$ denote the number of summands in the expression of the crossing, {\eg} $M=2$ for the Jones polynomial or Khovanov homology and $M=3$ for the B1 invariant.
As before, let $p(n)$ denote some polynomial in $n$ and $s\in\R_{\geq 1}$ be some scalar. Then
\begin{gather*}
Q_{\mathfrak{g},V,\epsilon}\in O(p(n)M^{sn}).
\end{gather*}
In particular, since, to the best of our knowledge, there is no R-matrix approach to Khovanov homology, the computation of it is exponential in the number of crossings. 

\begin{Remark}
There is evidence that the computation of the categorification
should be of the same complexity as the computation of the $\epsilon=0$ version; see, {\eg}, \cite{PrSi-complexity}.
\end{Remark}

\smallskip
\noindent{\color{mygray}\rule[1.0ex]{\linewidth}{0.5pt}}

\textit{Average runtime.}
Arguably, the \emph{worst case} analysis as above is not the right thing to consider. For example, with the calculations we run 
it seemed that there are few knots that took a very long time to compute, like 1000 times longer than the others with the same number of crossings. So, the correct measurement might be \emph{average runtime}.

We sadly do not know of any sources that computed the average runtime (and, indeed, this is usually much harder than the worst case). We here rather list our experimental data for the average runtime (the second plot is a \emph{box plot}) all of which were run with KnotTheory \cite{KnotTheory} on the same machine: the laptop of the second author. Additionally, we include the data for the Javascript program for B1, run separately on UNSW's Katana servers, which uses the Skein theoretic approach as opposed to the R-matrix approach of the KnotTheory library.
The details of these can be found in \cite{TuZh-quantum-big-data-code}.

\begin{gather*}
\begin{tikzpicture}[anchorbase]
\node at (0,0) {\includegraphics[width=0.6\textwidth]{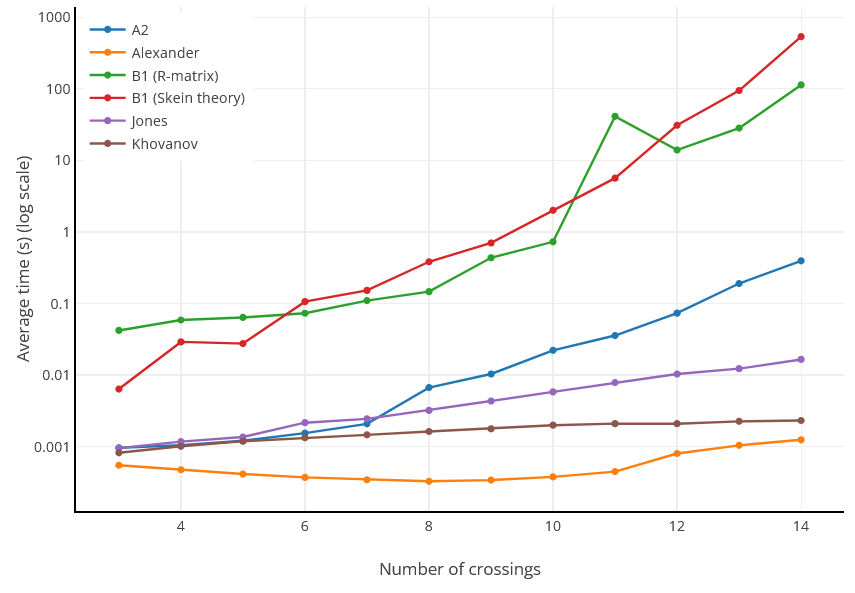}};
\end{tikzpicture}
\end{gather*}
\begin{gather*}
\begin{tikzpicture}[anchorbase]
\node at (0,0) {\includegraphics[width=0.6\textwidth]{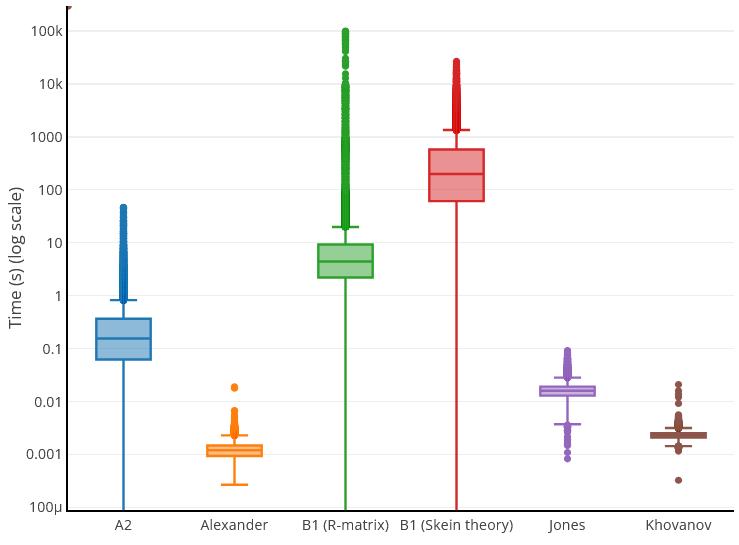}};
\end{tikzpicture}
\end{gather*}
The precise numbers can be found in \autoref{figure:0}.

We have three comments on this:
\begin{enumerate}[label=(\roman*)]

\item The variance is quite high for A2 and B1. Since the runtime of the R-matrix approach has leading factor $N=\dim_{\C}V$, and $N=3$ for these two invariants, the issue might be an not optimal conversion from a PD presentation to a Morse presentation for some knots. Indeed, 
there was knot number 775 on the list (11 crossings) where the computation took 21703.436821 seconds, much longer than any other knot.
This is one reason we wrote a skein theory program for B1 (that can be found in \cite{TuZh-quantum-big-data-code}) and it typically did better on knots that were difficult in the R-matrix approach but was slower on average; see the data above.

\item The computation of A was so fast that the program was still loading when the first few knots were computed. Hence, the drop below $n=8$.

\item The computation of K is not run in Mathematica itself, but instead starts a Java program. This could explain why K was computed faster than J.

\end{enumerate}

\section{Comparison -- distinct values}\label{S:Compare1}

For the first comparison, we are interested in how many distinct values quantum invariants take on our list of knots, and we measure this as a percentage. In other words, we want to know
\begin{gather*}
Q(n)^{\%}=
\#\{Q(K) \mid K\in\mathcal{K}_{n}\}/
\#\mathcal{K}_{n}.
\end{gather*}
These are the \emph{distinct values} $Q$ takes.

Recall that A, J, and K stand for Alexander, Jones, and Khovanov, respectively. Let ``All''
mean that we take all of them together, and we use J+KT1 to take Jones and KT1 together. 
The data is as follows.
\begin{gather*}
\begin{tikzpicture}[anchorbase]
\node at (0,0) {\includegraphics[width=0.6\textwidth]{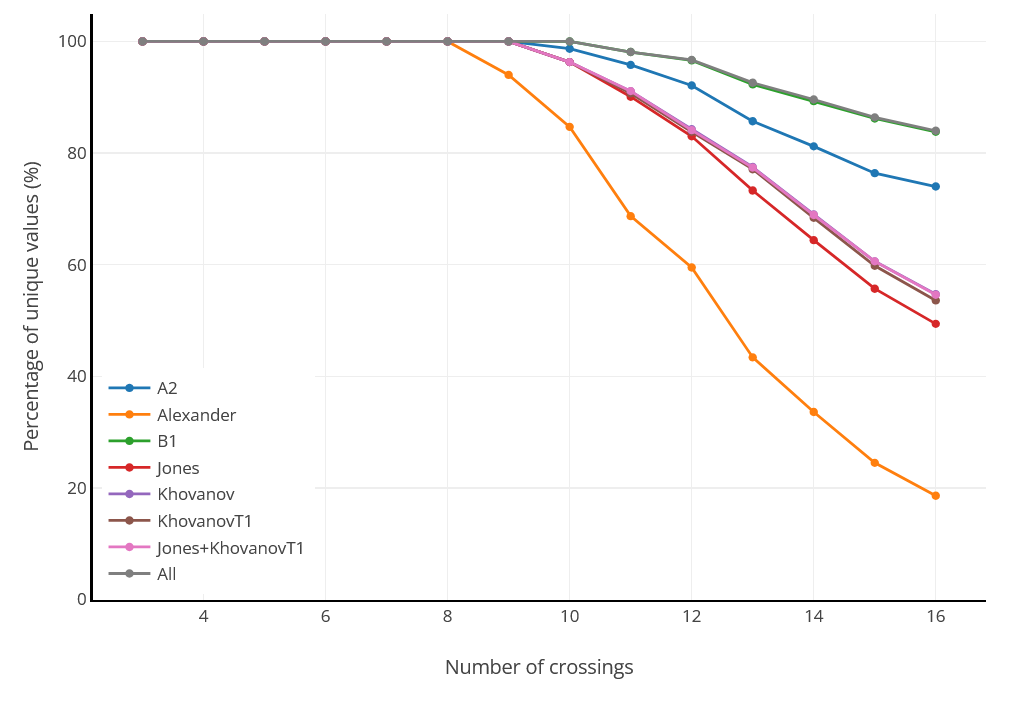}};
\end{tikzpicture}
,
\end{gather*}
with the precise values listed in \autoref{figure:comparison1}.
Here are a few observations.
\begin{enumerate}

\item From the data, \emph{four different classes} are visible: All/B1, A2, J/K/KT1, and A.

\item They all seem to drop to zero, but the rate of convergence to zero seems to depend on the class.

\item There is essentially no difference between K and J+KT1, so categorification seems to be as strong as two evaluations of it. Additionally, there is no huge difference between K and J or KT1.

\item There is essentially no difference between All and B1; in particular, coloring seems to be the preferred method to distinguish knots.

\end{enumerate}

In general, some \emph{zero-one law} should apply, so $Q(n)^{\%}$ should either converge to zero or one. The data suggest that the limit for all of them is zero. From a more detailed analysis of the data we even conclude an \emph{exponential decay}.

\begin{Conjecture}\label{C:Drops}
For $Q\in\{A2,A,B1,K\}$ (and therefore also for $Q=J$ or $Q=KT1$) we have
\begin{gather*}
Q(n)^{\%}\in O(\delta^{n})
\text{ for some }\delta=\delta(Q)\in(0,1).
\end{gather*}
Moreover, we have
\begin{gather*}
\delta(B1)>\delta(A2)>\delta(K)>\delta(A).
\end{gather*}
In other words, they all drop to zero exponentially fast, but the exponential factor depends on the class.
\end{Conjecture}

\begin{Speculation}\label{S:Drops}
For all quantum invariants $Q$ we have
\begin{gather*}
Q(n)^{\%}\in O(\delta^{n})
\text{ for some }\delta=\delta(Q)\in(0,1).
\end{gather*}
In other words, they all drop to zero exponentially fast.
\end{Speculation}

\begin{Speculation}
Assume \autoref{S:Drops} holds and fix $Q_{\mathfrak{g},V,0}$ for $\mathfrak{g}$ of classical type. Then:
\begin{gather*}
\delta(color)>\delta(rank)>\delta(cat).
\end{gather*}
That is, coloring is better than rank increase, which is better than categorification.
\end{Speculation}

\begin{Remark}
Detecting certain knots, such as K detects the unknot \cite{KrMo-unknot-detector}, seems to be only possible if the knots are very special (unknot, torus knots, {\etc}).
\end{Remark}

Factoring in computational complexity, one might argue that A is the best invariant. Furthermore, assuming that 
categorification is actually of the same complexity as its $\epsilon=0$ counterpart, one might argue that 
$cat>color>rank$, since coloring and rank increase have a bigger exponential factor in their complexity analysis, {\cf} \autoref{SS:Algorithm}. This however needs a deeper analysis of $\delta$ which is beyond this paper.

\subsection*{Comments on how to potentially prove parts of \autoref{C:Drops}}\label{SS:Expect}

The rest of the section is dedicated to proof a version of \autoref{C:Drops}.
A \emph{skein relation} is a relation of the following form. Let $a,b,c$ be elements in a ring, then
\begin{gather*}
a\cdot
\begin{tikzpicture}[anchorbase,scale=1]
\draw[usual,directed=1] (0.5,0) to[out=90,in=270] (0,0.5);
\draw[usual,crossline,directed=1] (0,0) to[out=90,in=270] (0.5,0.5);
\end{tikzpicture}
+b\cdot
\begin{tikzpicture}[anchorbase,scale=1]
\draw[usual,directed=1] (0,0) to[out=90,in=270] (0.5,0.5);
\draw[usual,crossline,directed=1] (0.5,0) to[out=90,in=270] (0,0.5);
\end{tikzpicture}
+c\cdot
\begin{tikzpicture}[anchorbase,scale=1]
\draw[usual,directed=1] (0,0) to[out=90,in=270] (0,0.5);
\draw[usual,directed=1] (0.5,0) to[out=90,in=270] (0.5,0.5);
\end{tikzpicture}
=
0
,\quad \text{for } a,b\text{ invertible}.
\end{gather*}
\emph{Multiplicity freeness} of a quantum invariant $Q_{\mathfrak{g},V,0}$ is the property that $V\otimes V$ is multiplicity free.
Examples of quantum invariants that satisfy a skein relation or multiplicity freeness are A2, A, B1, and J.

Let $\mathcal{AL}_{n}$ be the set of alternating links of $\leq n$ crossings. 
Similarly to $Q(n)^{\%}$ define
\begin{gather*}
Q(n)^{\%}_{AL}=
\#\{Q(L) \mid L\in\mathcal{AL}_{n}\}/
\#\mathcal{AL}_{n}.
\end{gather*}
These are the \emph{distinct values} $Q$ takes on \emph{alternating links}. For the next statement, the reader may want to recall \emph{Conway mutation} as, for example, in \cite[Section 2.3]{Ad-knots}. Let $K_{2}$ be Khovanov homology in characteristic $2$.

\begin{Theorem}[Exponential decay theorem]\label{T:Drop}
For any quantum invariant $Q$ that satisfies a skein relation, is multiplicity free or does not detect Conway mutation we have
\begin{gather*}
Q(n)^{\%}_{AL}\in O(\delta^{n})
\text{ for some }\delta=\delta(Q)\in(0,0.996).
\end{gather*}
This applies to $Q\in\{A2,A,B1,J,K_{2}\}$.
\end{Theorem}

\begin{proof}
Any invariant that satisfies a skein relations or is multiplicity free, this invariant cannot now detect Conway mutation. This is pointed out in \cite{We-mutation,110646}. 

Now, let $\mathcal{ALCM}_{n}$ denote $\mathcal{AL}_{n}$ modulo Conway mutation.
It is known that
\begin{gather*}
\lim_{n\to\infty}
\sqrt[n]{\#\mathcal{AL}_{n}}
=\frac{\sqrt{21001}+101}{40},
\end{gather*}
see \cite[Theorem 1]{SuTh-growth}. By \cite[Theorem 2]{St-number-polynomials}, it is further known that 
\begin{gather*}
\limsup_{n\to\infty}
\sqrt[n]{\#\mathcal{ALCM}_{n}}
\leq
\frac{\sqrt{24584873929}+109417}{43334}
<\frac{\sqrt{21001}+101}{40}-0.004.
\end{gather*}
Thus, the theorem follows for $\delta$ being $0.996$ or smaller, and we just have to comment on the various special cases.
Over $\mathbb{F}_{2}$ Khovanov homology fails to detect Conway mutation \cite{We-mutation-2}, while A satisfies a skein relation, and A2, B1, and J are multiplicity free.
\end{proof}

The reader might have noticed that the proof of \autoref{T:Drop} also works for variations of quantum invariants. Explicitly, 
the HOMFLYPT polynomial satisfies a skein relation, $k$-colored Jones polynomials, for all $k\in\N$, are multiplicity free and odd Khovanov homology does not detect Conway mutation, so they all satisfy the same exponential decay.

\begin{Remark}
It seems likely that \autoref{T:Drop} also works for generalized mutation as in \cite{AnPrRo-mutation} and other relations on the set of knots or links. Essentially, ``any relation'' imposed on the set of knots or links will make the quotient set an order of magnitude smaller, and quantum invariants will satisfy some relation.
\end{Remark}

\begin{Remark}
Khovanov homology and actually KT1 do detect some Conway mutations \cite{We-mutation}, but they still do poorly. So, one should expect that there is another form of ``mutation'' that is not detected by Khovanov homology and other categorification, and this mutation will eventually dominate the set of knots. In fact, Khovanov homology satisfies a skein relation with an error term \cite{ChGoKoLa-skein}, which is probably enough to force exponential decay. However, this would require some linearization (or similar) of the set-based proof of \autoref{T:Drop}. While we do not expect it to be difficult, we do anticipate it being somewhat tedious, which is why we have omitted it from this work.
\end{Remark}

\begin{Remark}
A stronger condition than $V\otimes V$ being multiplicity free is $V$ being strongly multiplicity free, see \cite[Theorem 3.4]{LeZh-strongly-multiplicity-free} for a classification. Surprisingly, this stronger property, which implies \autoref{T:Drop}, is also related to faithful representations of the braid group, see \cite{LaTuVa-verma-howe}. This is another example of why braids are considered easier than knots and links.
\end{Remark}

\section{Comparison -- distinguishing pairs}\label{S:Compare2}

In \autoref{S:Compare1} we asked how likely one can distinguish
a knot from all others. An alternative way to measure how good an invariant is 
would be to asked how likely one can distinguish two knots from one another.
That is, we define a function 
\begin{gather*}
Q(n)^{\%,\%}=E(\%,\%),
\end{gather*}
where $E(\%,\%)$ is the expected number of times it is necessary to randomly select $K,L\in\mathcal{K}_{n}$, for $K\neq L$, until $Q(K)=Q(L)$.

\begin{Remark}
For the following data we used JavaScript's \texttt{Math.random()} to select two distinct knots (pseudo) randomly from the list of all prime knots $\mathcal{K}_{n}$. In particular, each data point is the average over 100,000 trials and is an approximation of $Q(n)^{\%,\%}$. As noted in ECMAScript specification (2025), \texttt{Math.random()} uses an implementation-defined pseudo-random number generation algorithm and thus may differ between machines. By design, the generator cannot be seeded.
\end{Remark}

With the same notation as in \autoref{S:Compare1}, the collected data is as follows. The second plot is a magnification of the tail of the first plot, without a log scale, to accentuate potential classes of invariants.
\begin{gather*}
\begin{tikzpicture}[anchorbase]
\node at (0,0) {\includegraphics[width=0.6\textwidth]{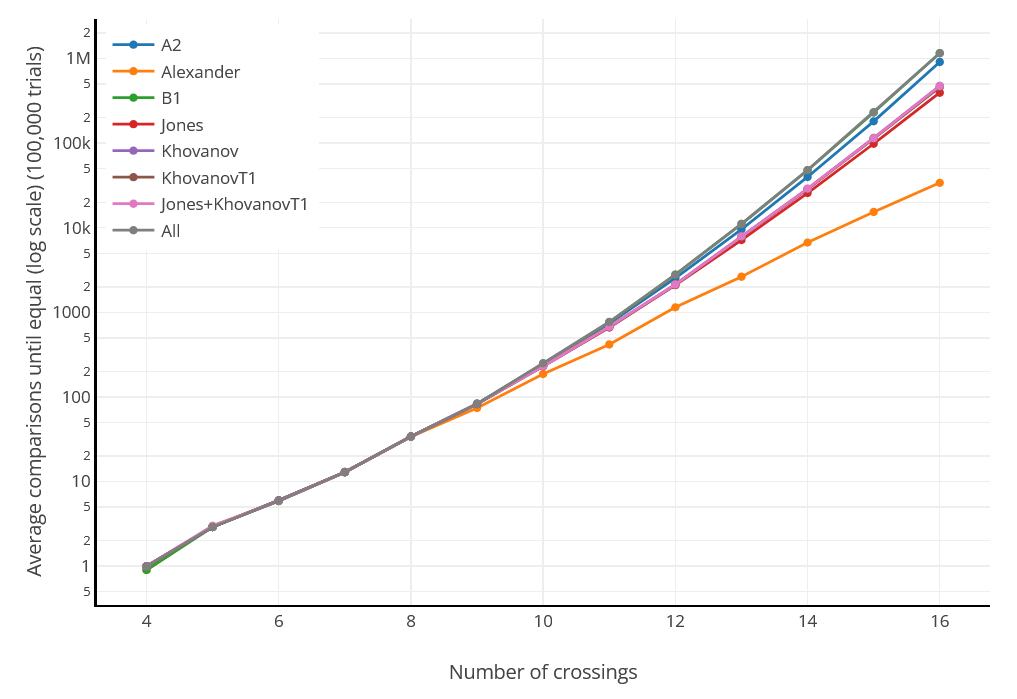}};
\end{tikzpicture}
,
\end{gather*}
\begin{gather*}
\begin{tikzpicture}[anchorbase]
\node at (0,0) {\includegraphics[width=0.6\textwidth]{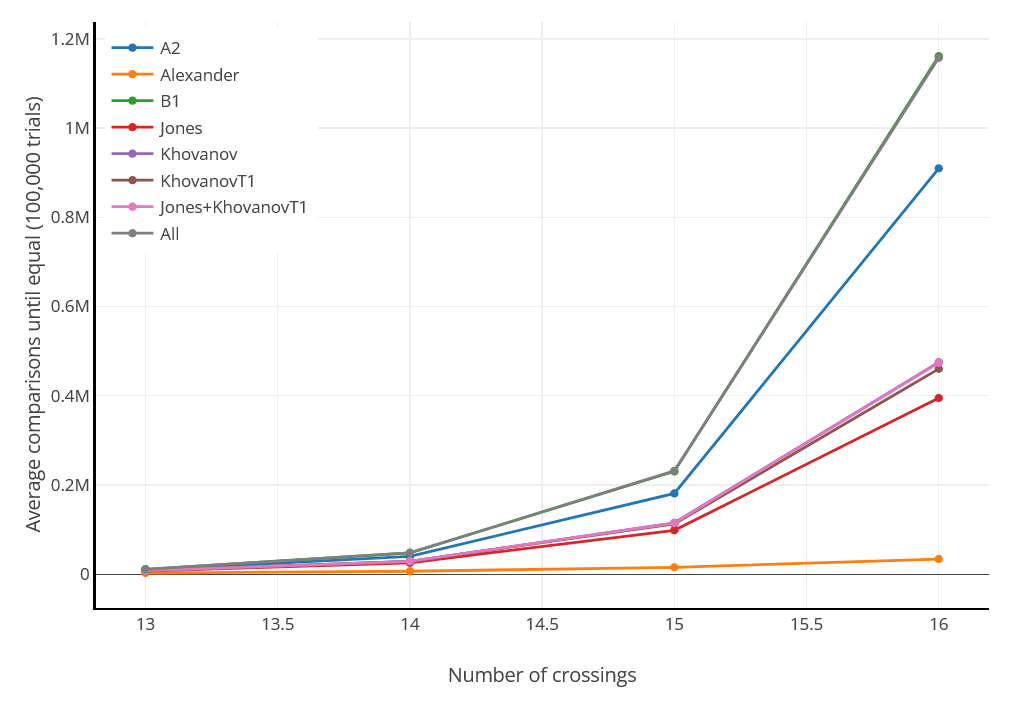}};
\end{tikzpicture}
.
\end{gather*}
We have listed the numbers in \autoref{figure:comparison2}.
Here are again a few observations:
\begin{enumerate}

\item The only change compared to \autoref{S:Compare1} is that now the probability
of detection goes to one instead of zero.

\item We stress that the preferred method again appears to be coloring.

\end{enumerate}

As before in \autoref{S:Compare1}, but slightly adjusted by changing to \emph{exponential growth}:

\begin{Conjecture}\label{C:GoesUp}
For $Q\in\{A2,A,B1,J,KT1\}$ (and therefore also for $Q=K$) we have
\begin{gather*}
Q(n)^{\%,\%}\in\Omega(\gamma^{n})
\text{ for some }\gamma=\gamma(Q)\in\R_{>1}.
\end{gather*}
Moreover, we have
\begin{gather*}
\gamma(B1)>\gamma(A2)>\gamma(K)>\gamma(A).
\end{gather*}
In other words, they all grow exponentially fast, but the exponential factor depends on the class.
\end{Conjecture}

We even speculate the following.

\begin{Speculation}\label{S:GoesUp}
For all quantum invariants $Q$ we have
\begin{gather*}
Q(n)^{\%,\%}\in\Omega(\gamma^{n})
\text{ for some }\gamma=\gamma(Q)\in\R_{>1}.
\end{gather*}
In other words, they all grow exponentially fast.
\end{Speculation}

\begin{Speculation}\label{S:ExpGrowth}
Assume \autoref{S:GoesUp} holds and fix $Q_{\mathfrak{g},V,0}$ for $\mathfrak{g}$ of classical type. Then:
\begin{gather*}
\gamma(color)>\gamma(rank)>\gamma(cat).
\end{gather*}
That is, coloring is better than rank increase is better than categorification.
\end{Speculation}

\begin{Remark}
As in \autoref{S:Compare1}, after cleaning for computational complexity, one might argue that the ordering in \autoref{S:ExpGrowth} needs to be adjusted.
\end{Remark}

In order to prove \autoref{C:GoesUp}, or any of the others, it is likely beneficial to study random knots (see {\eg} \cite{EvZo-random-knots} for a summary) and their values on quantum invariants. For example, the knots generated by the random model in
\cite{EvZoHaLiNo-random-knots} seem to be prime with high probability, and it would be interesting to consider whether one can compute, for example, a Jones polynomial in this model.

\section{Comparison -- Big data approach}\label{S:Compare4}

The following is inspired by \cite{LaTuVa-big-data}. 
Fix a quantum invariant $Q=Q_{\mathfrak{g},V,\epsilon}$. 
Examples of questions one could try to address are the following:
\begin{enumerate}

\item What is the asymptotic behavior, for $n\to\infty$, of
\begin{gather*}
\op{ev}_{n}=\optt{max}\{{\textstyle\sum} \mathtt{abs} \big(Q(K)\big) \mid K \in \mathcal{K}_{n}\}?
\end{gather*}
This corresponds to taking the sum of the coefficients. However, merely doing so would be dequantization, which is trivial, and so instead we take the sum of the absolute values $\optt{abs}$ of the coefficients.

\item \label{Q:GrowthB} Similarly, what is the asymptotic behavior, for $\sank\to\infty$, of
\begin{gather*}
\op{coeff}_{n} = \optt{max}\big\{ \optt{max}\big(\optt{abs}\big(Q(K)\big)\big) \mid K\in\mathcal{K}_{n}\big\}?
\end{gather*}
That is, we ask how fast the (absolute values of the) coefficients grow.

\item We can ask the same questions as in the previous two points, but for the \emph{average} instead of the maximum. That is,
\begin{gather*}
\op{ev}^{av}_{\sank} =
{\textstyle\sum}\big\{ {\textstyle\sum} \optt{abs}\big(Q(K)\big) \mid K\in\mathcal{K}_{n}\}/\#\{Q(K) \mid K\in\mathcal{K}_{n}\big\}
\end{gather*}
where we take the average over sums, and
\begin{gather*}
\op{coeff}^{av}_{\sank} =
\optt{max} \big\{ {\textstyle\sum} \optt{abs}\big(Q(K)\big) / \optt{span}\big(Q(K)\big) \mid K\in\mathcal{K}_{n}\} \big\}.
\end{gather*}
where we take the maximum over ``average coefficients''. Note that for $\op{coeff}^{av}_{n}$ we do not count padded zeros.

\item One could also ask about maximal and average \emph{span} (maximal minus minimal degree), which we denote as $\op{span}^{av}_{n}$ and $\op{span}_{n}$.
\item Note that it does not make sense to ask for a span of a two-variable polynomial, so invariant K is excluded from data that involves spans or average coefficients.

\end{enumerate}

Now we show the data. First, $\op{coeff}_n$, the maximal coefficient.

\begin{gather*}
\begin{tikzpicture}[anchorbase]
\node at (0,0) {\includegraphics[width=0.6\textwidth]{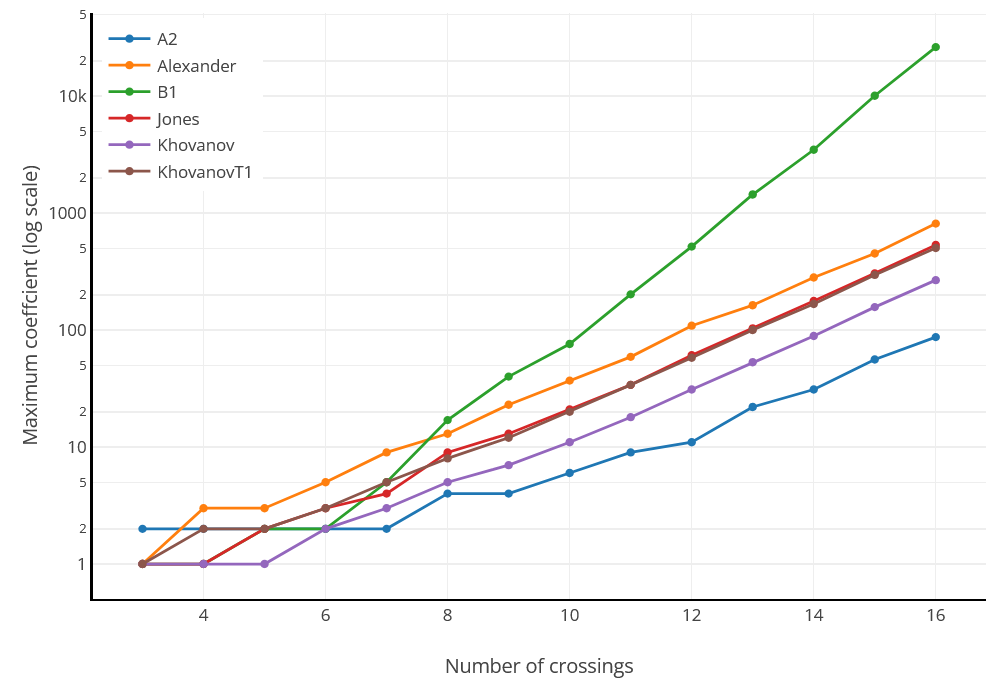}};
\end{tikzpicture}
,
\end{gather*}
\begin{gather*}
\begin{tikzpicture}[anchorbase]
\node at (0,0) {\includegraphics[width=0.6\textwidth]{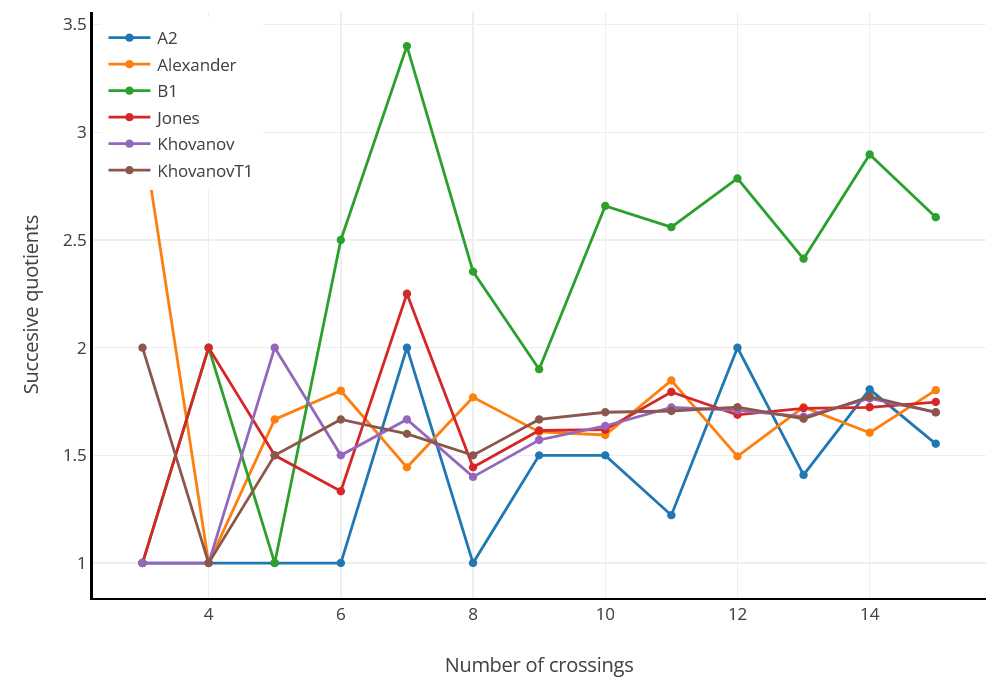}};
\end{tikzpicture}
,
\end{gather*}
and the precise figures can be found in \autoref{figure:3}.
Next, $\op{ev}_n$, the maximal coefficient sum.
\begin{gather*}
\begin{tikzpicture}[anchorbase]
\node at (0,0) {\includegraphics[width=0.6\textwidth]{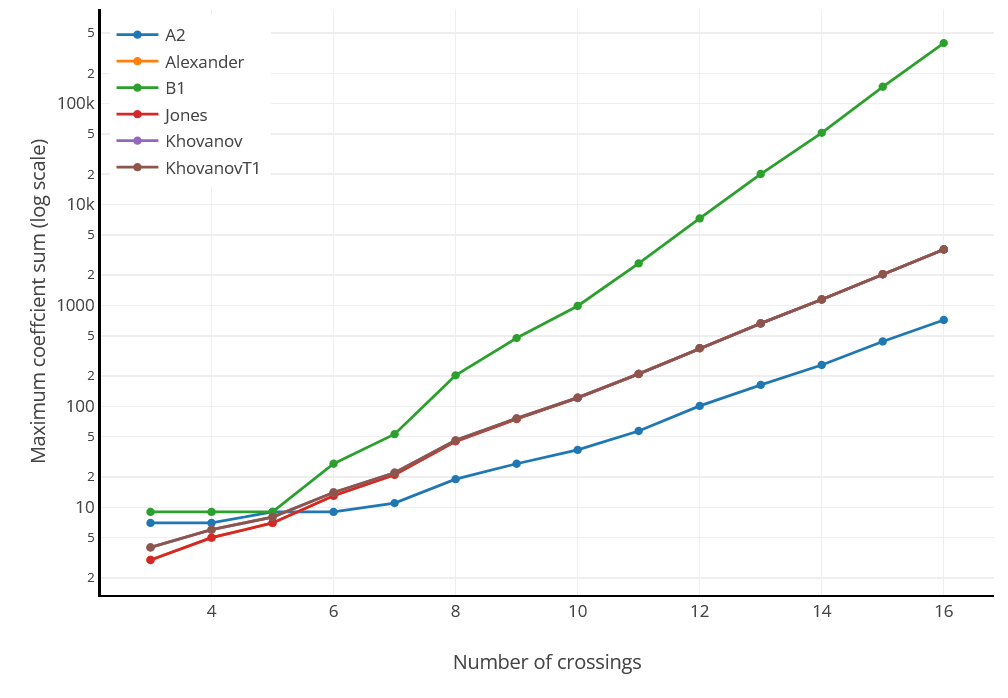}};
\end{tikzpicture}
,
\end{gather*}
\begin{gather*}
\begin{tikzpicture}[anchorbase]
\node at (0,0) {\includegraphics[width=0.6\textwidth]{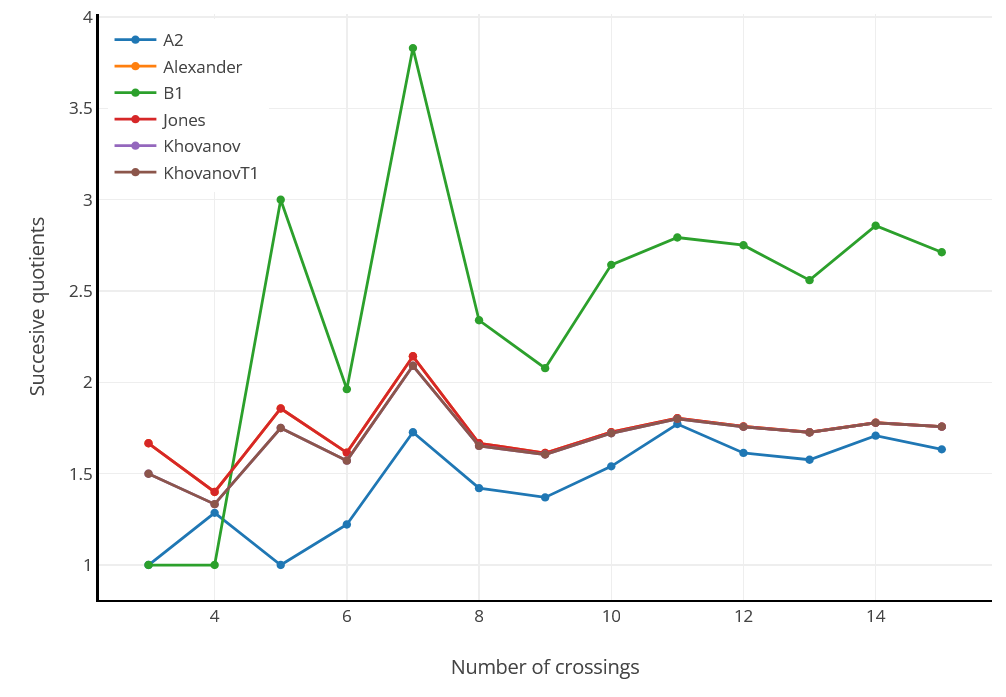}};
\end{tikzpicture}
,
\end{gather*}
with data in \autoref{figure:4}.
Then, taking averages, $\op{ev}^{av}_{\sank}$ and $\op{coeff}^{av}_{\sank}$.
\begin{gather*}
\begin{tikzpicture}[anchorbase]
\node at (0,0) {\includegraphics[width=0.6\textwidth]{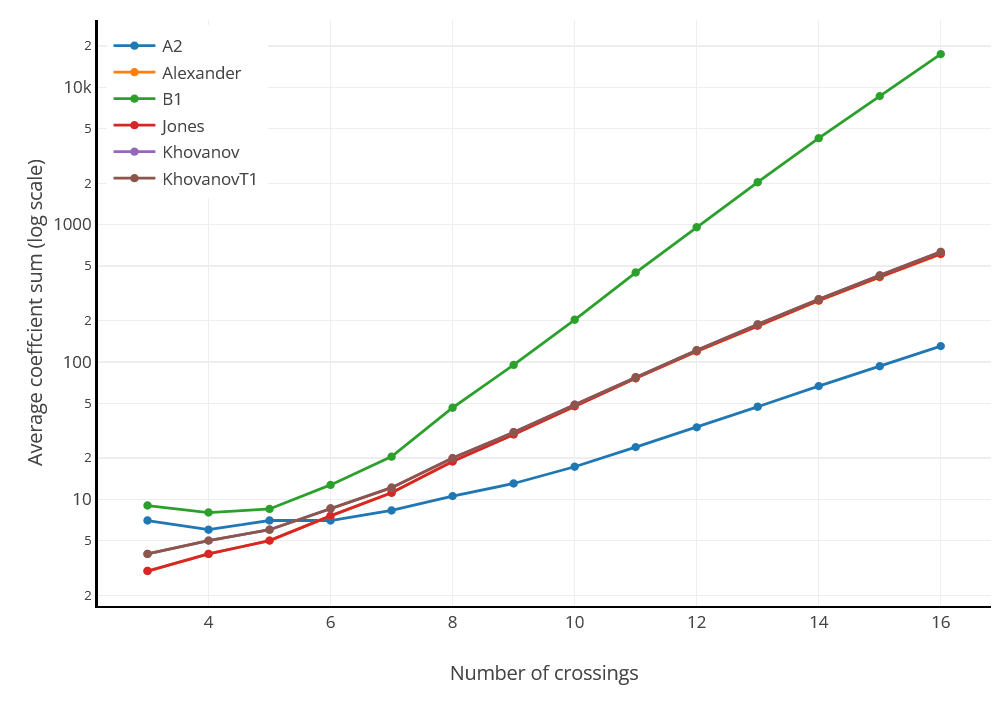}};
\end{tikzpicture}
,
\end{gather*}
\begin{gather*}
\begin{tikzpicture}[anchorbase]
\node at (0,0) {\includegraphics[width=0.6\textwidth]{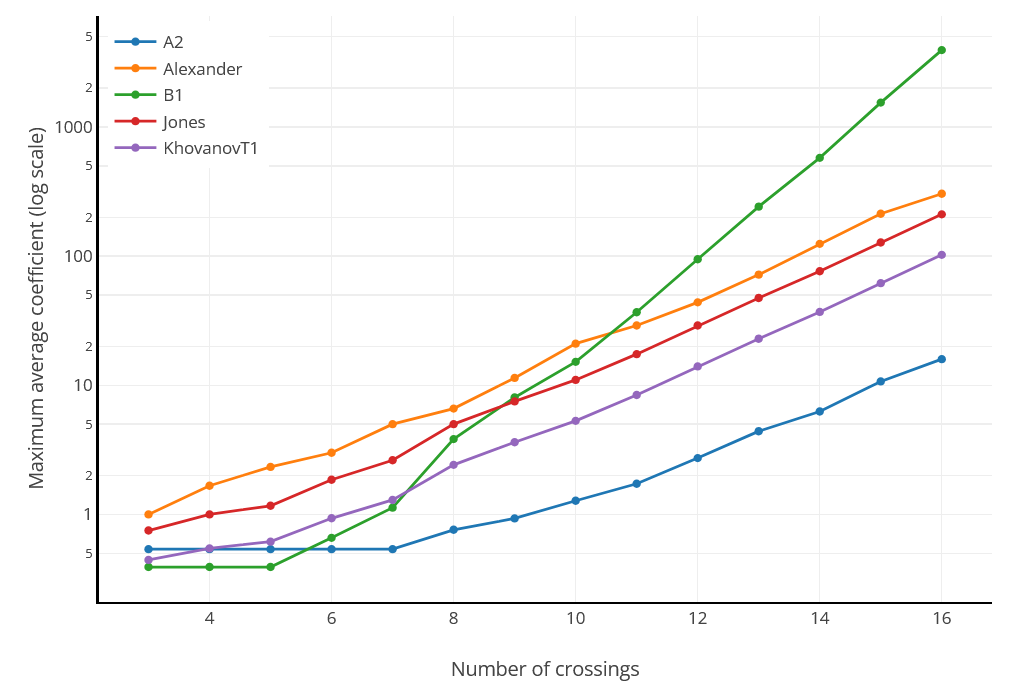}};
\end{tikzpicture}
,
\end{gather*}
with data in \autoref{figure:5} and \autoref{figure:6}.
Then we look at spans, $\op{span}_{n}$ and $\op{span}^{av}_{n}$.
\begin{gather*}
\begin{tikzpicture}[anchorbase]
\node at (0,0) {\includegraphics[width=0.6\textwidth]{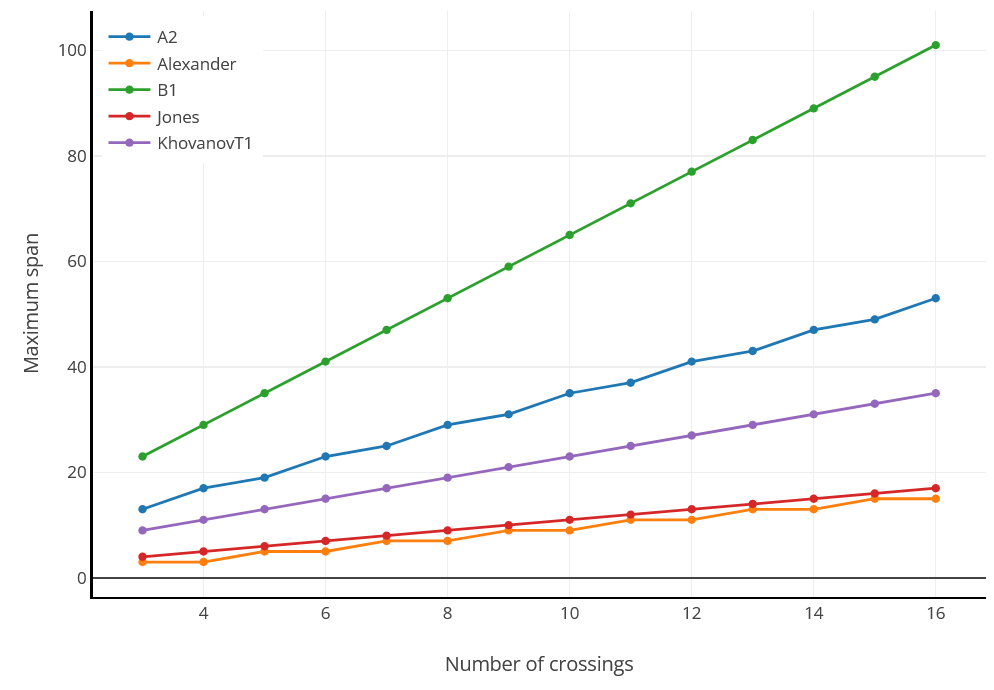}};
\end{tikzpicture}
,
\end{gather*}
\begin{gather*}
\begin{tikzpicture}[anchorbase]
\node at (0,0) {\includegraphics[width=0.6\textwidth]{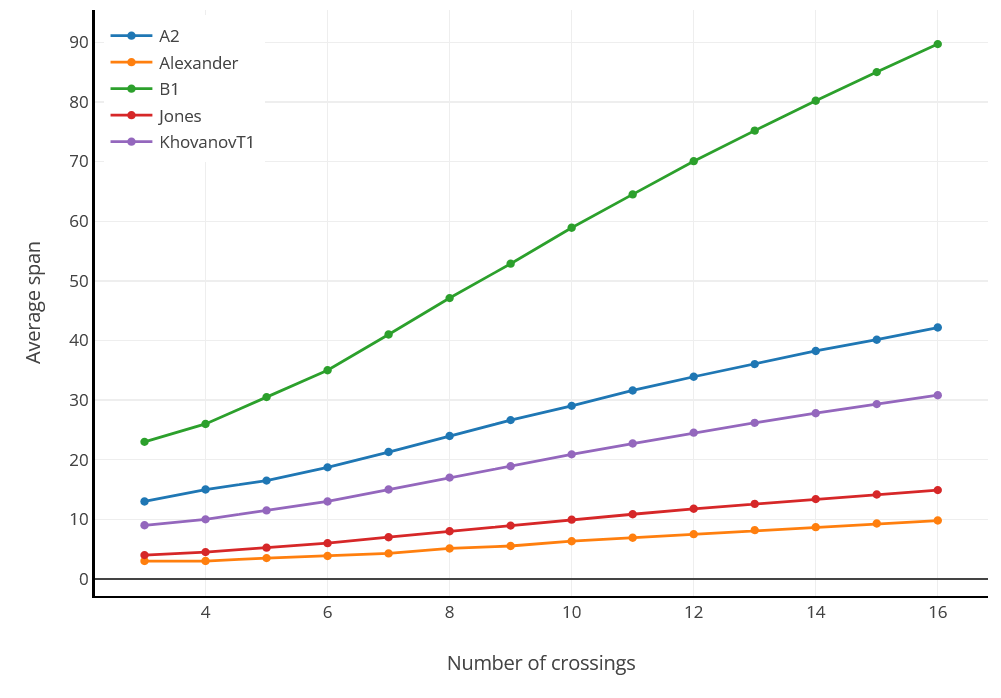}};
\end{tikzpicture}
,
\end{gather*}
where we refer to \autoref{figure:7} and \autoref{figure:8}.

From the data we conjecture:

\begin{Conjecture}[Exponential growth]\label{C:Growth}
For $Q\in\{A2,A,B1,Kh\}$ (and therefore also for $Q=J$ or $Q=KT1$),
we have
\begin{gather*}
\mathrm{coeff}_{n}\in\Omega(\gamma^{n})\text{ for some }\gamma\in\R_{>1}.
\end{gather*}
Since $\mathrm{ev}_{n}\geq\mathrm{coeff}_{n}$, the same holds for
$\mathrm{ev}_{n}$. We also have
\begin{gather*}
\gamma(B1)>\gamma(Q)\text{ for $Q\in\{A2,A,J,KT1\}$}
\end{gather*}
(and therefore also for $Q=K$).
\end{Conjecture}

We now compare pictures of roots, following the ideas in \cite{BaChDe-roots,WuWa-jones-roots,LaTuVa-big-data} and the references to various blogs in \cite{BaChDe-roots}.

The next few pictures show
the multiset of roots of the set of $Q\in\{A2,A,B1,J,KT1\}$, with the pictures zoomed in on the right. In formulas, we plot
\begin{gather*}
\Big(\optt{roots}(p) \mid p\in\{Q(K) \mid K\in\mathcal{K}_{16}\} \Big)
\end{gather*}
where high brightness indicates the high density.
Each picture has the same scaling in the axes and are centered on the origin. A faint green circle is drawn to indicate the unit circle. The left pictures are zoomed such that all roots are displayed with nothing more, and the right pictures are zoomed to the same scale as each other.
\begin{gather*}
A1:
\begin{tikzpicture}[anchorbase]
\node at (0,0) {\includegraphics[width=0.45\textwidth]{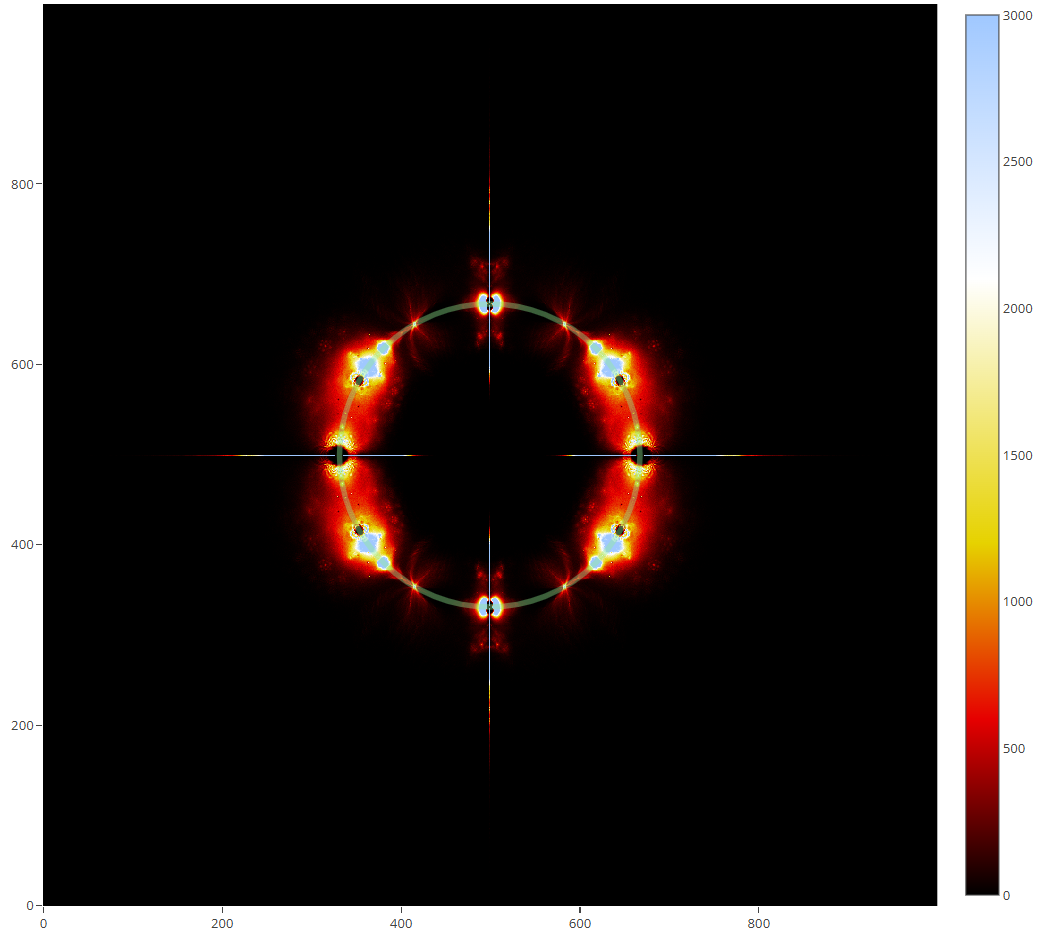}};
\end{tikzpicture}
,
\begin{tikzpicture}[anchorbase]
\node at (0,0) {\includegraphics[width=0.45\textwidth]{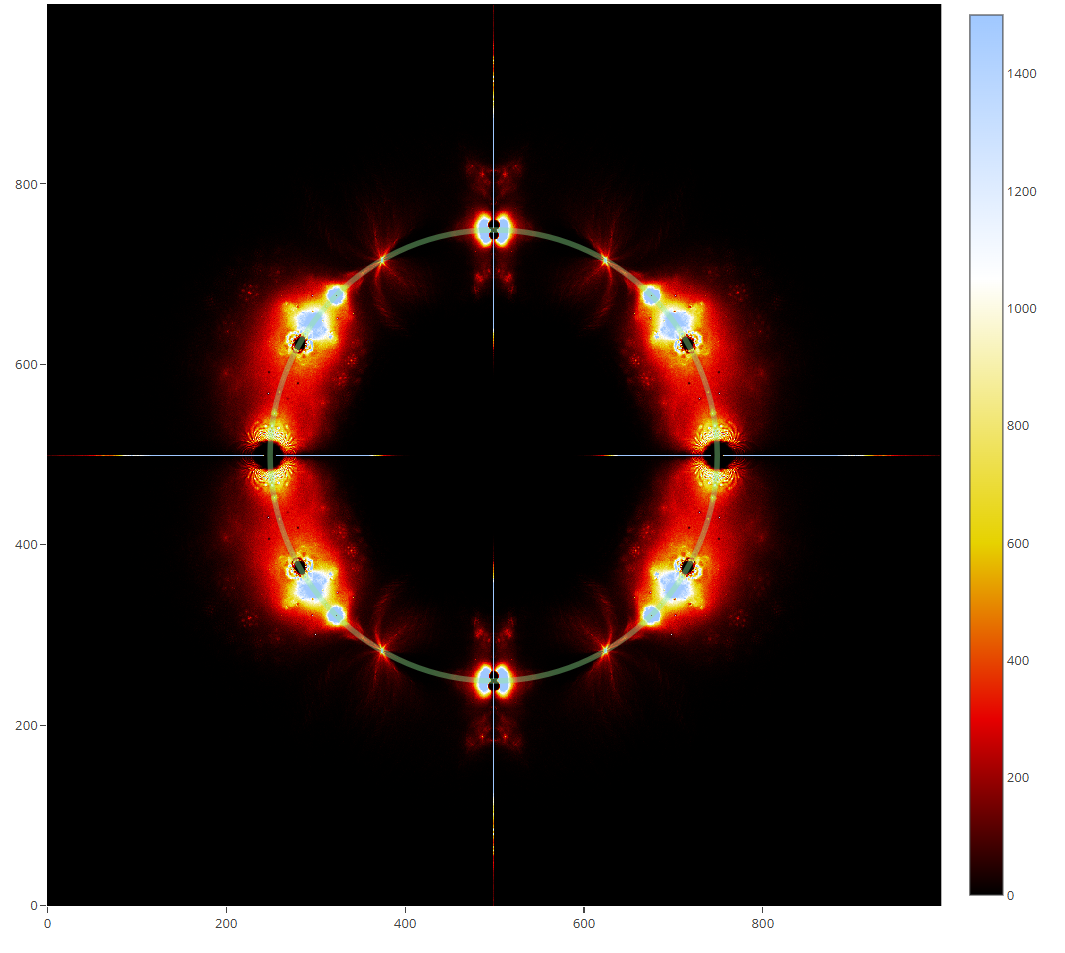}};
\end{tikzpicture}
,
\end{gather*}
\begin{gather*}
A:
\begin{tikzpicture}[anchorbase]
\node at (0,0) {\includegraphics[width=0.45\textwidth]{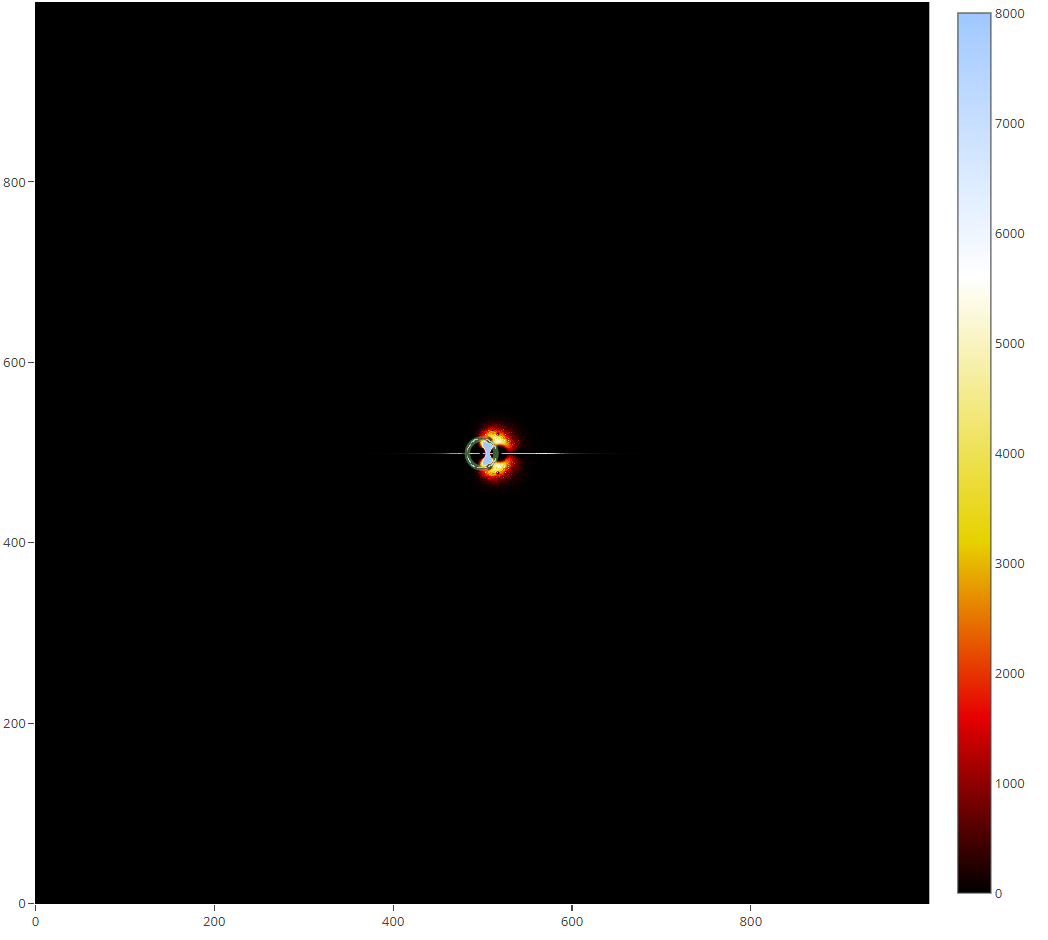}};
\end{tikzpicture}
,
\begin{tikzpicture}[anchorbase]
\node at (0,0) {\includegraphics[width=0.45\textwidth]{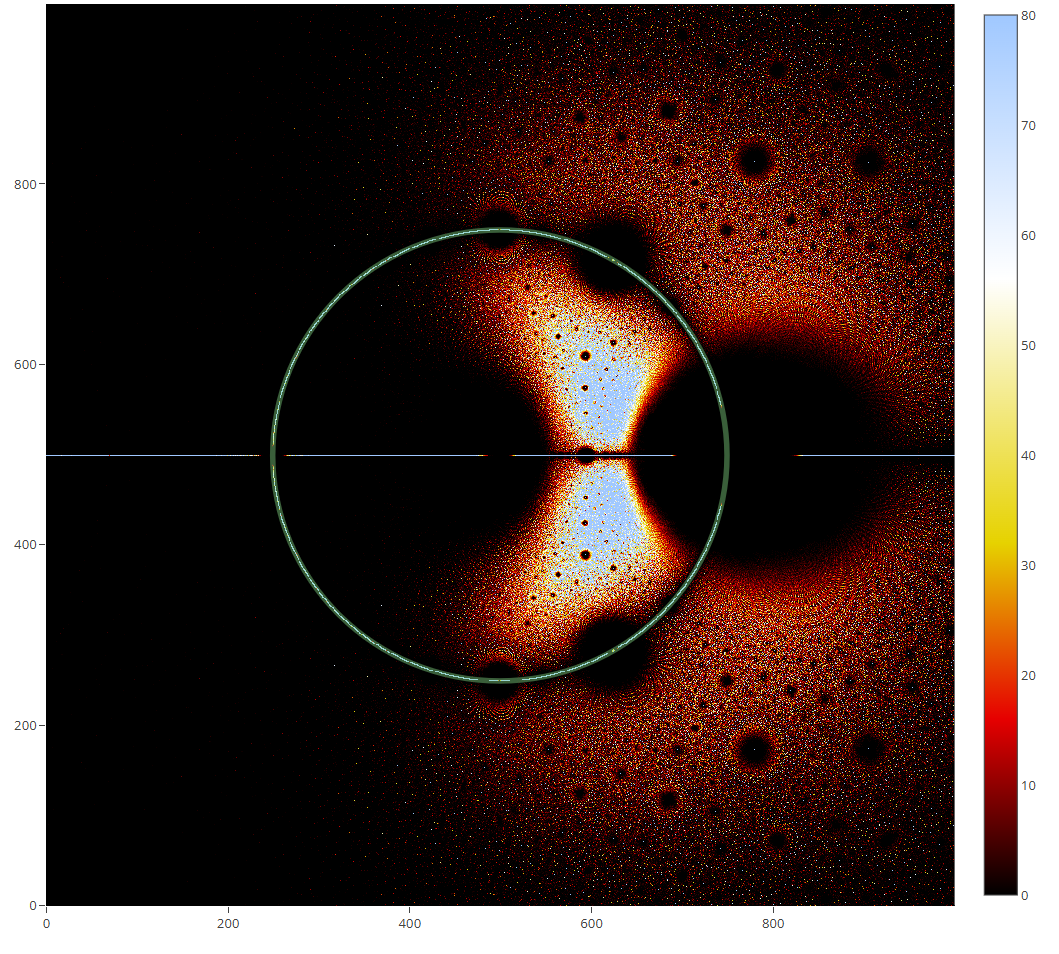}};
\end{tikzpicture}
,
\end{gather*}
\begin{gather*}
B1:
\begin{tikzpicture}[anchorbase]
\node at (0,0) {\includegraphics[width=0.45\textwidth]{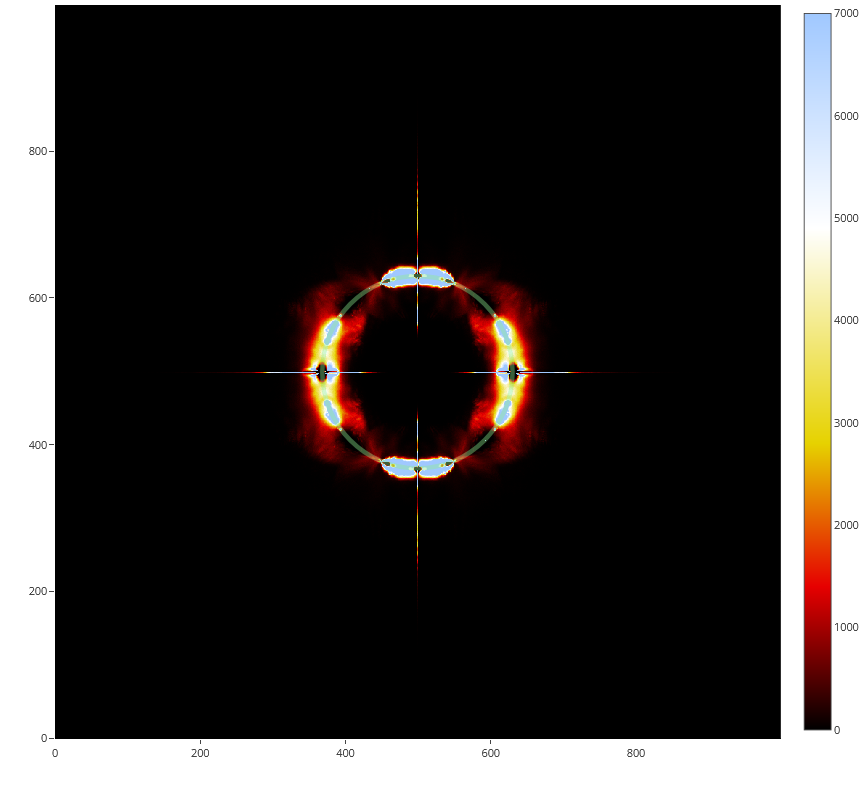}};
\end{tikzpicture}
,
\begin{tikzpicture}[anchorbase]
\node at (0,0) {\includegraphics[width=0.45\textwidth]{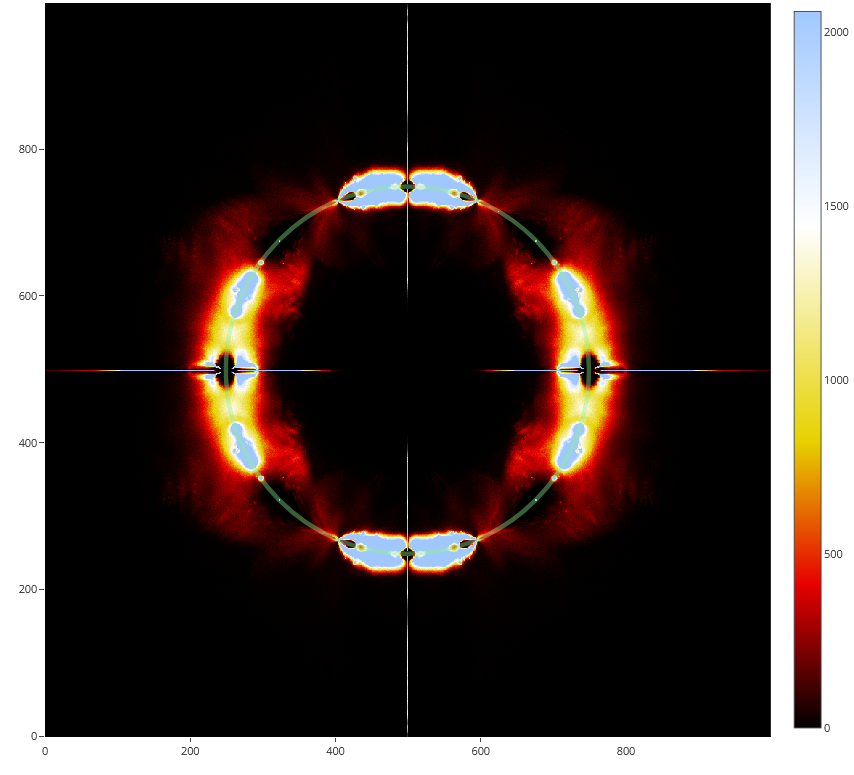}};
\end{tikzpicture}
,
\end{gather*}
\begin{gather*}
J:
\begin{tikzpicture}[anchorbase]
\node at (0,0) {\includegraphics[width=0.45\textwidth]{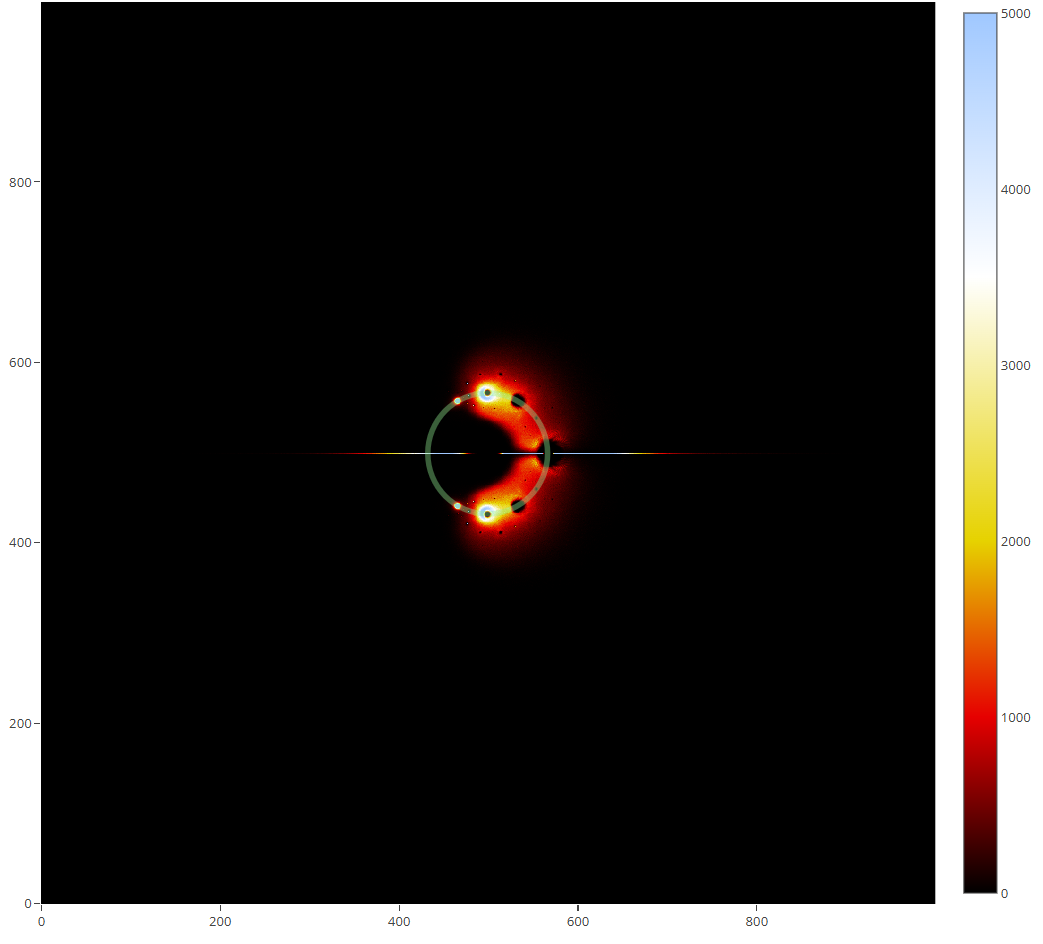}};
\end{tikzpicture}
,
\begin{tikzpicture}[anchorbase]
\node at (0,0) {\includegraphics[width=0.45\textwidth]{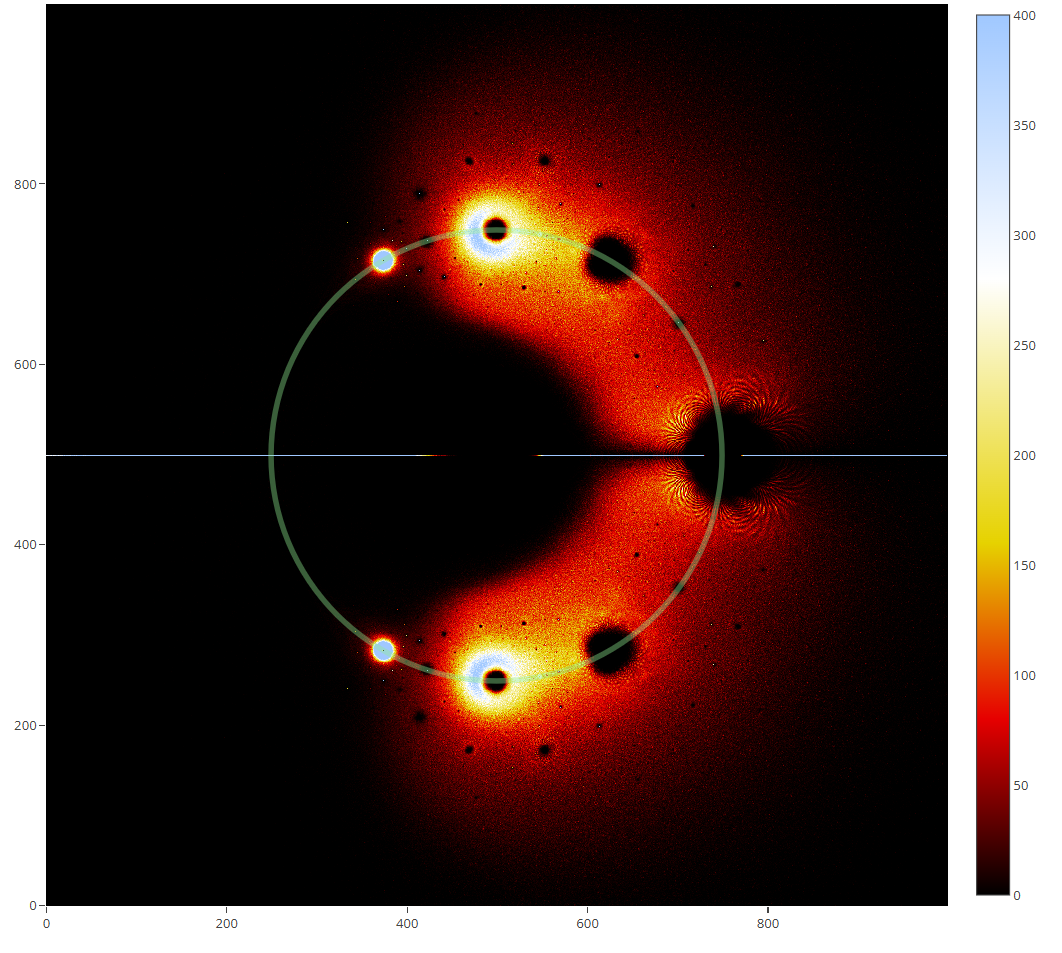}};
\end{tikzpicture}
,
\end{gather*}
\begin{gather*}
KT1:
\begin{tikzpicture}[anchorbase]
\node at (0,0) {\includegraphics[width=0.45\textwidth]{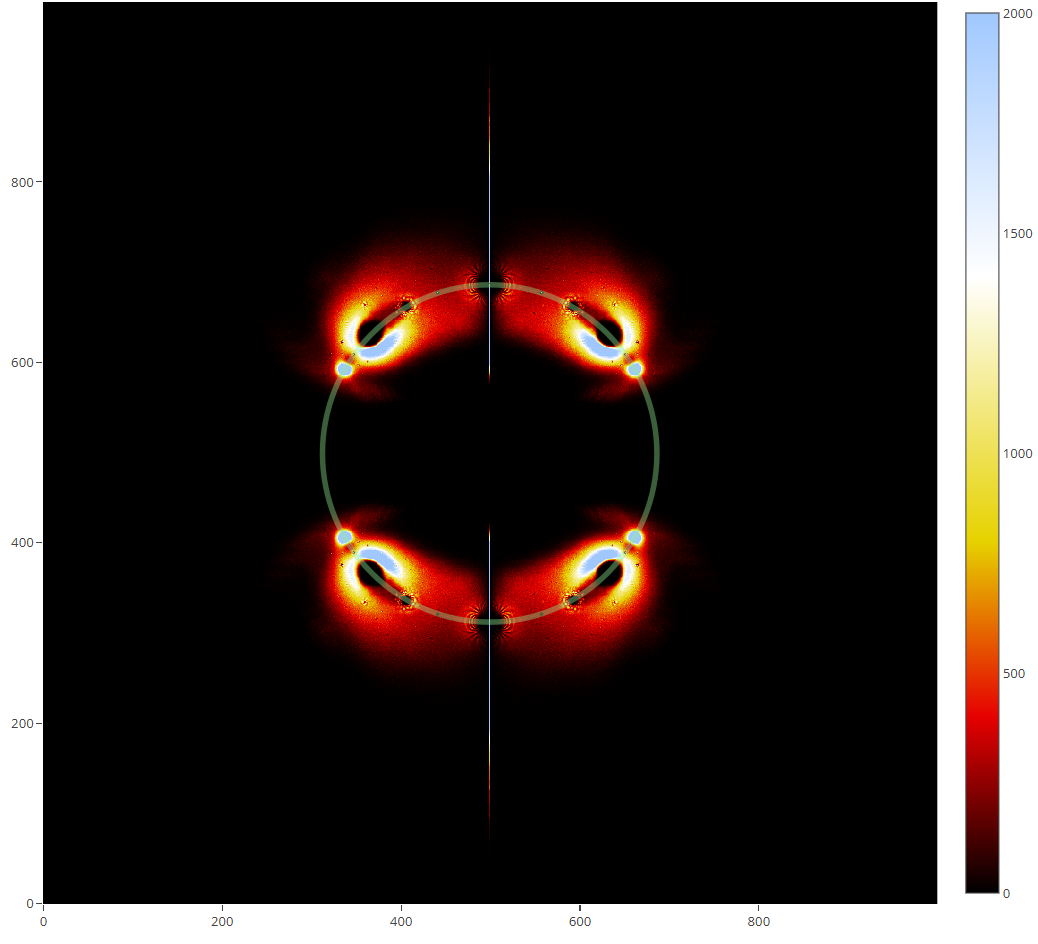}};
\end{tikzpicture}
,
\begin{tikzpicture}[anchorbase]
\node at (0,0) {\includegraphics[width=0.45\textwidth]{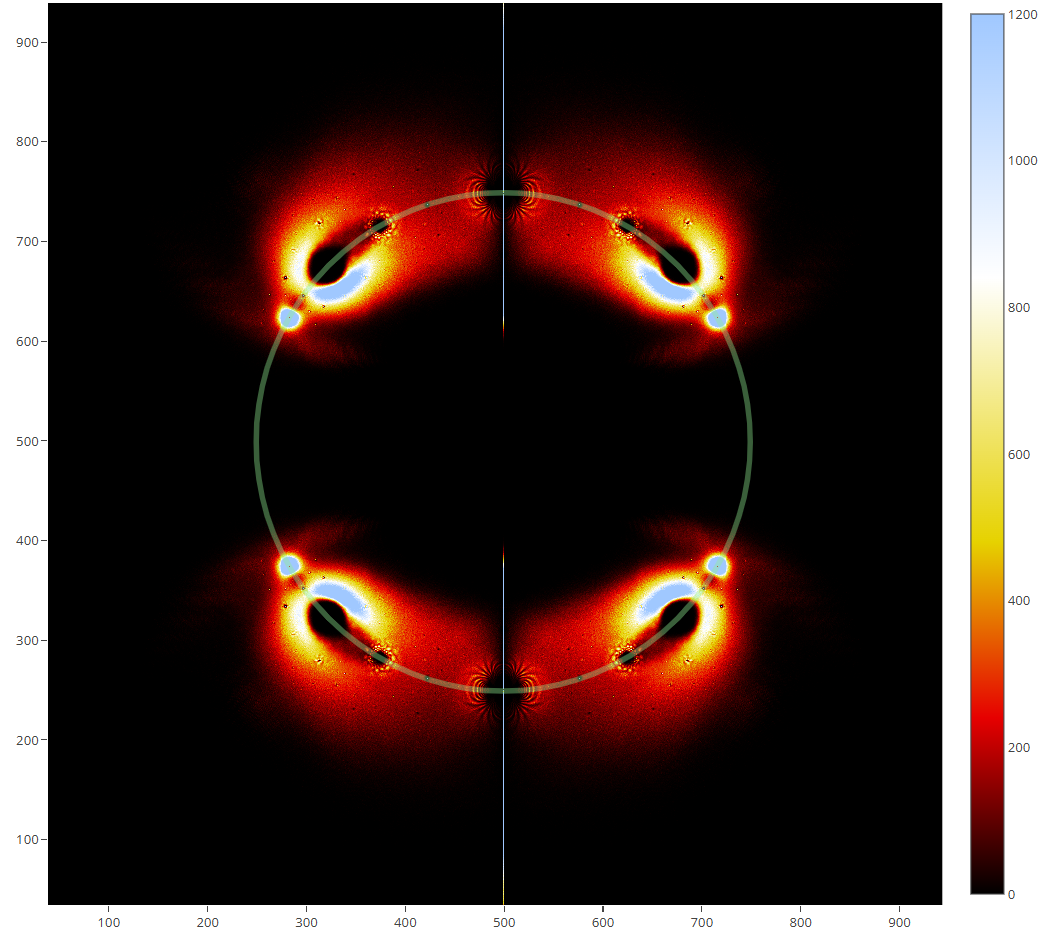}};
\end{tikzpicture}
.
\end{gather*}
Here is the associated data.
\begin{gather*}
\begin{tikzpicture}[anchorbase]
\node at (0,0) {\includegraphics[width=0.6\textwidth]{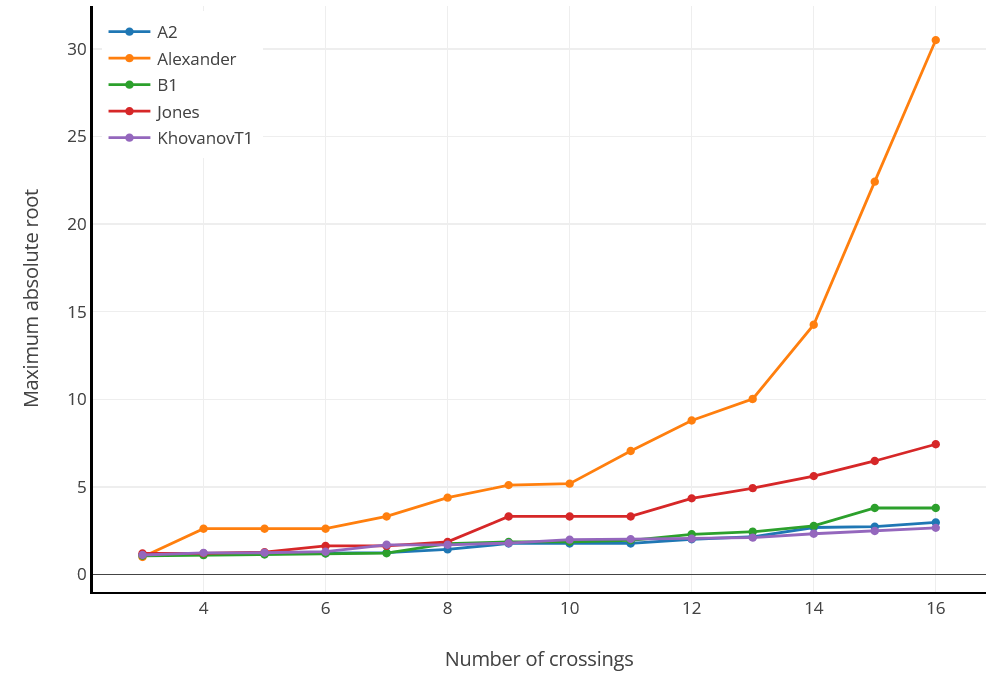}};
\end{tikzpicture}
,
\end{gather*}
see \autoref{figure:9} for details. Continued:
\begin{gather*}
\begin{tikzpicture}[anchorbase]
\node at (0,0) {\includegraphics[width=0.6\textwidth]{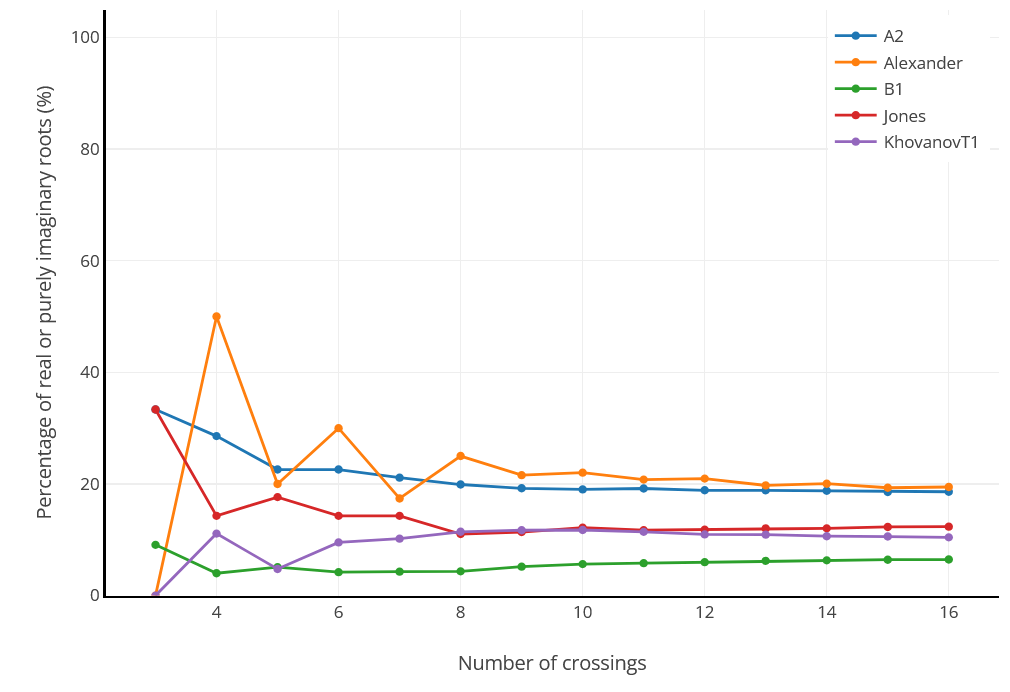}};
\end{tikzpicture}
,
\end{gather*}
\begin{gather*}
\begin{tikzpicture}[anchorbase]
\node at (0,0) {\includegraphics[width=0.6\textwidth]{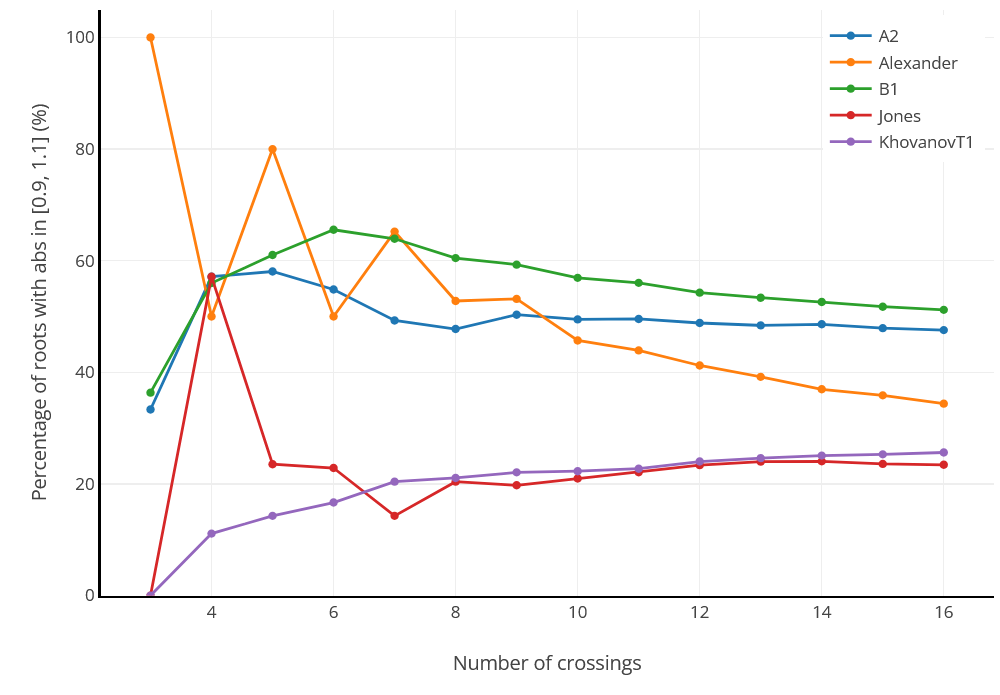}};
\end{tikzpicture}
,
\end{gather*}
which are \autoref{figure:10} and \autoref{figure:11}.

Before we analyze the pictures, we need a bit of terminology.

\begin{Notation}
A root is \emph{purely real or imaginary}, or simply pure, if it is on either of the coordinate axes.

The so-called \emph{Perron--Frobenius (PF) root} $\lambda$ is the pure maximal root when taking absolute values; this might not exist, but it almost always exists in our setting. The maximal one, when varying over all polynomials, can be read-off from the left pictures: the edges of the bounding squares in these pictures are of length twice $|\max\{\lambda|\lambda\text{ is a PF root of }Q(K)\text{ for }K\in\mathcal{K}_{16}\}|$.
\end{Notation}

We observe the following.
\begin{enumerate}

\item For A2, B1 and KT1 the roots are concentrated around the unit circle. The PF root is fairly small.

\item For A and J the roots are quite spread in the first and fourth quadrants. The PF root is much larger.

\end{enumerate}
For comparison, let us recall a few facts about the distribution of
roots of random polynomials. Notable are:
\begin{enumerate}[label=(\roman*)]

\item There is an expected tail of pure roots, with one large PF root $\mathrm{PFmax}(\text{random})$. In fact, for random polynomials, papers such as \cite{Ge-spectral,Ba-circular-law} implies that the PF root is expected to be very large. That is, letting $M$ be the bound of the coefficients, then, almost surely $\mathrm{PFmax}(\text{random})\geq\sqrt{M}$.

\item By \cite{Ka-roots-random}, the number of expected real roots of a polynomial of degree $k$ is $2 \ln (k)/\pi$. To simplify our calculation, assume that the average degree of our sample is $e^{2}$ for the usual $e\approx2.71...$\hspace{2px}. Then the expected percentage of real roots is $4\pi/e^{2}\approx 17.23\%$.

\item The clustering of roots around the unit circle is expected for random polynomials; see, for example, \cite{ShVa-roots-random}. Fairly explicit formulas for the distribution are known, for example, see
\cite{MeBeFoMaAa}, but the only thing we notice here is that, in the limit, almost all roots will have an absolute value in $[1-\epsilon,1+\epsilon]$ for all $\epsilon\in\R_{>0}$.

\item We do not know a general statement about holes in the plots of random polynomials, but see \cite{BaChDe-roots} for similar patterns (this should be true for other integer-valued polynomials as well). The case of coefficients in $\{-1,0,1\}$ is addressed, for example, in \cite{CaKoWa}.

\end{enumerate}

\section{Comparison -- ballmapper}\label{S:Compare3}

Mapper algorithms are fundamental tools in TDA and EDA, first introduced in \cite{mapper}, and are used for exploring and visualizing data. These algorithms integrate techniques such as dimensionality reduction, clustering, and graph construction to transform data into a graph. For our purposes, we will use a slight modification known as ball mapper from \cite{Dl-ballmapper}, which provides more aesthetically pleasing visualizations. Instead of a formal definition, here is a short summary:
\begin{enumerate}[label=(\roman*)]

\item Having a point cloud in $\R^{N}$ and a fixed $\epsilon\in\R_{>0}$, the ball mapper algorithm creates 
a graph $G=G(\epsilon)$ by choosing certain anchor points, and then using balls of radius $\epsilon$ around these points.

\item The vertices of $G$ are obtained by collapsing all points within one ball into a vertex. The size of the vertex corresponds to the number of points collapsed to that vertex (large \ = \ many points).

\item The edges of $G$ come from the intersections of the ball, as in this ChatGPT generated picture:
\begin{gather*}
\begin{tikzpicture}[anchorbase]
\foreach \i in {0,1,...,99} {
\fill[black] ({cos(\i*3.6)},{sin(\i*3.6)}) circle (0.02);
}
\foreach \i in {0,1,2,3,4} {
\coordinate (C\i) at ({cos(\i*72)},{sin(\i*72)});
}
\foreach \i in {0,1,2,3,4} {
\fill[green,opacity=0.3] (C\i) circle (0.8);
\draw[thick] (C\i) circle (0.8);
}
\foreach \i [count=\j from 1] in {0,1,2,3,4} {
\node[font=\large] at ({cos(\i*72)*1.6},{sin(\i*72)*1.6}) {\j};
}
\end{tikzpicture}
\Rightarrow
\begin{tikzpicture}[anchorbase]
\foreach \i in {0,1,2,3,4} {
\fill[green] (C\i) circle (0.1); 
}
\draw[thick]
(C0) -- (C1) -- (C2) -- (C3) -- (C4) -- cycle;
\foreach \i [count=\j from 1] in {0,1,2,3,4} {
\node[font=\small] at ({cos(\i*72)*1.3},{sin(\i*72)*1.3}) {\j};
}
\end{tikzpicture}
.
\end{gather*}

\end{enumerate}
We give only a sample of what can be found in \cite{TuZh-quantum-big-data-code}; in particular, the graphs are much more impressive in the interactive plot that can be found in \cite{TuZh-quantum-big-data-code}. On that page, we also explain how these were created using blueprint files provided by the webpage of the
Dioscuri Centre in Topological Data Analysis.

All plots below are for $n=15$, and the right picture is a zoomed in version.
\begin{gather*}
\text{A2}\colon
\begin{tikzpicture}[anchorbase]
\node at (0,0) {\includegraphics[height=5.2cm]{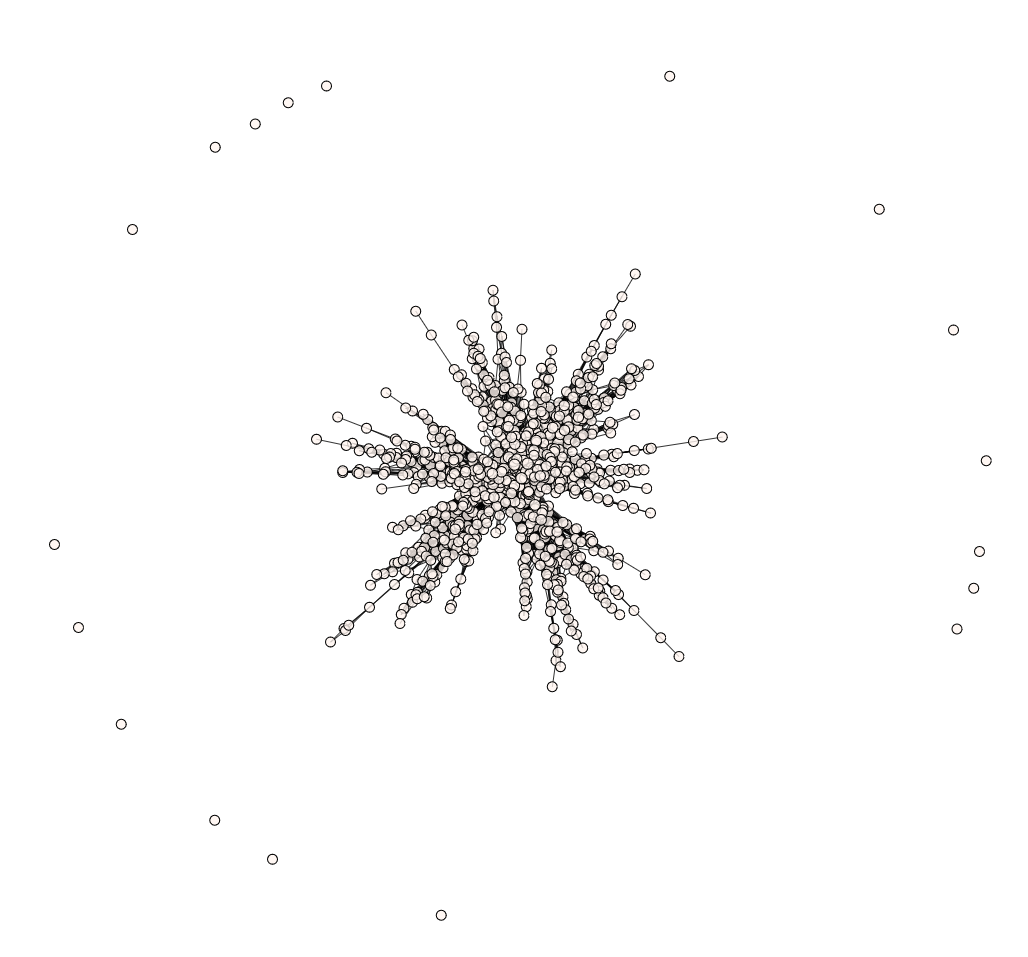}};
\end{tikzpicture}
,
\begin{tikzpicture}[anchorbase]
\node at (0,0) {\includegraphics[height=5.2cm]{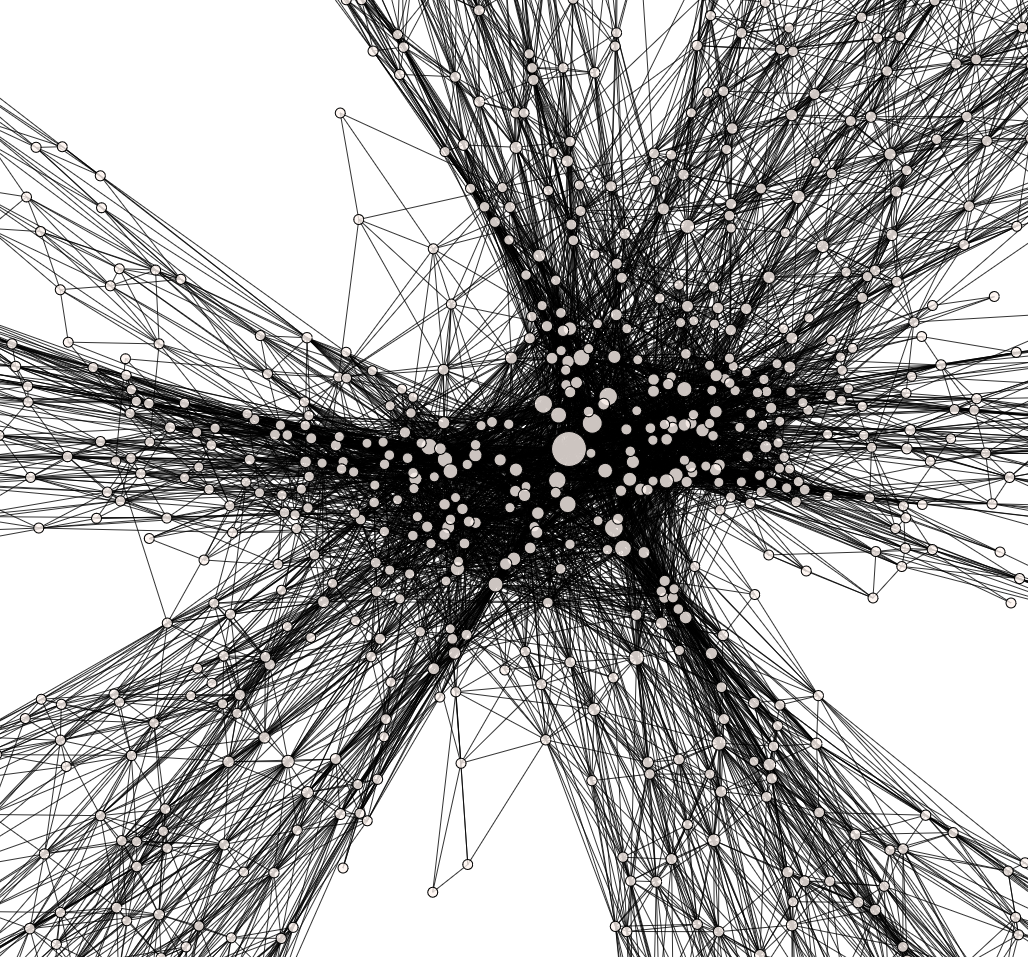}};
\end{tikzpicture}
,
\end{gather*}
\begin{gather*}
\text{A}\colon
\begin{tikzpicture}[anchorbase]
\node at (0,0) {\includegraphics[height=5.2cm]{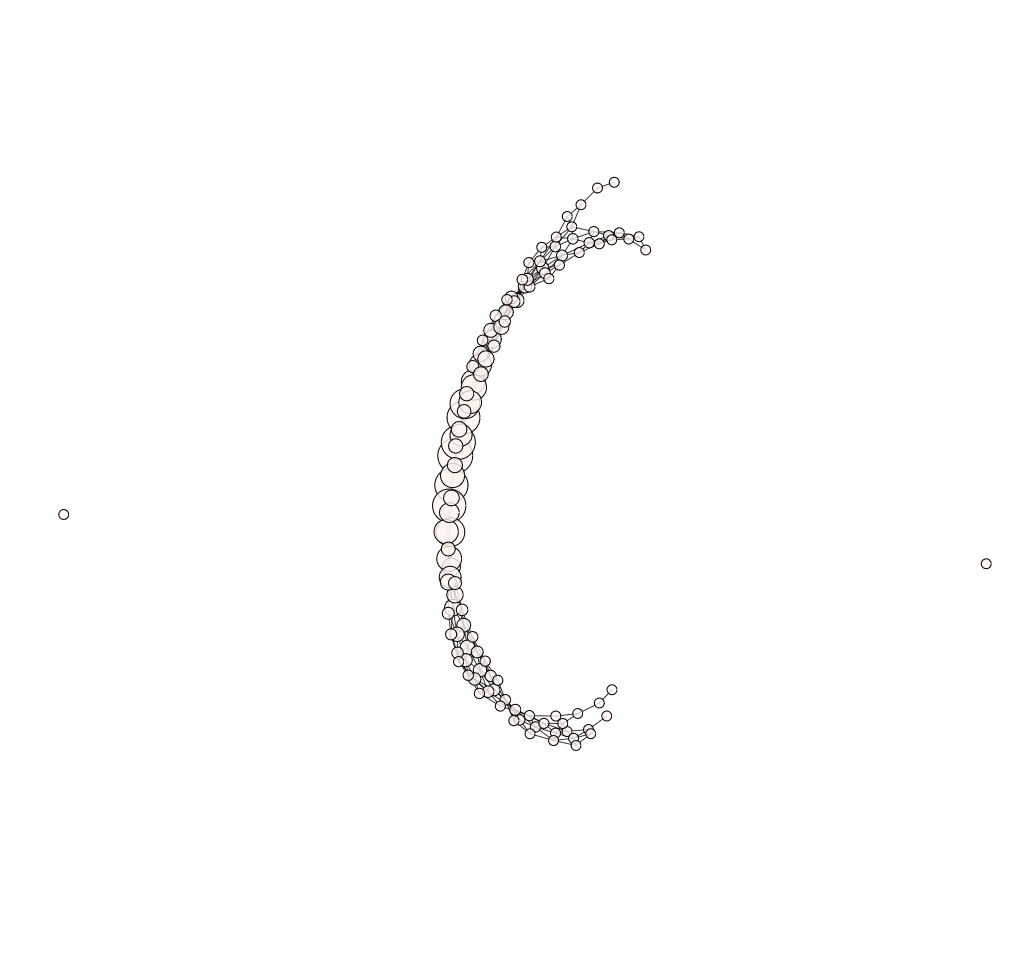}};
\end{tikzpicture}
,
\begin{tikzpicture}[anchorbase]
\node at (0,0) {\includegraphics[height=5.2cm]{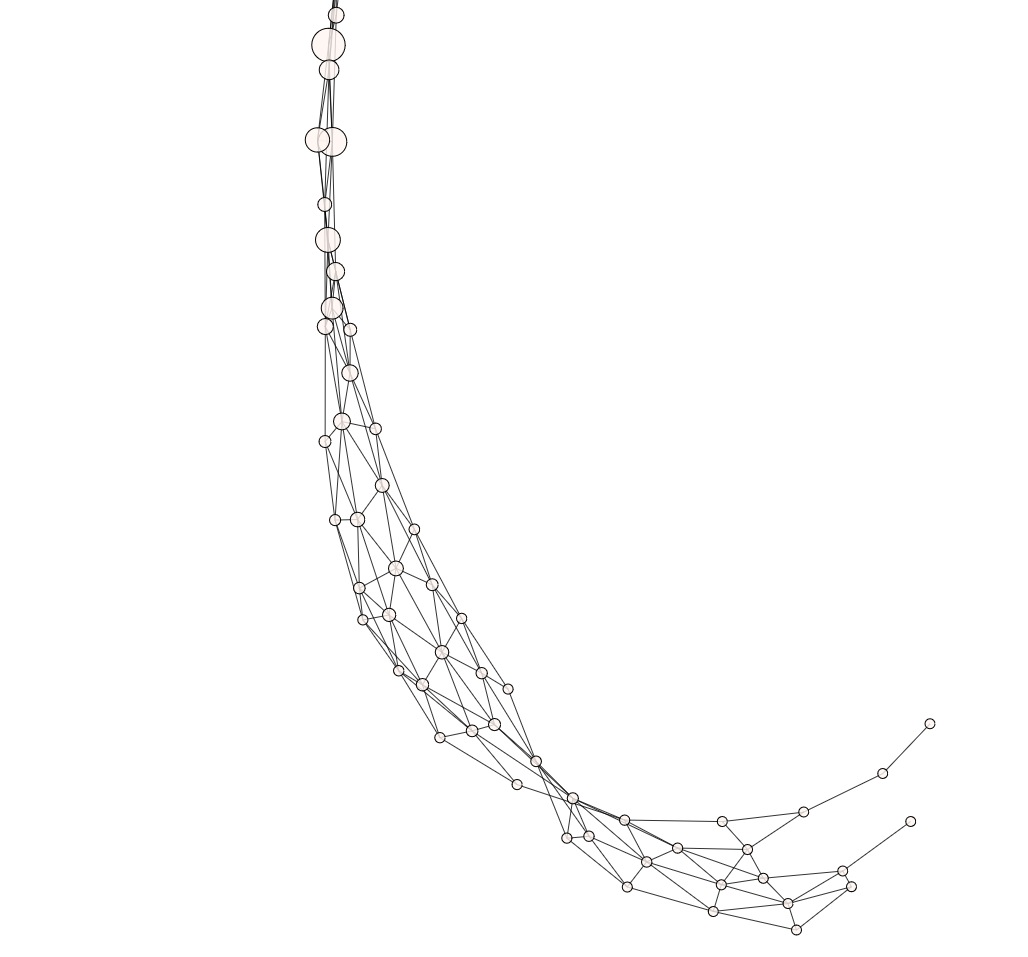}};
\end{tikzpicture}
,
\end{gather*}
\begin{gather*}
\text{B1}\colon
\begin{tikzpicture}[anchorbase]
\node at (0,0) {\includegraphics[height=5.2cm]{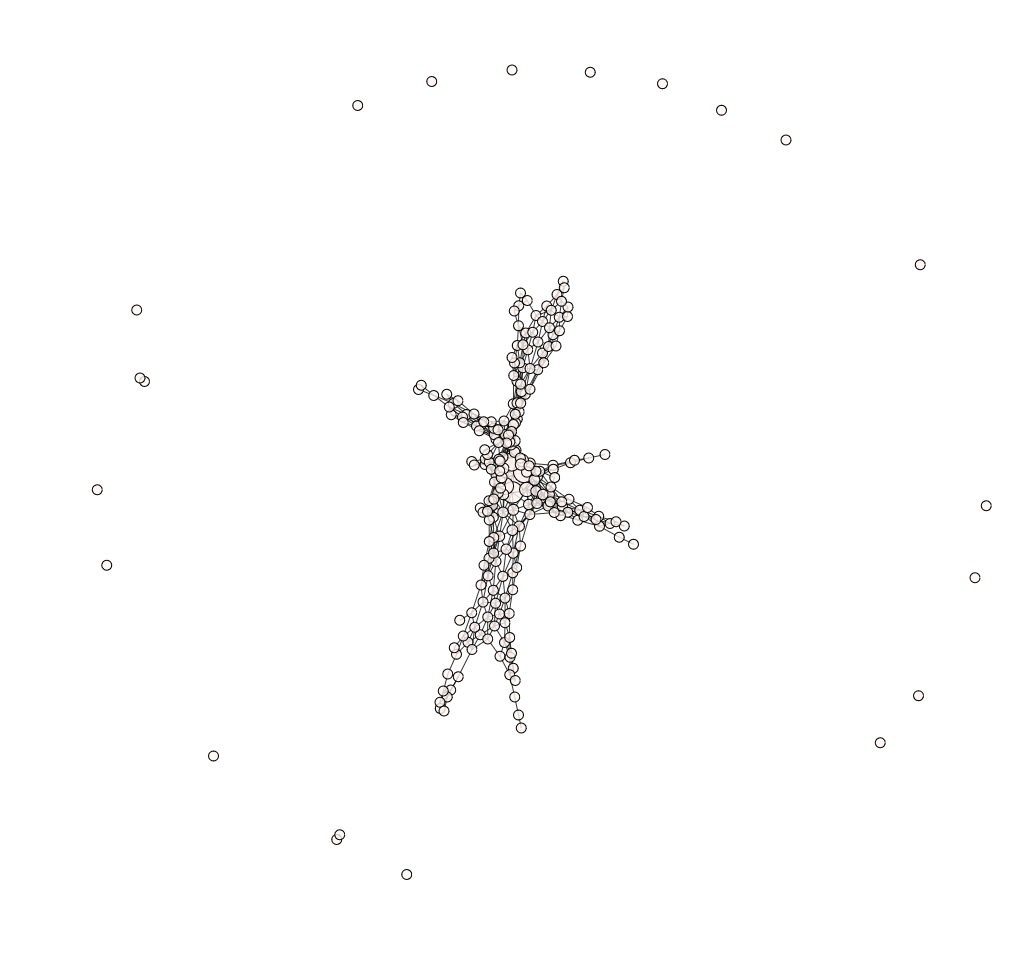}};
\end{tikzpicture}
,
\begin{tikzpicture}[anchorbase]
\node at (0,0) {\includegraphics[height=5.2cm]{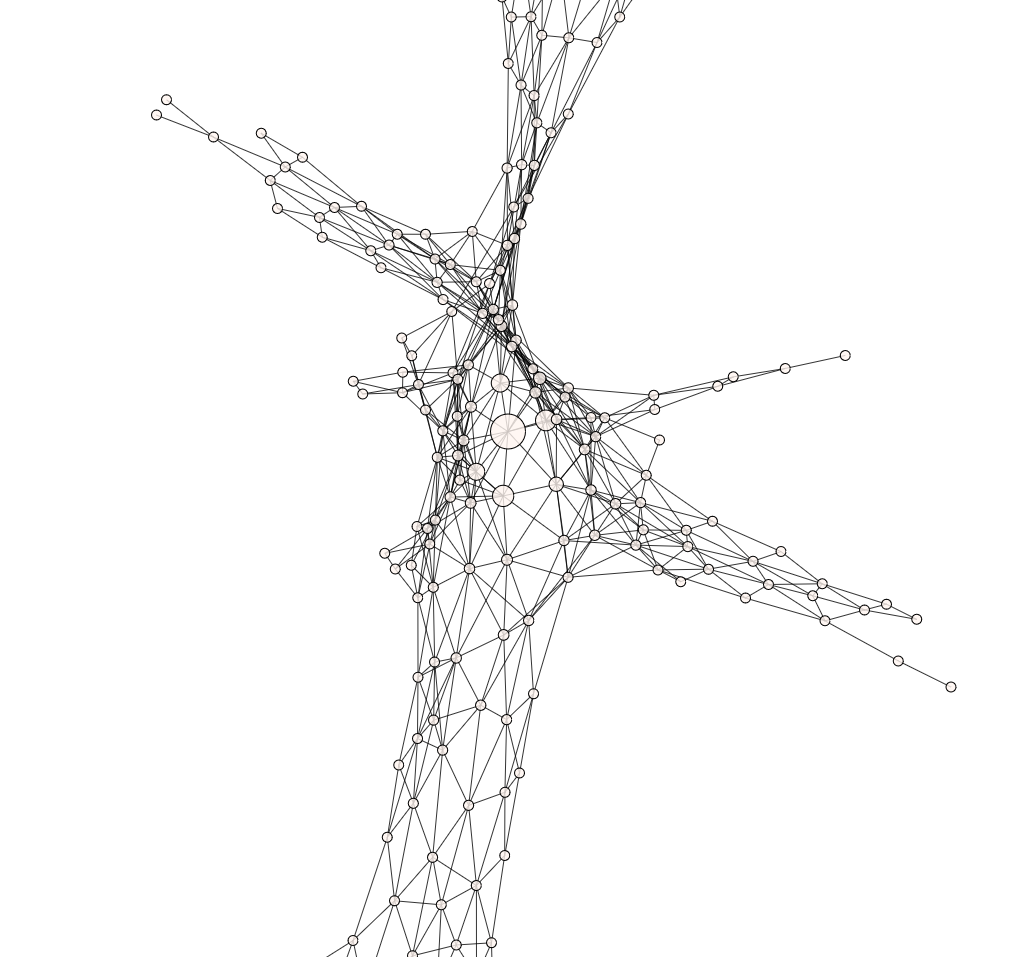}};
\end{tikzpicture}
,
\end{gather*}
\begin{gather*}
\text{J}\colon
\begin{tikzpicture}[anchorbase]
\node at (0,0) {\includegraphics[height=5.2cm]{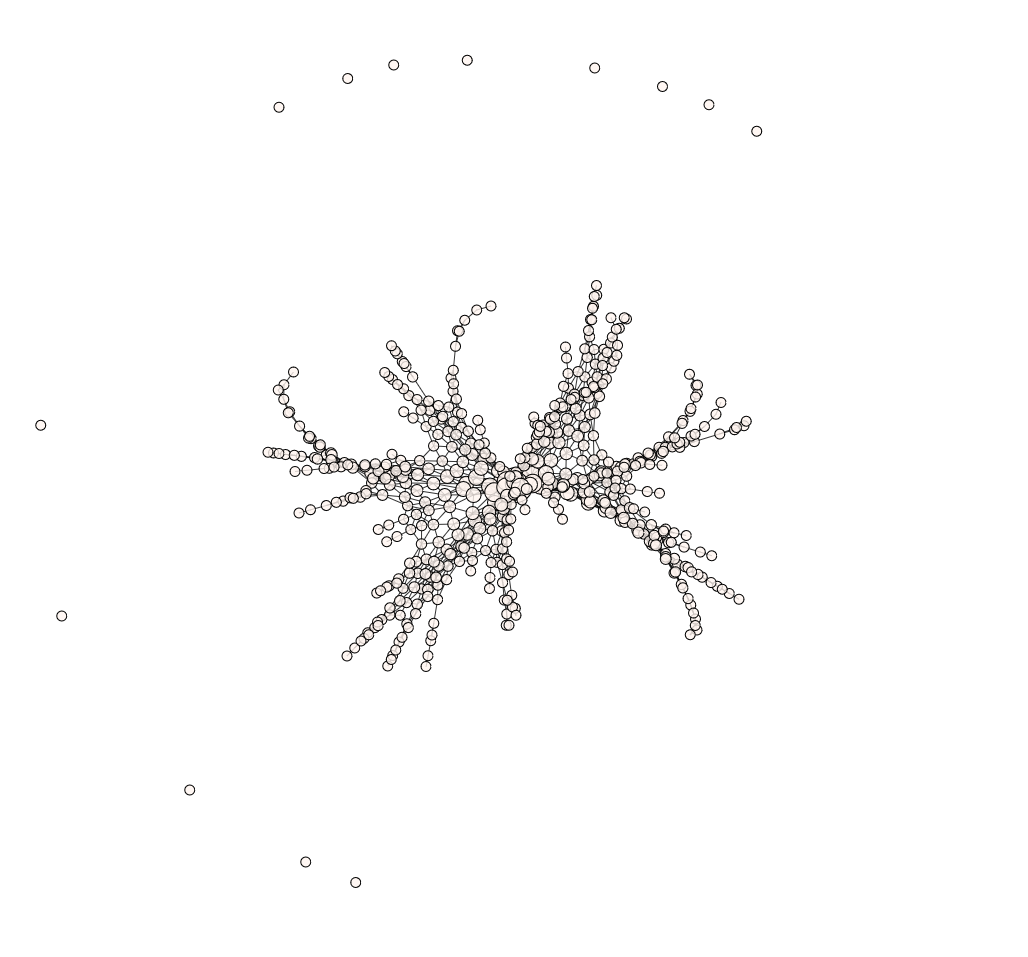}};
\end{tikzpicture}
,
\begin{tikzpicture}[anchorbase]
\node at (0,0) {\includegraphics[height=5.2cm]{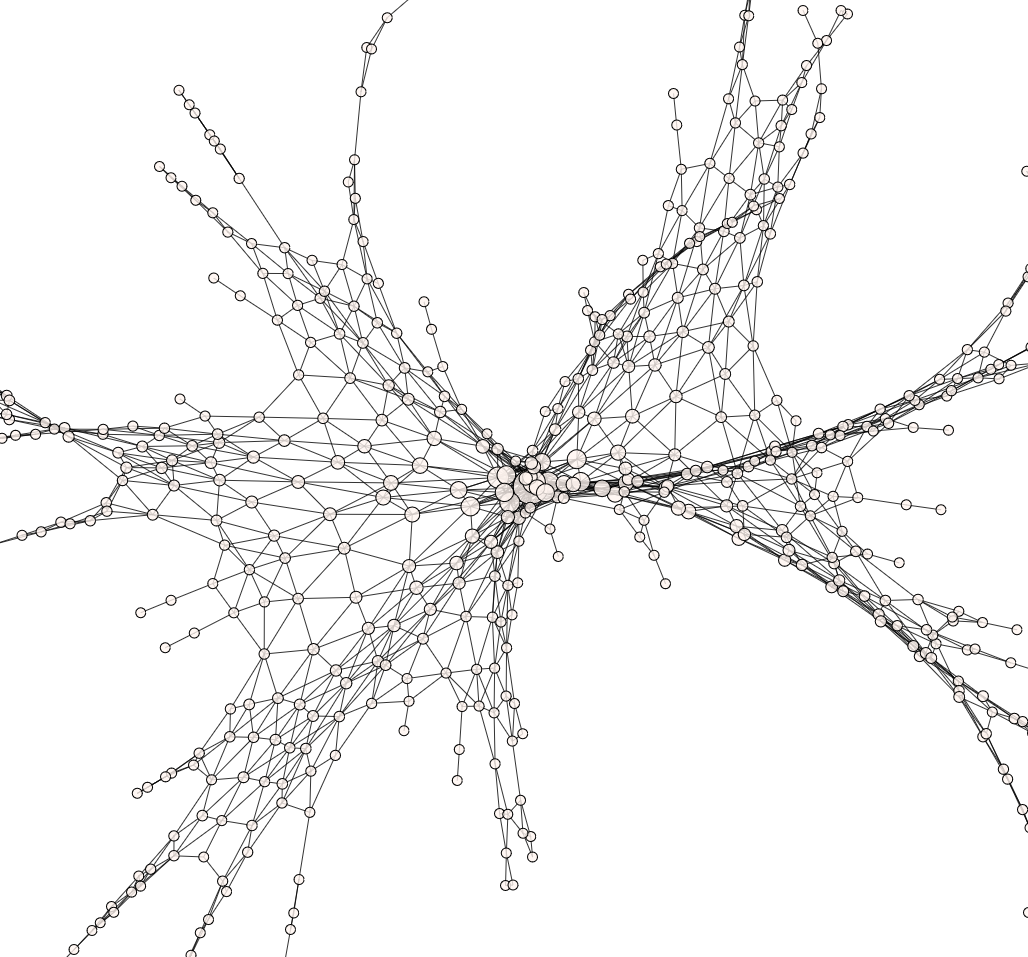}};
\end{tikzpicture}
,
\end{gather*}
\begin{gather*}
\text{K}\colon
\begin{tikzpicture}[anchorbase]
\node at (0,0) {\includegraphics[height=5.2cm]{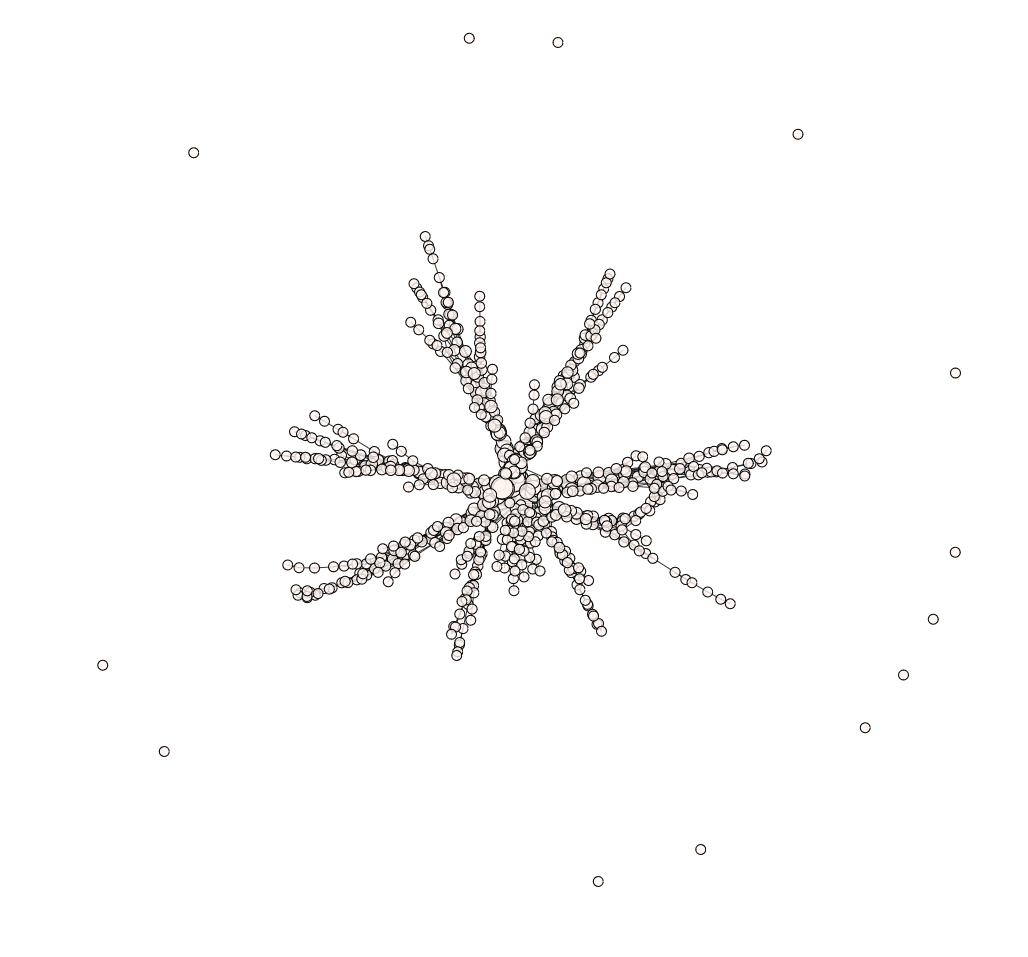}};
\end{tikzpicture}
,
\begin{tikzpicture}[anchorbase]
\node at (0,0) {\includegraphics[height=5.2cm]{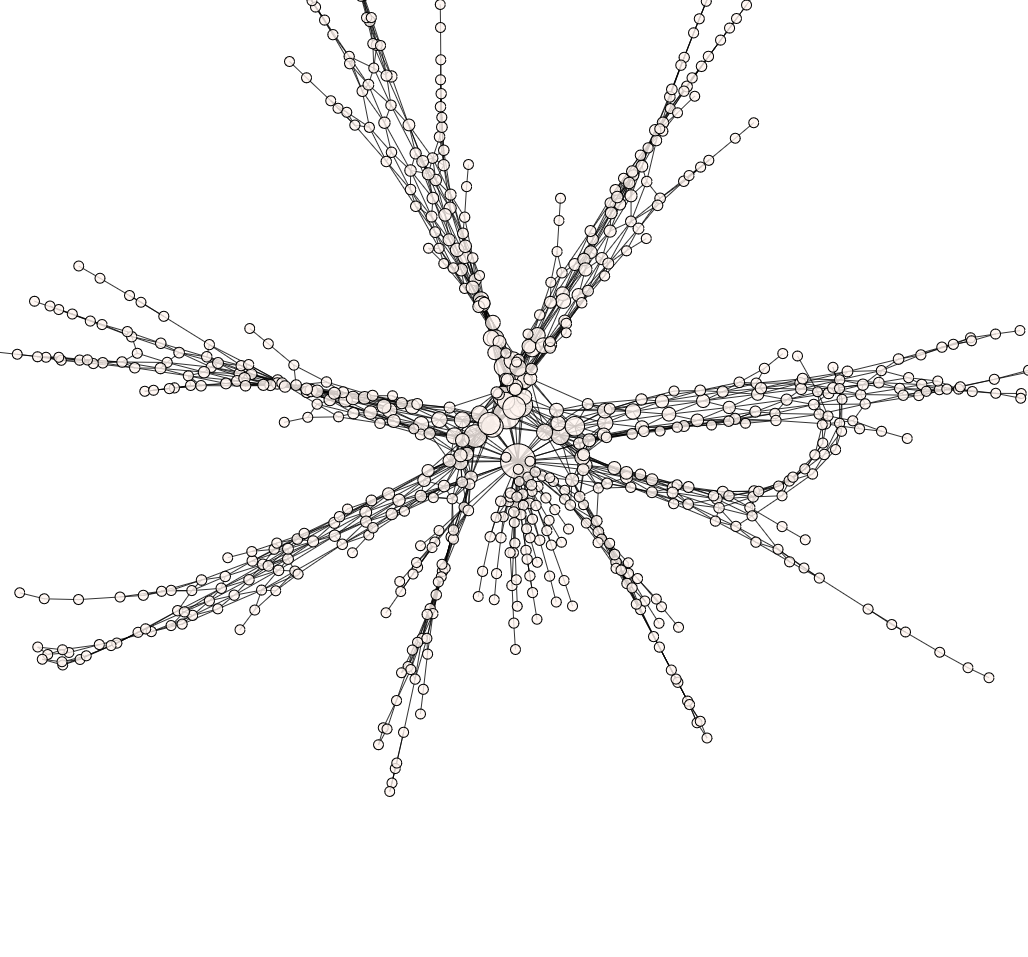}};
\end{tikzpicture}
.
\end{gather*}
These are the \emph{quantum aliens with their tentacles}.

Additionally, and much nicer in the interactive plot, there is a comparison 
of the invariants, {\eg}:
\begin{gather*}
\text{J}\colon
\begin{tikzpicture}[anchorbase]
\node at (0,0) {\includegraphics[height=5.2cm]{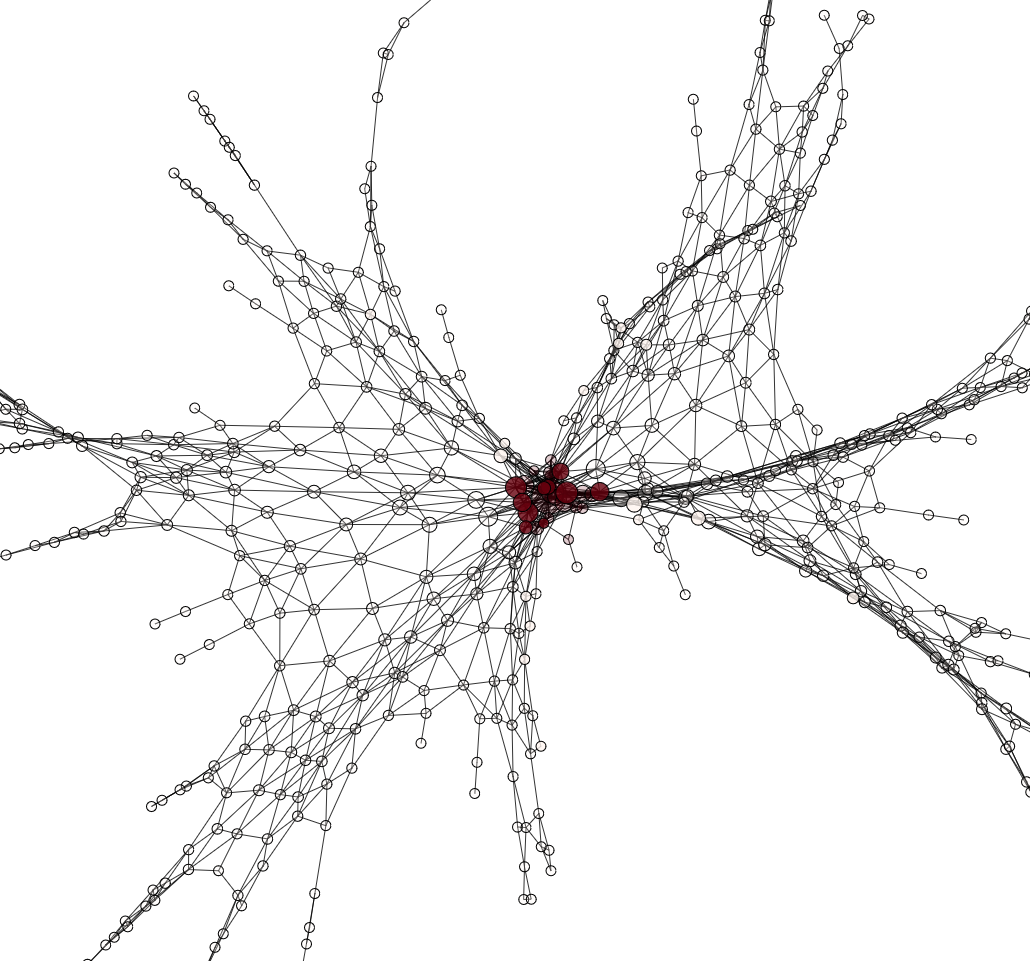}};
\end{tikzpicture}
,
\text{K}\colon
\begin{tikzpicture}[anchorbase]
\node at (0,0) {\includegraphics[height=5.2cm]{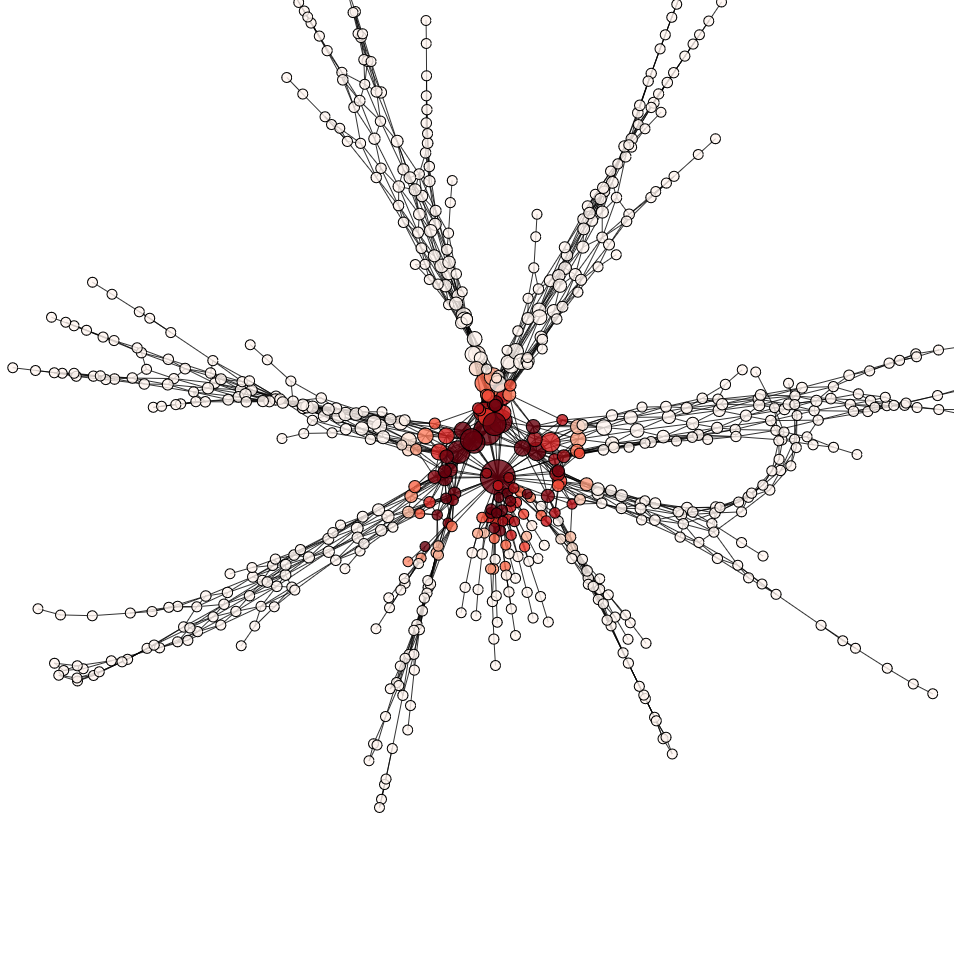}};
\end{tikzpicture}
.
\end{gather*}
This should be read as follows. The selected (and colored) vertices in J appear in the colored range in K.

One observes the following:
\begin{enumerate}

\item All graphs have a dense center. That means that there are a lot of knots with similar quantum invariants.

\item The vertices on the end of the tentacles are small and mostly very special knots such as alternating or torus knots. These are the knots with very unique quantum invariants.

\item The only invariant that is different from the others is B1. All other invariants are roughly similarly spread.
B1 seems to detect something different from the others, while being more closely packed at the center.

\end{enumerate}

\section{Copyable data tables}\label{S:Data}

The title of this section is spot on.

\begin{figure}[ht]
\fcolorbox{tomato!50}{white}{
\rowcolors{2}{spinach!25}{}
\begin{tabular}{c|c|c|c|c|c|c}
\rowcolor{orchid!50}
n & A2 & A & B1 (R) & B1 (Skein) & J & K \\ \hline
3  & 0.000961 & 0.000548 &   0.04195  &   0.006364 & 0.000952 & 0.000819 \\ \hline
4  & 0.001044 & 0.000474 &   0.058764 &   0.029009 & 0.001174 & 0.001009 \\ \hline
5  & 0.001204 & 0.000413 &   0.063757 &   0.027541 & 0.001359 & 0.001187 \\ \hline
6  & 0.001540 & 0.000371 &   0.073052 &   0.106084 & 0.002156 & 0.001333 \\ \hline
7  & 0.002078 & 0.000344 &   0.109800 &   0.152202 & 0.002442 & 0.001458 \\ \hline
8  & 0.006668 & 0.000327 &   0.146327 &   0.382624 & 0.003204 & 0.001624 \\ \hline
9  & 0.010328 & 0.000340 &   0.434684 &   0.700274 & 0.004329 & 0.001772 \\ \hline
10 & 0.022115 & 0.000377 &   0.727308 &   2.001462 & 0.005811 & 0.001989 \\ \hline
11 & 0.035645 & 0.000447 &  41.247445 &   5.641350 & 0.007834 & 0.002090 \\ \hline
12 & 0.073140 & 0.000800 &  13.954301 &  30.957447 & 0.010320 & 0.002089 \\ \hline
13 & 0.189876 & 0.001039 &  28.231526 &  94.590263 & 0.012261 & 0.002250 \\ \hline
14 & 0.393864 & 0.001244 & 113.210082 & 535.154620 & 0.016506 & 0.002312 \\
\end{tabular}
}
\caption{Average time; copyable data.}
\label{figure:0}
\end{figure}

\begin{figure}[ht]
\fcolorbox{tomato!50}{white}{
\rowcolors{2}{spinach!25}{}
\begin{tabular}{c|P{1.2cm}|P{1.2cm}|P{1.2cm}|P{1.2cm}|P{1.2cm}|P{1.2cm}|P{1.2cm}|P{1.2cm}}
\rowcolor{orchid!50}
n & A2 & A & B1 & J & K & KT1 & J+KT1 & All \\ \hline
3 & 100.0 & 100.0 & 100.0 & 100.0 & 100.0 & 100.0 & 100.0 & 100.0 \\ \hline
4 & 100.0 & 100.0 & 100.0 & 100.0 & 100.0 & 100.0 & 100.0 & 100.0 \\ \hline
5 & 100.0 & 100.0 & 100.0 & 100.0 & 100.0 & 100.0 & 100.0 & 100.0 \\ \hline
6 & 100.0 & 100.0 & 100.0 & 100.0 & 100.0 & 100.0 & 100.0 & 100.0 \\ \hline
7 & 100.0 & 100.0 & 100.0 & 100.0 & 100.0 & 100.0 & 100.0 & 100.0 \\ \hline
8 & 100.0 & 100.0 & 100.0 & 100.0 & 100.0 & 100.0 & 100.0 & 100.0 \\ \hline
9 & 100.0 & 94.0 & 100.0 & 100.0 & 100.0 & 100.0 & 100.0 & 100.0 \\ \hline
10 & 98.7 & 84.7 & 100.0 & 96.3 & 96.3 & 96.3 & 96.3 & 100.0 \\ \hline
11 & 95.8 & 68.7 & 98.1 & 90.1 & 91.1 & 90.7 & 91.1 & 98.1 \\ \hline
12 & 92.1 & 59.5 & 96.6 & 83.0 & 84.3 & 83.8 & 84.1 & 96.7 \\ \hline
13 & 85.7 & 43.4 & 92.3 & 73.3 & 77.5 & 77.1 & 77.4 & 92.6 \\ \hline
14 & 81.2 & 33.6 & 89.3 & 64.4 & 69.0 & 68.4 & 68.9 & 89.6 \\ \hline
15 & 76.4 & 24.5 & 86.2 & 55.7 & 60.6 & 59.8 & 60.6 & 86.4 \\ \hline
16 & 74.0 & 18.6 & 83.8 & 49.4 & 54.7 & 53.6 & 54.6 & 84.0 \\
\end{tabular}
}
\caption{Percentages of unique values; copyable data.}
\label{figure:comparison1}
\end{figure}

\newpage

\begin{figure}[ht]
\fcolorbox{tomato!50}{white}{
\rowcolors{2}{spinach!25}{}
\begin{tabular}{c|P{1.5cm}|P{1.5cm}|P{1.5cm}|P{1.5cm}|P{1.5cm}|P{1.5cm}|P{1.5cm}|P{1.5cm}}
\rowcolor{orchid!50}
n & A2 & A & B1 & J & K & KT1 & J+KT1 & All \\ \hline
 4 &      1.0 &     0.9 &       0.9 &      1.0 &      1.0 &      1.0 &      1.0 &       1.0 \\ \hline
 5 &      2.9 &     3.0 &       2.9 &      2.9 &      2.9 &      2.9 &      3.0 &       2.9 \\ \hline
 6 &      6.0 &     5.9 &       6.0 &      6.0 &      6.0 &      5.9 &      5.9 &       5.9 \\ \hline
 7 &     13.0 &    12.9 &      12.9 &     13.0 &     12.9 &     12.9 &     12.9 &      12.9 \\ \hline
 8 &     34.1 &    33.9 &      34.0 &     34.0 &     33.9 &     34.1 &     34.0 &      34.0 \\ \hline
 9 &     83.2 &    74.0 &      82.9 &     82.8 &     83.2 &     82.3 &     83.1 &      82.7 \\ \hline
10 &    241.9 &   187.0 &     247.4 &    231.0 &    232.6 &    231.7 &    230.7 &     250.4 \\ \hline
11 &    738.4 &   417.1 &     764.5 &    662.1 &    673.2 &    671.5 &    676.0 &     768.4 \\ \hline
12 &   2532.7 &  1147.5 &    2791.9 &   2098.2 &   2141.2 &   2117.2 &   2134.0 &    2794.7 \\ \hline
13 &   9579.0 &  2638.0 &   11090.3 &   7145.8 &   7902.2 &   7846.9 &   7916.5 &   11077.2 \\ \hline
14 &  39876.6 &  6679.7 &   47869.6 &  25740.7 &  28845.6 &  28586.1 &  29036.8 &   48079.8 \\ \hline
15 & 181201.6 & 15376.9 &  230546.2 &  98445.1 & 115234.8 & 112986.8 & 114700.2 &  231711.7 \\ \hline
16 & 909889.2 & 34001.2 & 1161532.7 & 394876.7 & 475079.0 & 460340.4 & 473986.2 & 1157489.3 \\
\end{tabular}
}
\caption{Average comparisons until equal; copyable data.}
\label{figure:comparison2}
\end{figure}

\begin{figure}[ht]
\fcolorbox{tomato!50}{white}{
\rowcolors{2}{spinach!25}{}
\begin{tabular}{c|c|c|c|c|c|c}
\rowcolor{orchid!50}
n & A2 & Alexander & B1 & Jones & Khovanov & KhovanovT1 \\ \hline
3  & 2  & 1   & 1     & 1   & 1   & 1   \\ \hline
4  & 2  & 3   & 1     & 1   & 1   & 2   \\ \hline
5  & 2  & 3   & 2     & 2   & 1   & 2   \\ \hline
6  & 2  & 5   & 2     & 3   & 2   & 3   \\ \hline
7  & 2  & 9   & 5     & 4   & 3   & 5   \\ \hline
8  & 4  & 13  & 17    & 9   & 5   & 8   \\ \hline
9  & 4  & 23  & 40    & 13  & 7   & 12  \\ \hline
10 & 6  & 37  & 76    & 21  & 11  & 20  \\ \hline
11 & 9  & 59  & 202   & 34  & 18  & 34  \\ \hline
12 & 11 & 109 & 517   & 61  & 31  & 58  \\ \hline
13 & 22 & 163 & 1440  & 103 & 53  & 100 \\ \hline
14 & 31 & 281 & 3474  & 177 & 89  & 167 \\ \hline
15 & 56 & 451 & 10063 & 305 & 157 & 296 \\ \hline
16 & 87 & 813 & 26219 & 533 & 267 & 503 \\
\end{tabular}
}
\caption{Maximal coefficient; copyable data.}
\label{figure:3}
\end{figure}

\begin{figure}[ht]
\fcolorbox{tomato!50}{white}{
\rowcolors{2}{spinach!25}{}
\begin{tabular}{c|c|c|c|c|c|c}
\rowcolor{orchid!50}
n & A2 & Alexander & B1 & Jones & Khovanov & KhovanovT1 \\ \hline
3 & 7 & 3 & 9 & 3 & 4 & 4 \\ \hline
4 & 7 & 5 & 9 & 5 & 6 & 6 \\ \hline
5 & 9 & 7 & 9 & 7 & 8 & 8 \\ \hline
6 & 9 & 13 & 27 & 13 & 14 & 14 \\ \hline
7 & 11 & 21 & 53 & 21 & 22 & 22 \\ \hline
8 & 19 & 45 & 203 & 45 & 46 & 46 \\ \hline
9 & 27 & 75 & 475 & 75 & 76 & 76 \\ \hline
10 & 37 & 121 & 987 & 121 & 122 & 122 \\ \hline
11 & 57 & 209 & 2609 & 209 & 210 & 210 \\ \hline
12 & 101 & 377 & 7287 & 377 & 378 & 378 \\ \hline
13 & 163 & 663 & 20047 & 663 & 664 & 664 \\ \hline
14 & 257 & 1145 & 51311 & 1145 & 1146 & 1146 \\ \hline
15 & 439 & 2037 & 146629 & 2037 & 2038 & 2038 \\ \hline
16 & 717 & 3581 & 397707 & 3581 & 3582 & 3582 \\
\end{tabular}
}
\caption{Maximal coefficient sum; copyable data.}
\label{figure:4}
\end{figure}

\newpage

\begin{figure}[ht]
\fcolorbox{tomato!50}{white}{
\rowcolors{2}{spinach!25}{}
\begin{tabular}{c|c|c|c|c|c|c}
\rowcolor{orchid!50}
n & A2 & Alexander & B1 & Jones & Khovanov & Khovanovt1 \\ \hline
3  &   7.0000 &   3.0000 &    9.0000 &   3.0000 &   4.0000 &   4.0000 \\ \hline
4  &   6.0000 &   4.0000 &    8.0000 &   4.0000 &   5.0000 &   5.0000 \\ \hline
5  &   7.0000 &   5.0000 &    8.5000 &   5.0000 &   6.0000 &   6.0000 \\ \hline
6  &   7.0000 &   7.5714 &   12.7142 &   7.5714 &   8.5714 &   8.5714 \\ \hline
7  &   8.2857 &  11.1428 &   20.4285 &  11.1428 &  12.1428 &  12.1428 \\ \hline
8  &  10.5428 &  18.8857 &   46.3714 &  18.8285 &  19.9428 &  19.9428 \\ \hline
9  &  13.0238 &  29.6666 &   95.0476 &  29.6428 &  30.7142 &  30.7142 \\ \hline
10 &  17.2409 &  47.6345 &  202.6867 &  47.5461 &  48.8192 &  48.8192 \\ \hline
11 &  23.9762 &  75.9538 &  447.5218 &  75.9837 &  77.5205 &  77.5205 \\ \hline
12 &  33.3903 & 119.2499 &  955.8572 & 119.3829 & 121.4605 & 121.4605 \\ \hline
13 &  47.1625 & 183.4964 & 2034.3612 & 184.1180 & 187.5322 & 187.5322 \\ \hline
14 &  66.6552 & 278.7659 & 4256.5323 & 280.3088 & 286.0564 & 286.0564 \\ \hline
15 &  92.9985 & 412.9504 & 8620.9249 & 416.3381 & 426.4781 & 426.4781 \\ \hline
16 & 130.3453 & 607.3186 & 17436.7151 & 614.2326 & 632.2439 & 632.2439 \\
\end{tabular}
}
\caption{Average coefficient sum; copyable data.}
\label{figure:5}
\end{figure}

\begin{figure}[ht]
\fcolorbox{tomato!50}{white}{
\rowcolors{2}{spinach!25}{}
\begin{tabular}{c|c|c|c|c|c}
\rowcolor{orchid!50}
n & A2 & Alexander & B1 & Jones & Khovanovt1 \\ \hline
3  &  0.5384 &   1.0000 &   0.3913 &   0.7500 &   0.4444 \\ \hline
4  &  0.5384 &   1.6666 &   0.3913 &   1.0000 &   0.5454 \\ \hline
5  &  0.5384 &   2.3333 &   0.3913 &   1.1666 &   0.6153 \\ \hline
6  &  0.5384 &   3.0000 &   0.6585 &   1.8571 &   0.9333 \\ \hline
7  &  0.5384 &   5.0000 &   1.1276 &   2.6250 &   1.2941 \\ \hline
8  &  0.7600 &   6.6000 &   3.8301 &   5.0000 &   2.4210 \\ \hline
9  &  0.9310 &  11.4000 &   8.0508 &   7.5000 &   3.6190 \\ \hline
10 &  1.2758 &  21.0000 &  15.1846 &  11.0000 &   5.3043 \\ \hline
11 &  1.7272 &  29.0000 &  36.7464 &  17.4166 &   8.4000 \\ \hline
12 &  2.7297 &  43.8571 &  94.6363 &  29.0000 &  14.0000 \\ \hline
13 &  4.4054 &  71.8571 & 241.5301 &  47.3571 &  22.8965 \\ \hline
14 &  6.2682 & 124.4285 & 576.5280 &  76.3333 &  36.9677 \\ \hline
15 & 10.7073 & 213.0000 & 1543.4631 & 127.3125 &  61.7575 \\ \hline
16 & 15.9333 & 304.1111 & 3937.6930 & 210.6470 & 102.3428 \\
\end{tabular}
}
\caption{Maximum average coefficient; copyable data.}
\label{figure:6}
\end{figure}

\begin{figure}[ht]
\fcolorbox{tomato!50}{white}{
\rowcolors{2}{spinach!25}{}
\begin{tabular}{c|c|c|c|c|c}
\rowcolor{orchid!50}
n & A2 & Alexander & B1 & Jones & Khovanovt1 \\ \hline
3  & 13 & 3  & 23  & 4  & 9  \\ \hline
4  & 17 & 3  & 29  & 5  & 11 \\ \hline
5  & 19 & 5  & 35  & 6  & 13 \\ \hline
6  & 23 & 5  & 41  & 7  & 15 \\ \hline
7  & 25 & 7  & 47  & 8  & 17 \\ \hline
8  & 29 & 7  & 53  & 9  & 19 \\ \hline
9  & 31 & 9  & 59  & 10 & 21 \\ \hline
10 & 35 & 9  & 65  & 11 & 23 \\ \hline
11 & 37 & 11 & 71  & 12 & 25 \\ \hline
12 & 41 & 11 & 77  & 13 & 27 \\ \hline
13 & 43 & 13 & 83  & 14 & 29 \\ \hline
14 & 47 & 13 & 89  & 15 & 31 \\ \hline
15 & 49 & 15 & 95  & 16 & 33 \\ \hline
16 & 53 & 15 & 101 & 17 & 35 \\
\end{tabular}
}
\caption{Maximal span; copyable data.}
\label{figure:7}
\end{figure}

\newpage

\begin{figure}[ht]
\fcolorbox{tomato!50}{white}{
\rowcolors{2}{spinach!25}{}
\begin{tabular}{c|c|c|c|c|c}
\rowcolor{orchid!50}
n & A2 & Alexander & B1 & Jones & Khovanovt1 \\ \hline
3  & 13.0000 & 3.0000 & 23.0000 &  4.0000 &  9.0000 \\ \hline
4  & 15.0000 & 3.0000 & 26.0000 &  4.5000 & 10.0000 \\ \hline
5  & 16.5000 & 3.5000 & 30.5000 &  5.2500 & 11.5000 \\ \hline
6  & 18.7142 & 3.8571 & 35.0000 &  6.0000 & 13.0000 \\ \hline
7  & 21.2857 & 4.2857 & 41.0000 &  7.0000 & 15.0000 \\ \hline
8  & 23.9714 & 5.1142 & 47.1142 &  8.0000 & 17.0000 \\ \hline
9  & 26.6428 & 5.5238 & 52.8809 &  8.9523 & 18.9047 \\ \hline
10 & 29.0321 & 6.3253 & 58.9036 &  9.9437 & 20.8955 \\ \hline
11 & 31.6117 & 6.9001 & 64.5056 & 10.8514 & 22.7153 \\ \hline
12 & 33.9183 & 7.5072 & 70.0755 & 11.7638 & 24.5391 \\ \hline
13 & 36.0324 & 8.1801 & 75.2020 & 12.5885 & 26.1881 \\ \hline
14 & 38.2418 & 8.6834 & 80.2214 & 13.3948 & 27.8008 \\ \hline
15 & 40.1237 & 9.2937 & 85.0202 & 14.1587 & 29.3301 \\ \hline
16 & 42.1702 & 9.7961 & 89.7204 & 14.9030 & 30.8202 \\
\end{tabular}
}
\caption{Average span; copyable data.}
\label{figure:8}
\end{figure}

\begin{figure}[ht]
\fcolorbox{tomato!50}{white}{
\rowcolors{2}{spinach!25}{}
\begin{tabular}{c|c|c|c|c|c}
\rowcolor{orchid!50}
n & A2 & Alexander & B1 & Jones & Khovanovt1 \\ \hline
3  & 1.1278 &  1.0000 & 1.0731 & 1.2106 & 1.1168 \\ \hline
4  & 1.1837 &  2.6180 & 1.1080 & 1.2106 & 1.2406 \\ \hline
5  & 1.1837 &  2.6180 & 1.1388 & 1.2837 & 1.2451 \\ \hline
6  & 1.2196 &  2.6180 & 1.1990 & 1.6355 & 1.3042 \\ \hline
7  & 1.2547 &  3.3165 & 1.2197 & 1.6355 & 1.6993 \\ \hline
8  & 1.4348 &  4.3902 & 1.7655 & 1.8692 & 1.6993 \\ \hline
9  & 1.7765 &  5.1069 & 1.8668 & 3.3168 & 1.7980 \\ \hline
10 & 1.7765 &  5.1903 & 1.8668 & 3.3168 & 1.9942 \\ \hline
11 & 1.7765 &  7.0507 & 1.9518 & 3.3168 & 2.0254 \\ \hline
12 & 2.0162 &  8.7946 & 2.2990 & 4.3519 & 2.0504 \\ \hline
13 & 2.1548 & 10.0233 & 2.4413 & 4.9314 & 2.1074 \\ \hline
14 & 2.6858 & 14.2600 & 2.7764 & 5.6205 & 2.3300 \\ \hline
15 & 2.7322 & 22.4258 & 3.7987 & 6.4831 & 2.4871 \\ \hline
16 & 2.9737 & 30.5071 & 3.7987 & 7.4391 & 2.6661 \\
\end{tabular}
}
\caption{Maximal absolute root; copyable data.}
\label{figure:9}
\end{figure}

\begin{figure}[ht]
\fcolorbox{tomato!50}{white}{
\rowcolors{2}{spinach!25}{}
\begin{tabular}{c|c|c|c|c|c}
\rowcolor{orchid!50}
n & A2 & Alexander & B1 & Jones & Khovanovt1 \\ \hline
3  & 33.3333 &  0.0000 & 9.0909 & 33.3333 &  0.0000 \\ \hline
4  & 28.5714 & 50.0000 & 4.0000 & 14.2857 & 11.1111 \\ \hline
5  & 22.5806 & 20.0000 & 5.0847 & 17.6470 &  4.7619 \\ \hline
6  & 22.5806 & 30.0000 & 4.2016 & 14.2857 &  9.5238 \\ \hline
7  & 21.1267 & 17.3913 & 4.2857 & 14.2857 & 10.2040 \\ \hline
8  & 19.9004 & 25.0000 & 4.3370 & 11.0204 & 11.4285 \\ \hline
9  & 19.2200 & 21.5789 & 5.1858 & 11.3772 & 11.7021 \\ \hline
10 & 19.0257 & 22.0211 & 5.6457 & 12.1688 & 11.7480 \\ \hline
11 & 19.1843 & 20.7786 & 5.8111 & 11.7222 & 11.4177 \\ \hline
12 & 18.8636 & 20.9580 & 5.9901 & 11.7900 & 10.9652 \\ \hline
13 & 18.8596 & 19.7378 & 6.2311 & 12.0096 & 10.9197 \\ \hline
14 & 18.7889 & 20.0468 & 6.3172 & 12.0427 & 10.6347 \\ \hline
15 & 18.5564 & 19.3164 & 6.4275 & 12.2916 & 10.5783 \\ \hline
16 & 18.5997 & 19.4536 & 6.4558 & 12.3533 & 10.4381 \\
\end{tabular}
}
\caption{Percentage of pure roots; copyable data.}
\label{figure:10}
\end{figure}

\newpage

\begin{figure}[ht]
\fcolorbox{tomato!50}{white}{
\rowcolors{2}{spinach!25}{}
\begin{tabular}{c|c|c|c|c|c}
\rowcolor{orchid!50}
n & A2 & Alexander & B1 & Jones & Khovanovt1 \\ \hline
3  & 33.3333 & 100.0000 & 36.3636 &  0.0000 &  0.0000 \\ \hline
4  & 57.1428 &  50.0000 & 56.0000 & 57.1428 & 11.1111 \\ \hline
5  & 58.0645 &  80.0000 & 61.0169 & 23.5294 & 14.2857 \\ \hline
6  & 54.8387 &  50.0000 & 65.5462 & 22.8571 & 16.6666 \\ \hline
7  & 49.2957 &  65.2173 & 63.9285 & 14.2857 & 20.4081 \\ \hline
8  & 47.7611 &  52.7777 & 60.4708 & 20.4081 & 21.0714 \\ \hline
9  & 50.3249 &  53.1578 & 59.2932 & 19.7604 & 22.0744 \\ \hline
10 & 49.4842 &  45.7013 & 56.9288 & 20.9699 & 22.2850 \\ \hline
11 & 49.5595 &  43.9272 & 56.0352 & 22.0884 & 22.7319 \\ \hline
12 & 48.8356 &  41.2244 & 54.2710 & 23.3678 & 24.0139 \\ \hline
13 & 48.4072 &  39.1900 & 53.3642 & 23.9913 & 24.6212 \\ \hline
14 & 48.6044 &  36.9571 & 52.5937 & 24.0480 & 25.0827 \\ \hline
15 & 47.9160 &  35.8780 & 51.7600 & 23.5892 & 25.3112 \\ \hline
16 & 47.5652 &  34.3872 & 51.1914 & 23.4344 & 25.6379 \\
\end{tabular}
}
\caption{Percentage of roots in $[0.9,1.1]$; copyable data.}
\label{figure:11}
\end{figure}

\newcommand{\etalchar}[1]{$^{#1}$}


\begin{thebibliography}{MMM{\etalchar{+}}23}

\bibitem[Ada94]{Ad-knots}
C.C. Adams.
\newblock {\em The knot book}.
\newblock W. H. Freeman and Company, New York, 1994.
\newblock An elementary introduction to the mathematical theory of knots.

\bibitem[AK92]{AnWe-mixed-qgroup}
H.H. Andersen and W.X. Kexin.
\newblock Representations of quantum algebras. {T}he mixed case.
\newblock {\em J. Reine Angew. Math.}, 427:35--50, 1992.
\newblock \href {https://doi.org/10.1515/crll.1992.427.35}
  {\path{doi:10.1515/crll.1992.427.35}}.

\bibitem[APR89]{AnPrRo-mutation}
R.P. Anstee, J.H. Przytycki, and D.~Rolfsen.
\newblock Knot polynomials and generalized mutation.
\newblock {\em Topology Appl.}, 32(3):237--249, 1989.
\newblock URL: \url{https://arxiv.org/abs/math/0405382}, \href
  {https://doi.org/10.1016/0166-8641(89)90031-X}
  {\path{doi:10.1016/0166-8641(89)90031-X}}.

\bibitem[Atl24]{KnotTheory}
K.~Atlas.
\newblock The {M}athematica package {KnotTheory`}, version september 27, 2024.
\newblock 2024.
\newblock URL:
  \url{https://katlas.org/wiki/The_Mathematica_Package_KnotTheory%60}.

\bibitem[BCD23]{BaChDe-roots}
J.C. Baez, J.D. Christensen, and S.~Derbyshire.
\newblock The beauty of roots.
\newblock {\em Notices Amer. Math. Soc.}, 70(9):1495--1497, 2023.
\newblock URL: \url{https://arxiv.org/abs/2310.00326}, \href
  {https://doi.org/10.1090/noti2789} {\path{doi:10.1090/noti2789}}.

\bibitem[Bai97]{Ba-circular-law}
Z.D. Bai.
\newblock Circular law.
\newblock {\em Ann. Probab.}, 25(1):494--529, 1997.
\newblock \href {https://doi.org/10.1214/aop/1024404298}
  {\path{doi:10.1214/aop/1024404298}}.

\bibitem[BN02]{BaNa-categorification-jones}
D.~Bar-Natan.
\newblock On {K}hovanov's categorification of the {J}ones polynomial.
\newblock {\em Algebr. Geom. Topol.}, 2:337--370, 2002.
\newblock URL: \url{https://arxiv.org/abs/math/0201043}, \href
  {https://doi.org/10.2140/agt.2002.2.337} {\path{doi:10.2140/agt.2002.2.337}}.

\bibitem[BN07]{BaNa-fast}
D.~Bar-Natan.
\newblock Fast {K}hovanov homology computations.
\newblock {\em J. Knot Theory Ramifications}, 16(3):243--255, 2007.
\newblock URL: \url{https://arxiv.org/abs/math/0606318}, \href
  {https://doi.org/10.1142/S0218216507005294}
  {\path{doi:10.1142/S0218216507005294}}.

\bibitem[BNvdV19]{BaNaVe-polynomial-time-invariant}
D.~Bar-Natan and R.~van~der Veen.
\newblock A polynomial time knot polynomial.
\newblock {\em Proc. Amer. Math. Soc.}, 147(1):377--397, 2019.
\newblock URL: \url{https://arxiv.org/abs/1708.04853}, \href
  {https://doi.org/10.1090/proc/14166} {\path{doi:10.1090/proc/14166}}.

\bibitem[CKW17]{CaKoWa}
D.~Calegari, S.~Koch, and A.~Walker.
\newblock Roots, {S}chottky semigroups, and a proof of {B}andt's conjecture.
\newblock {\em Ergodic Theory Dynam. Systems}, 37(8):2487--2555, 2017.
\newblock URL: \url{https://arxiv.org/abs/1410.8542}, \href
  {https://doi.org/10.1017/etds.2016.17} {\path{doi:10.1017/etds.2016.17}}.

\bibitem[CGKL21]{ChGoKoLa-skein}
M.~Chlouveraki, D.~Goundaroulis, A.~Kontogeorgis, and S.~Lambropoulou.
\newblock A generalized skein relation for {K}hovanov homology and a
  categorification of the {$\theta$}-invariant.
\newblock {\em Proc. Roy. Soc. Edinburgh Sect. A}, 151(6):1731--1757, 2021.
\newblock URL: \url{https://arxiv.org/abs/1904.07794}, \href
  {https://doi.org/10.1017/prm.2020.78} {\path{doi:10.1017/prm.2020.78}}.

\bibitem[COT24]{CoOsTu-growth}
K.~Coulembier, V.~Ostrik, and D.~Tubbenhauer.
\newblock Growth rates of the number of indecomposable summands in tensor
  powers.
\newblock {\em Algebr. Represent. Theory}, 27(2):1033--1062, 2024.
\newblock URL: \url{https://arxiv.org/abs/2301.00885}, \href
  {https://doi.org/10.1007/s10468-023-10245-7}
  {\path{doi:10.1007/s10468-023-10245-7}}.

\bibitem[D{\l}o19]{Dl-ballmapper}
P.~D{\l}otko.
\newblock Ball mapper: a shape summary for topological data analysis.
\newblock 2019.
\newblock URL: \url{https://arxiv.org/abs/1901.07410}.

\bibitem[DGS25]{XYZ}
P.~D{\l}otko, D.~Gurnari, and R.~Sazdanovic.
\newblock Data-driven perspectives on knot invariants.
\newblock 2025.
\newblock URL: Not yet available, but probably available when you read this.

\bibitem[DGS24]{DlGuSa-mapper}
P.~D{\l}otko, D.~Gurnari, and R.~Sazdanovic.
\newblock Mapper--{T}ype {A}lgorithms for {C}omplex {D}ata and {R}elations.
\newblock {\em J. Comput. Graph. Statist.}, 33(4):1383--1396, 2024.
\newblock URL: \url{https://arxiv.org/abs/2109.00831}, \href
  {https://doi.org/10.1080/10618600.2024.2343321}
  {\path{doi:10.1080/10618600.2024.2343321}}.

\bibitem[Eli17]{El-q-satake}
B.~Elias.
\newblock Quantum {S}atake in type {$A$}. {P}art {I}.
\newblock {\em J. Comb. Algebra}, 1(1):63--125, 2017.
\newblock URL: \url{https://arxiv.org/abs/1403.5570}, \href
  {https://doi.org/10.4171/JCA/1-1-4} {\path{doi:10.4171/JCA/1-1-4}}.

\bibitem[ES87]{ErSu-growth-knots}
C.~Ernst and D.W. Sumners.
\newblock The growth of the number of prime knots.
\newblock {\em Math. Proc. Cambridge Philos. Soc.}, 102(2):303--315, 1987.
\newblock \href {https://doi.org/10.1017/S0305004100067323}
  {\path{doi:10.1017/S0305004100067323}}.

\bibitem[EGNO15]{EtGeNiOs-tensor-categories}
P.~Etingof, S.~Gelaki, D.~Nikshych, and V.~Ostrik.
\newblock {\em Tensor categories}, volume 205 of {\em Mathematical Surveys and
  Monographs}.
\newblock American Mathematical Society, Providence, RI, 2015.
\newblock \href {https://doi.org/10.1090/surv/205}
  {\path{doi:10.1090/surv/205}}.

\bibitem[EZ17]{EvZo-random-knots}
C.~Even-Zohar.
\newblock Models of random knots.
\newblock {\em J. Appl. Comput. Topol.}, 1(2):263--296, 2017.
\newblock URL: \url{https://arxiv.org/pdf/1711.10470.pdf}, \href
  {https://doi.org/10.1007/s41468-017-0007-8}
  {\path{doi:10.1007/s41468-017-0007-8}}.

\bibitem[EZHLN16]{EvZoHaLiNo-random-knots}
C.~Even-Zohar, J.~Hass, N.~Linial, and T.~Nowik.
\newblock Invariants of random knots and links.
\newblock {\em Discrete Comput. Geom.}, 56(2):274--314, 2016.
\newblock URL: \url{https://arxiv.org/abs/1411.3308}, \href
  {https://doi.org/10.1007/s00454-016-9798-y}
  {\path{doi:10.1007/s00454-016-9798-y}}.

\bibitem[FL21]{FlLa-knots}
J.~Flake and R.~Laugwitz.
\newblock On the monoidal center of {D}eligne's category {$\underline {\rm
  Re}{\rm p}(S_t)$}.
\newblock {\em J. Lond. Math. Soc. (2)}, 103(3):1153--1185, 2021.
\newblock URL: \url{https://arxiv.org/abs/1901.08657}, \href
  {https://doi.org/10.1112/jlms.12403} {\path{doi:10.1112/jlms.12403}}.

\bibitem[GKPM11]{GeKuPaMi-gen-traces-modified-dimensions}
N.~Geer, J.~Kujawa, and B.~Patureau-Mirand.
\newblock Generalized trace and modified dimension functions on ribbon
  categories.
\newblock {\em Selecta Math. (N.S.)}, 17(2):453--504, 2011.
\newblock URL: \url{https://arxiv.org/pdf/1001.0985.pdf}, \href
  {https://doi.org/10.1007/s00029-010-0046-7}
  {\path{doi:10.1007/s00029-010-0046-7}}.

\bibitem[GPMT09]{GePaMiTu-modified-qdims}
N.~Geer, B.~Patureau-Mirand, and V.~Turaev.
\newblock Modified quantum dimensions and re-normalized link invariants.
\newblock {\em Compos. Math.}, 145(1):196--212, 2009.
\newblock URL: \url{https://arxiv.org/abs/0711.4229}, \href
  {https://doi.org/10.1112/S0010437X08003795}
  {\path{doi:10.1112/S0010437X08003795}}.

\bibitem[Gem86]{Ge-spectral}
S.~Geman.
\newblock The spectral radius of large random matrices.
\newblock {\em Ann. Probab.}, 14(4):1318--1328, 1986.
\newblock \href {https://doi.org/10.1214/aop/1176992372}
  {\path{doi:10.1214/aop/1176992372}}.

\bibitem[Guk05]{Gu-cs-and-the-a-polynomial}
S.~Gukov.
\newblock Three-dimensional quantum gravity, {C}hern--{S}imons theory, and the
  {A}-polynomial.
\newblock {\em Comm. Math. Phys.}, 255(3):577--627, 2005.
\newblock URL: \url{https://arxiv.org/abs/hep-th/0306165}, \href
  {https://doi.org/10.1007/s00220-005-1312-y}
  {\path{doi:10.1007/s00220-005-1312-y}}.

\bibitem[Jon85]{Jo-jones-polynomial}
V.F.R. Jones.
\newblock A polynomial invariant for knots via von {N}eumann algebras.
\newblock {\em Bull. Amer. Math. Soc. (N.S.)}, 12(1):103--111, 1985.
\newblock \href {https://doi.org/10.1090/S0273-0979-1985-15304-2}
  {\path{doi:10.1090/S0273-0979-1985-15304-2}}.

\bibitem[Kac49]{Ka-roots-random}
M.~Kac.
\newblock On the average number of real roots of a random algebraic equation.
  {II}.
\newblock {\em Proc. London Math. Soc. (2)}, 50:390--408, 1949.
\newblock \href {https://doi.org/10.1112/plms/s2-50.5.390}
  {\path{doi:10.1112/plms/s2-50.5.390}}.

\bibitem[Kac77]{Ka-lie-super}
V.G. Kac.
\newblock Lie superalgebras.
\newblock {\em Advances in Math.}, 26(1):8--96, 1977.
\newblock \href {https://doi.org/10.1016/0001-8708(77)90017-2}
  {\path{doi:10.1016/0001-8708(77)90017-2}}.

\bibitem[KTB15]{MR3431029}
A.~Kawauchi, I.~Tayama, and B.~Burton.
\newblock Tabulation of 3-manifolds of lengths up to 10.
\newblock {\em Topology Appl.}, 196:937--975, 2015.
\newblock \href {https://doi.org/10.1016/j.topol.2015.05.036}
  {\path{doi:10.1016/j.topol.2015.05.036}}.

\bibitem[Kho00]{Kh-cat-jones}
M.~Khovanov.
\newblock A categorification of the {J}ones polynomial.
\newblock {\em Duke Math. J.}, 101(3):359--426, 2000.
\newblock URL: \url{http://arxiv.org/abs/math/9908171}, \href
  {https://doi.org/10.1215/S0012-7094-00-10131-7}
  {\path{doi:10.1215/S0012-7094-00-10131-7}}.

\bibitem[KST24]{KhSiTu-monoidal-cryptography}
M.~Khovanov, M.~Sitaraman, and D.~Tubbenhauer.
\newblock Monoidal categories, representation gap and cryptography.
\newblock {\em Trans. Amer. Math. Soc. Ser. B}, 11:329--395, 2024.
\newblock URL: \url{https://arxiv.org/abs/2201.01805}, \href
  {https://doi.org/10.1090/btran/151} {\path{doi:10.1090/btran/151}}.

\bibitem[KM11]{KrMo-unknot-detector}
P.B. Kronheimer and T.S. Mrowka.
\newblock Khovanov homology is an unknot-detector.
\newblock {\em Publ. Math. Inst. Hautes {\'E}tudes Sci.}, (113):97--208, 2011.
\newblock URL: \url{https://arxiv.org/abs/1005.4346}, \href
  {https://doi.org/10.1007/s10240-010-0030-y}
  {\path{doi:10.1007/s10240-010-0030-y}}.

\bibitem[Kup96]{Ku-spiders-rank-2}
G.~Kuperberg.
\newblock Spiders for rank {$2$} {L}ie algebras.
\newblock {\em Comm. Math. Phys.}, 180(1):109--151, 1996.
\newblock URL: \url{https://arxiv.org/abs/q-alg/9712003}.

\bibitem[LTV24a]{LaTuVa-big-data}
A.~Lacabanne, D.~Tubbenhauer, and P.~Vaz.
\newblock Big data approach to {K}azhdan--{L}usztig polynomials.
\newblock 2024.
\newblock URL: \url{https://arxiv.org/abs/2412.01283}.

\bibitem[LTV24b]{LaTuVa-nhedral}
A.~Lacabanne, D.~Tubbenhauer, and P.~Vaz.
\newblock On {H}ecke and asymptotic categories for complex reflection groups.
\newblock 2024.
\newblock URL: \url{https://arxiv.org/abs/2409.01005}.

\bibitem[LTV22]{LaTuVa-verma-howe}
A.~Lacabanne, D.~Tubbenhauer, and P.~Vaz.
\newblock Verma {H}owe duality and {LKB} representations.
\newblock 2022.
\newblock URL: \url{https://arxiv.org/abs/2207.09124}.

\bibitem[LZ06]{LeZh-strongly-multiplicity-free}
G.I. Lehrer and R.B. Zhang.
\newblock Strongly multiplicity free modules for {L}ie algebras and quantum
  groups.
\newblock {\em J. Algebra}, 306(1):138--174, 2006.
\newblock \href {https://doi.org/10.1016/j.jalgebra.2006.03.043}
  {\path{doi:10.1016/j.jalgebra.2006.03.043}}.

\bibitem[LHS22]{LeHaSa-jones}
J.S.F. Levitt, M.~Hajij, and R.~Sazdanovic.
\newblock Big data approaches to knot theory: understanding the structure of
  the {J}ones polynomial.
\newblock {\em J. Knot Theory Ramifications}, 31(13):Paper No. 2250095, 20,
  2022.
\newblock URL: \url{https://arxiv.org/abs/1912.10086}, \href
  {https://doi.org/10.1142/s021821652250095x}
  {\path{doi:10.1142/s021821652250095x}}.

\bibitem[LM25]{knotinfo}
C.~Livingston and A.H. Moore.
\newblock Knotinfo: Table of knot invariants.
\newblock URL: \url{knotinfo.math.indiana.edu}, 2025.

\bibitem[MMMT20]{MaMaMiTu-trihedral}
M.~Mackaay, V.~Mazorchuk, V.~Miemietz, and D.~Tubbenhauer.
\newblock Trihedral {S}oergel bimodules.
\newblock {\em Fund. Math.}, 248(3):219--300, 2020.
\newblock URL: \url{https://arxiv.org/abs/1804.08920}, \href
  {https://doi.org/10.4064/fm566-3-2019} {\path{doi:10.4064/fm566-3-2019}}.

\bibitem[MMM{\etalchar{+}}23]{MaMaMiTuZh-soergel-2reps}
M.~Mackaay, V.~Mazorchuk, V.~Miemietz, D.~Tubbenhauer, and X.~Zhang.
\newblock Simple transitive {$2$}-representations of {S}oergel bimodules for
  finite {C}oxeter types.
\newblock {\em Proc. Lond. Math. Soc. (3)}, 126(5):1585--1655, 2023.
\newblock URL: \url{https://arxiv.org/abs/1906.11468}, \href
  {https://doi.org/10.1112/plms.12515} {\path{doi:10.1112/plms.12515}}.

\bibitem[MT19]{MaTu-soergel}
M.~Mackaay and D.~Tubbenhauer.
\newblock Two-color {S}oergel calculus and simple transitive
  {$2$}-representations.
\newblock {\em Canad. J. Math.}, 71(6):1523--1566, 2019.
\newblock URL: \url{https://arxiv.org/abs/1609.00962}, \href
  {https://doi.org/10.4153/CJM-2017-061-2} {\path{doi:10.4153/CJM-2017-061-2}}.

\bibitem[Mar21]{Ma-complexity-quantum}
C.~Maria.
\newblock Parameterized complexity of quantum knot invariants.
\newblock In {\em 37th {I}nternational {S}ymposium on {C}omputational
  {G}eometry}, volume 189 of {\em LIPIcs. Leibniz Int. Proc. Inform.}, pages
  Art. No. 53, 17. Schloss Dagstuhl. Leibniz-Zent. Inform., Wadern, 2021.
\newblock URL: \url{https://arxiv.org/abs/1910.00477}.

\bibitem[MBF{\etalchar{+}}97]{MeBeFoMaAa}
G.~Mezincescu, D.~Bessis, J.-D. Fournier, G.~Mantica, and F.D. Aaron.
\newblock Distribution of roots of random real generalized polynomials.
\newblock {\em J. Statist. Phys.}, 86(3-4):675--705, 1997.
\newblock URL: \url{https://arxiv.org/pdf/chao-dyn/9606012}, \href
  {https://doi.org/10.1007/BF02199115} {\path{doi:10.1007/BF02199115}}.

\bibitem[MPS{\etalchar{+}}18]{MuPrSiWaYa-torsion}
S.~Mukherjee, J.H. Przytycki, M.~Silvero, X.~Wang, and S.Y. Yang.
\newblock Search for torsion in {K}hovanov homology.
\newblock {\em Exp. Math.}, 27(4):488--497, 2018.
\newblock URL: \url{https://arxiv.org/abs/1701.04924}, \href
  {https://doi.org/10.1080/10586458.2017.1320242}
  {\path{doi:10.1080/10586458.2017.1320242}}.

\bibitem[{OEI}23]{Oeis}
{OEIS Foundation Inc.}
\newblock The {O}n-{L}ine {E}ncyclopedia of {I}nteger {S}equences, 2023.
\newblock Published electronically at \url{http://oeis.org}.

\bibitem[ORS13]{OzRaSz-odd}
P.S. Ozsv\'ath, J.~Rasmussen, and Z.~Szab\'o.
\newblock Odd {K}hovanov homology.
\newblock {\em Algebr. Geom. Topol.}, 13(3):1465--1488, 2013.
\newblock URL: \url{https://arxiv.org/abs/0710.4300}, \href
  {https://doi.org/10.2140/agt.2013.13.1465}
  {\path{doi:10.2140/agt.2013.13.1465}}.

\bibitem[OS04]{HFK1}
P.~Ozsv\'ath and Z.~Szab\'o.
\newblock Holomorphic disks and knot invariants.
\newblock {\em Adv. Math.}, 186(1):58--116, 2004.
\newblock URL: \url{hhttps://arxiv.org/abs/math/0209056}, \href
  {https://doi.org/10.1016/j.aim.2003.05.001}
  {\path{doi:10.1016/j.aim.2003.05.001}}.

\bibitem[PS24]{PrSi-complexity}
J.H. Przytycki and M.~Silvero.
\newblock Khovanov homology, wedges of spheres and complexity.
\newblock {\em Rev. R. Acad. Cienc. Exactas F\'is. Nat. Ser. A Mat. RACSAM},
  118(3):Paper No. 102, 34, 2024.
\newblock URL: \url{https://arxiv.org/abs/2305.18648}, \href
  {https://doi.org/10.1007/s13398-024-01594-z}
  {\path{doi:10.1007/s13398-024-01594-z}}.

\bibitem[Ras03]{HFK2}
J.~Rasmussen.
\newblock Floer homology and knot complements.
\newblock 2003.
\newblock Ph.D. thesis, Harvard University.
\newblock URL: \url{https://arxiv.org/abs/math/0306378}.

\bibitem[RT91]{ReTu-invariants-3-manifolds-qgroups}
N.~Reshetikhin and V.G. Turaev.
\newblock Invariants of {$3$}-manifolds via link polynomials and quantum
  groups.
\newblock {\em Invent. Math.}, 103(3):547--597, 1991.
\newblock \href {https://doi.org/10.1007/BF01239527}
  {\path{doi:10.1007/BF01239527}}.

\bibitem[RT16]{RoTu-symmetric-howe}
D.E.V. Rose and D.~Tubbenhauer.
\newblock Symmetric webs, {J}ones--{W}enzl recursions, and {$q$}-{H}owe
  duality.
\newblock {\em Int. Math. Res. Not. IMRN}, (17):5249--5290, 2016.
\newblock URL: \url{https://arxiv.org/abs/1501.00915}, \href
  {https://doi.org/10.1093/imrn/rnv302} {\path{doi:10.1093/imrn/rnv302}}.

\bibitem[Sar15]{Sa-alexander}
A.~Sartori.
\newblock The {A}lexander polynomial as quantum invariant of links.
\newblock {\em Ark. Mat.}, 53(1):177--202, 2015.
\newblock URL: \url{https://arxiv.org/abs/1308.2047}, \href
  {https://doi.org/10.1007/s11512-014-0196-5}
  {\path{doi:10.1007/s11512-014-0196-5}}.

\bibitem[SW24]{SaWe-quantum-f4}
A.~Savage and B.W. Westbury.
\newblock Quantum diagrammatics for {$F_4$}.
\newblock {\em J. Pure Appl. Algebra}, 228(11):Paper No. 107731, 35, 2024.
\newblock URL: \url{https://arxiv.org/abs/2204.11976}, \href
  {https://doi.org/10.1016/j.jpaa.2024.107731}
  {\path{doi:10.1016/j.jpaa.2024.107731}}.

\bibitem[ST94]{SeTh-ratcatcher}
P.D. Seymour and R.~Thomas.
\newblock Call routing and the ratcatcher.
\newblock {\em Combinatorica}, 14(2):217--241, 1994.
\newblock \href {https://doi.org/10.1007/BF01215352}
  {\path{doi:10.1007/BF01215352}}.

\bibitem[SV95]{ShVa-roots-random}
L.A. Shepp and R.J. Vanderbei.
\newblock The complex zeros of random polynomials.
\newblock {\em Trans. Amer. Math. Soc.}, 347(11):4365--4384, 1995.
\newblock \href {https://doi.org/10.2307/2155041} {\path{doi:10.2307/2155041}}.

\bibitem[SMC07]{mapper}
G.~Singh, F.~Memoli, and G.~Carlsson.
\newblock {Topological Methods for the Analysis of High Dimensional Data Sets
  and 3D Object Recognition}.
\newblock In M.~Botsch, R.~Pajarola, B.~Chen, and M.~Zwicker, editors, {\em
  Eurographics Symposium on Point-Based Graphics}. The Eurographics
  Association, 2007.
\newblock \href {https://doi.org/10.2312/SPBG/SPBG07/091-100}
  {\path{doi:10.2312/SPBG/SPBG07/091-100}}.

\bibitem[Sto04]{St-number-polynomials}
A.~Stoimenow.
\newblock On the number of links and link polynomials.
\newblock {\em Q. J. Math.}, 55(1):87--98, 2004.
\newblock \href {https://doi.org/10.1093/qjmath/55.1.87}
  {\path{doi:10.1093/qjmath/55.1.87}}.

\bibitem[ST98]{SuTh-growth}
C.~Sundberg and M.~Thistlethwaite.
\newblock The rate of growth of the number of prime alternating links and
  tangles.
\newblock {\em Pacific J. Math.}, 182(2):329--358, 1998.
\newblock \href {https://doi.org/10.2140/pjm.1998.182.329}
  {\path{doi:10.2140/pjm.1998.182.329}}.

\bibitem[STWZ23]{SuTuWeZh-mixed-tilting}
L.~Sutton, D.~Tubbenhauer, P.~Wedrich, and J.~Zhu.
\newblock {SL2} tilting modules in the mixed case.
\newblock {\em Selecta Math. (N.S.)}, 29(3):39, 2023.
\newblock URL: \url{https://arxiv.org/abs/2105.07724}, \href
  {https://doi.org/10.1007/s00029-023-00835-0}
  {\path{doi:10.1007/s00029-023-00835-0}}.

\bibitem[Tub24]{Tu-web-reps}
D.~Tubbenhauer.
\newblock On rank one 2-representations of web categories.
\newblock {\em Algebr. Comb.}, 7(6):1813--1843, 2024.
\newblock URL: \url{https://arxiv.org/abs/2307.00785}, \href
  {https://doi.org/10.5802/alco.389} {\path{doi:10.5802/alco.389}}.

\bibitem[Tub22]{Tu-qt}
D.~Tubbenhauer.
\newblock Quantum topology without topology.
\newblock 2022.
\newblock URL: \url{https://www.dtubbenhauer.com/qinvariants.pdf}.

\bibitem[TZ25]{TuZh-quantum-big-data-code}
D.~Tubbenhauer and V.~Zhang.
\newblock {C}ode and more for the paper {B}ig data comparison of quantum
  invariants.
\newblock 2025.
\newblock \url{https://github.com/dtubbenhauer/quantumdata}.
\newblock \url{https://dustbringer.github.io/web--knot-invariant-comparison}.

\bibitem[Tur94]{Tu-qgroups-3mfds}
V.G. Turaev.
\newblock {\em Quantum invariants of knots and {$3$}-manifolds}, volume~18 of
  {\em De Gruyter Studies in Mathematics}.
\newblock Walter de Gruyter \& Co., Berlin, 1994.

\bibitem[Web17]{We-knot-invariants}
B.~Webster.
\newblock Knot invariants and higher representation theory.
\newblock {\em Mem. Amer. Math. Soc.}, 250(1191):v+141, 2017.
\newblock URL: \url{https://arxiv.org/abs/1309.3796}, \href
  {https://doi.org/10.1090/memo/1191} {\path{doi:10.1090/memo/1191}}.

\bibitem[Weh03]{We-mutation}
S.M. Wehrli.
\newblock {K}hovanov homology and {C}onway mutation.
\newblock 2003.
\newblock URL: \url{https://arxiv.org/abs/math/0301312}.

\bibitem[Weh10]{We-mutation-2}
S.M. Wehrli.
\newblock Mutation invariance of {K}hovanov homology over {$\mathbb{F}_2$}.
\newblock {\em Quantum Topol.}, 1(2):111--128, 2010.
\newblock URL: \url{https://arxiv.org/abs/0904.3401}, \href
  {https://doi.org/10.4171/QT/3} {\path{doi:10.4171/QT/3}}.

\bibitem[Wel92]{We-number-knots}
D.J.A. Welsh.
\newblock On the number of knots and links.
\newblock In {\em Sets, graphs and numbers ({B}udapest, 1991)}, volume~60 of
  {\em Colloq. Math. Soc. J\'anos Bolyai}, pages 713--718. North-Holland,
  Amsterdam, 1992.

\bibitem[Wel93]{We-complexity}
D.J.A. Welsh.
\newblock {\em Complexity: knots, colourings and counting}, volume 186 of {\em
  London Mathematical Society Lecture Note Series}.
\newblock Cambridge University Press, Cambridge, 1993.
\newblock \href {https://doi.org/10.1017/CBO9780511752506}
  {\path{doi:10.1017/CBO9780511752506}}.

\bibitem[Wes12]{110646}
B.~Westbury.
\newblock Links with same {J}ones polynomial, comment on that page.
\newblock MathOverflow (version: 2012-10-25), 2012.
\newblock URL: \url{https://mathoverflow.net/q/110646}.

\bibitem[WW01]{WuWa-jones-roots}
F.~Y. Wu and J.~Wang.
\newblock Zeroes of the {J}ones polynomial.
\newblock {\em Phys. A}, 296(3-4):483--494, 2001.
\newblock URL: \url{https://arxiv.org/pdf/cond-mat/0105013}, \href
  {https://doi.org/10.1016/S0378-4371(01)00189-3}
  {\path{doi:10.1016/S0378-4371(01)00189-3}}.

\end{thebibliography}
\end{document}